\pdfoutput=1

\documentclass[aos]{imsart}

\RequirePackage{amsthm,amsmath,amsfonts,amssymb}
\RequirePackage[numbers]{natbib}
\RequirePackage[colorlinks,citecolor=blue,urlcolor=blue,hypertexnames=false]{hyperref}
\RequirePackage{graphicx}

\RequirePackage{bbm,stmaryrd,enumitem}

\startlocaldefs
\theoremstyle{plain}

\newtheorem{theorem}{Theorem}[section]
\newtheorem{lemma}[theorem]{Lemma}
\newtheorem{proposition}[theorem]{Proposition}
\theoremstyle{remark}
\newtheorem{definition}[theorem]{Definition}
\newtheorem{assumption}{Assumption}

\newtheorem{remark}[theorem]{Remark}
\newtheorem*{notation}{Notation}
\newcommand{\dimu}{\underline{\Delta}^{i \mu}}
\newcommand{\uSig}{\underline{\Sigma}}
\newcommand{\qand}{\quad \text{ and } \quad}

\newcommand{\oprec}{\mathcal{O}_{\prec}}

\newcommand{\bfk}[1]{\boldsymbol{\mathfrak{#1}}}

\endlocaldefs

\begin{document}

\begin{frontmatter}
\title{Asymptotic distribution of spiked eigenvalues in the large signal-plus-noise model}
\runtitle{Large Signal-plus-noise Model}

\begin{aug}
\author[A]{\fnms{Zeqin}~\snm{Lin}\ead[label=e1]{ZEQIN001@e.ntu.edu.sg}},
\author[A]{\fnms{Guangming}~\snm{Pan}\ead[label=e2]{GMPAN@ntu.edu.sg}\orcid{0000-0000-0000-0000}}
\author[B]{\fnms{Peng}~\snm{Zhao}\ead[label=e3]{zhaop@jsnu.edu.cn}}
\and
\author[C]{\fnms{Jia}~\snm{Zhou}\ead[label=e4]{tszhjia@mail.ustc.edu.cn }}
\address[A]{School of Physical and Mathematical Sciences,
Nanyang Technological University
\printead[presep={,\ }]{e1,e2}}

\address[B]{School of Mathematics and Statistics, 
Jiangsu Normal University
\printead[presep={,\ }]{e3}}

\address[C]{School of Economics, 
Hefei University of Technology 
\printead[presep={,\ }]{e4}}
\end{aug}

\begin{abstract}
Consider large signal-plus-noise data matrices of the form $S + \Sigma^{1/2} X$, where $S$ is a low-rank deterministic signal matrix and the noise covariance matrix $\Sigma$ can be anisotropic. We establish the asymptotic joint distribution of its spiked singular values when the dimensionality and sample size are comparably large and the signals are supercritical under general assumptions concerning the structure of $(S, \Sigma)$ and the distribution of the random noise $X$. It turns out that the asymptotic distributions exhibit nonuniversality in the sense of dependence on the distributions of the entries of $X$, which contrasts with what has previously been established for the spiked sample eigenvalues in the context of spiked population models. Such a result yields the asymptotic distribution of the sample spiked eigenvalues associated with mixture models. We also explore the application of these findings in detecting mean heterogeneity of data matrices.
\end{abstract}




\end{frontmatter}

\section{Introduction}
\label{sec-introduction}

In many scientific endeavors, the signal that researchers intend to measure or analyze is often contaminated by random noise originating from a variety of sources, including measurement error, environmental factors, or intrinsic variability in the system being studied. To analyze such signal-plus-noise data, a common mathematical model is given by a $M \times N$ measurement matrix formed by concatenating the $N$ observation vectors of dimension $M$ alongside each other. Specifically, consider the data matrix 
\begin{equation}
    \tilde{Y} = S + Y,
    \label{model}
\end{equation}
where $S$ denotes a deterministic signal matrix and $Y$ represents the random noise matrix.

We assumes that the signal $S$ admits a low-rank structure. Denote the singular value decomposition (SVD) of $S$ by
\begin{equation*}
    S = U D V^\top = \sum_{k = 1}^K d_k \bfk{u}_k \bfk{v}_k^\top,
\end{equation*}
where $D = \operatorname{diag} (d_1, \cdots, d_K)$ consists of the singular values of $S$, and
\begin{equation*}
    U = [\bfk{u}_1, \cdots, \bfk{u}_K] \in \mathbb{R}^{M \times K}
    \quad \text{and} \quad
    V = [\bfk{v}_1, \cdots, \bfk{v}_K] \in \mathbb{R}^{N \times K}
\end{equation*}
consist of the normalized left and right singular vectors of $S$, respectively. For the random noise $Y$, we consider a general covariance structure $\Sigma = \mathbb{E} YY^\top \in \mathbb{R}^{M \times M}$. To be specific, we assume
\begin{equation*}
    Y = \Sigma^{1/2} X
\end{equation*}
where $\Sigma=(\Sigma^{1/2})^2$ is a positive semi-definite matrix and $X = [x_{i \mu}] \in \mathbb{R}^{M \times N}$ has independent real random entries with zero mean and variance $1/N$.


This formulation appears ubiquitously in various fields. For instance, in wireless communication, wired communication, satellite communication and other systems, signals are disturbed by noise during transmission. By the above model \eqref{model}, one can evaluate the performance of the communication system, optimize the signal transmission strategy, and improve the reliability and stability of the signals \citep{edfors1998ofdm, Zeng2009}. In addition, in the field of signal processing, such as audio processing, image processing, noise is a common problem. By studying this model, different noise types can be processed to improve the quality and accuracy of the signals \citep{scharf1991statistical, Raychaudhuri1999, Konstantinides1997}.

Furthermore, in the field of high dimensional statistical inference and machine learning, the model also plays a very important role \citep{hastie2009elements,kannan2009spectral}. For example, it can be used for data reduction and feature extraction which are the core techniques of high dimensional data analysis.  By factorizing the data matrix into low-rank approximations, redundant information can be removed, computational complexity reduced, and performance of classification or clustering tasks may be improved \cite{JohnstoneandLu2009, Fanjianqing2011, Han2023}. Moreover, when each column of $S$ can be only chosen from finite number of distinct unknown deterministic vectors, \eqref{model} can be regarded as a mixture model. It is well known that mixture models are the foundation of heterogeneous data analysis and have gained widespread attention, especially in the current era of big data where heterogeneous data is more frequent \citep{MaHuang2017, WangFang2019, Shenhe2015}.

%


Motivated by the increasing prevalence of applications of this model, in this paper, we investigate the largest singular values of $\tilde{Y}$ in the regime where the data dimension $M$ and sample size $N$ are comparably large ($M/N \asymp 1$), while the rank of signal $K = \operatorname{rank} (S)$ remains bounded. Equivalently, we examine the largest eigenvalues of the sample covariance matrix
\begin{equation*}
    \tilde{Q}_{\mathrm{a}} := \tilde{Y} \tilde{Y}^\top
    = (S + Y) (S + Y)^\top.
\end{equation*}

From a theoretical standpoint, our formulation of signal-plus-noise model (\ref{model}) falls within the framework of the finite-rank deformation models in random matrix theory, which have received considerable attention over the past two decades. The subscript in $\tilde{Q}_{\mathrm{a}}$ indicates that the model (\ref{model}) is additive. As in many other finite-rank deformation models, the largest eigenvalues of $\tilde{Q}_{\mathrm{a}}$ (or the largest singular values of $\tilde{Y}$) also undergo a phase transition. This paper is devoted to the asymptotic distribution of the largest eigenvalues of $\tilde{Q}_{\mathrm{a}}$ that detaches from the bulk of its spectrum, which we refer to as \emph{spiked eigenvalues}. 

As pointed out in \cite{liuAsymptoticPropertiesSpiked2023} and demonstrated in the subsequent sections, the phase transition threshold and almost sure limits of the spiked eigenvalues obtained for the model (\ref{model}) are the same as those established for the popular sample covariance matrices 
\begin{equation}
    \tilde{Q}_{\mathrm{m}} = \tilde{\Sigma}^{1/2} X X^\top \tilde{\Sigma}^{1/2}
    \quad \text{ where } \quad
    \tilde{\Sigma}= \mathbb{E} \tilde{Q}_{\mathrm{a}} =SS^\top+\Sigma.
    \label{eqn-spiked-population}
\end{equation}
As a result, another motivation behind this paper is to ascertain whether this similarity extends to the level of the second order fluctuation. The main theorem of this paper dismisses this conjecture and provides insight about the distinction between the asymptotic distributions of the spiked eigenvalues of these two models. Let us proceed by offering a brief review of related works, after which we provide an overview of our results and the proof strategy employed in this paper.

\subsection{Related works}


In the field of random matrix theory, there are two well-established finite-rank deformation models that are closely related to our signal-plus-noise model: the spiked population model and deformed Wigner matrices. The sample covariance matrix $\tilde{Q}_{\mathrm{m}}$ associated with the former is given in (\ref{eqn-spiked-population}), where the underlying population covariance matrix $\tilde{\Sigma} \in \mathbb{R}^{M \times M}$ possesses a finite number of large eigenvalues detached from the bulk of its spectrum. Specifically, we can decompose
\begin{equation*}
    \tilde{\Sigma} = \Psi \operatorname{diag} (\tilde{\Lambda}_{\mathrm{s}}, \tilde{\Lambda}_{\mathrm{n}}) \Psi^\top,
\end{equation*}
so that $\tilde{\Lambda}_{\mathrm{s}}$ and $\tilde{\Lambda}_{\mathrm{n}}$ are diagonal matrices consisting of the spiked and non-spiked eigenvalues of $\tilde{\Sigma}$, respectively. This model, originally introduced by Johnstone \cite{johnstoneDistributionLargestEigenvalue2001}, can be viewed as the multiplicative analog of our signal-plus-noise matrix model (\ref{model}). Over the past two decades, many efforts have been devoted to the asymptotic behavior of the spiked eigenvalues of $\tilde{Q}_{\mathrm{m}}$ induced by $\tilde{\Lambda}_{\mathrm{s}}$. The seminal paper \cite{baikPhaseTransitionLargest2005} by Baik, Ben Arous and P\'{e}ch\'{e} established the well-known BBP transition for the largest eigenvalue $\lambda_1 (\tilde{Q}_{\mathrm{m}})$ in the case of complex Gaussian $X$ and diagonal $\tilde{\Sigma}$ with $\tilde{\Lambda}_{\mathrm{n}} = I$, i.e.
\begin{equation*}
    \tilde{\Sigma} = \operatorname{diag} (\tilde{\sigma}_1, \cdots, \tilde{\sigma}_K, 1, \cdots, 1).
\end{equation*}
Such a population covariance matrix is obtained when $SS^\top$ is diagonal and $\Sigma$ is the identity matrix in (\ref{eqn-spiked-population}).

Notably, the fluctuation of $\lambda_1 (\tilde{Q}_{\mathrm{m}})$ was proved to be asymptotic Gaussian in the supercritical regime. For real Gaussian $X$, the asymptotic Gaussian fluctuation of the spiked eigenvalues of $\tilde{Q}_{\mathrm{m}}$ was established by Paul \cite{paulAsymptoticsSampleEigenstructure2007}. For general $X$ whose entries are not necessarily Gaussian, Baik and Silverstein \cite{baikEigenvaluesLargeSample2006} investigated the almost sure convergence of extreme eigenvalues when $\tilde{\Lambda}_{\mathrm{n}} = I$. By assuming a block diagonal structure for $\tilde{\Sigma}$, Bai and Yao developed a central limit theorem (CLT) for spiked eigenvalues, initially for $\tilde{\Lambda}_{\mathrm{n}} = I$ in \cite{baiCentralLimitTheorems2008}, and later extended it to accommodate general $\tilde{\Lambda}_{\mathrm{n}}$ in \cite{baiSampleEigenvaluesGeneralized2012}. Recently, this block diagonal assumption on $\tilde{\Sigma}$ was removed by the works \cite{jiangGeneralizedFourMoment2021,zhangAsymptoticIndependenceSpiked2022}. Importantly, the derived asymptotic distribution is Gaussian and depends on the variables $\{ x_{i \mu} \}$ only via their first four cumulants (or moments). That is, we can consider this asymptotic distribution as being universal with respect to the distribution of entries in $X$. In contrast, as we will demonstrate in Section \ref{sec-main-results}, the fluctuations of the spiked eigenvalues associated with the signal-plus-noise model (\ref{model}) are nonuniversal and exhibit sensitivity to the distribution of entries in $X$.

The phenomenon of nonuniversality for outlier eigenvalues has previously been explored in the context of deformed Wigner matrices, which can be considered as the Hermitian analogue to our signal-plus-noise model. In this context, one examines the matrix $W + A$, where $W$ represents an $N \times N$ Wigner matrix, and $A$ is a deterministic matrix with a finite rank. Akin to what is observed in the spiked population model, the extreme eigenvalues of $W + A$ also undergo a BBP transition, which was initially demonstrated by P\'{e}ch\'{e} \cite{pecheLargestEigenvalueSmall2006} when $W$ is a Gaussian unitary ensemble (GUE). For general Wigner matrices, Capitaine, Donati-Martin and F\'{e}ral \cite{capitaineLargestEigenvaluesFinite2009} established the almost sure convergence of the extreme eigenvalues of $W + A$. Moreover, within the same study, the authors highlighted the nonuniversality phenomenon of the fluctuations of outlier eigenvalues by examining the specific case $A = d \mathbf{e}_1 \mathbf{e}_1^*$. This topic has been further investigated in \cite{capitaineCentralLimitTheorems2012,pizzoFiniteRankDeformations2013,renfrewFiniteRankDeformations2013}. In essence, it was found that the fluctuations of outlier eigenvalues in $W + A$ exhibit distinct behaviors depending on whether the eigenvectors of $A$ are localized or delocalized. These two regimes were unified by Knowles and Yin in \cite{knowlesIsotropicSemicircleLaw2013}, where they considered an arbitrary finite-rank deformation with bounded norm. Notably, in the case of rank-$1$ deformation $A = d \mathbf{v} \mathbf{v}^*$ with supercritical $d > 0$, it was shown that $\lambda_1(W + A)$ asymptotically follows the same distribution as a linear combination of $\mathbf{v}^* W \mathbf{v}$ and a Gaussian variable $\Phi$ that is independent of $W$. Moreover, the mean and variance of $\Phi$ can be explicitly computed from $A$ and the first four cumulants of $\{ w_{ij} \}$.

Compared to the extensive literature on extreme eigenvalues of the spiked population model and deformed Wigner matrices, there is a relative dearth of research on extreme singular values of additively deformed rectangular matrices, despite their widespread applications in various fields. Regarding the model examined in this paper, Benaych-Georges and Nadakuditi \cite{benaych-georgesSingularValuesVectors2012} characterized the BBP transition of the largest singular values of $\tilde{Y}$ by deriving their almost sure limit. In their work, $\bfk{u}_k$ and $\bfk{v}_k$ are generated by random vectors with i.i.d. entries or obtained through the Gram-Schmidt orthogonalization of these vectors. They also obtained a CLT for the largest singular value of $\tilde{Y}$ in cases where only a single signal is present ($K = 1$) and the common fourth cumulant of the entries of the random vectors vanishes. On the other hand, the complexity of the analysis escalates when the signals are permitted to have a more general structure. In such scenarios, Ding \cite{dingHighDimensionalDeformed2020} demonstrated that when $\Sigma = I$, the fluctuation of detached singular values of $\tilde{Y}$ around their almost sure limit is of order $N^{-1/2}$. More recently, \cite{liuAsymptoticPropertiesSpiked2023} determined the almost sure limit of the detached singular values, with only mild assumptions on the noise covariance matrix $\Sigma$. However, to the best of our knowledge, there are currently no theoretical results addressing the asymptotic distribution of the detached singular values of $\tilde{Y}$ under fairly weak assumptions on the signal $S$, even in the case where $\Sigma = I$. 

Finally, we remark that the recent work by Bao et al \cite{baoSingularVectorSingular2021} derived the asymptotic distribution of angles between the principal singular vectors of $\tilde{Y}$ and their deterministic counterparts $\bfk{u}_k$ and $\bfk{v}_k$ in the case of $\Sigma = I$. Our approach to characterizing the nonuniversality phenomenon associated with the fluctuation of the spiked eigenvalues of $\tilde{Q}_{\mathrm{a}} = \tilde{Y} \tilde{Y}^\top$ is inspired by their work.
 
\subsection{Overview of results}

The primary contribution of this paper is to establish the asymptotic joint distribution of the scaled fluctuations of the spiked eigenvalues $\lambda_k \equiv \lambda_k (\tilde{Q}_{\mathrm{a}})$ around their almost sure limit $\theta_k$. More precisely, we demonstrate that each fluctuation can be decomposed into three components,
\begin{equation}
    \sqrt{N} (\lambda_k - \theta_k) = \Phi_k + \Theta_k + \mathcal{L}_k.
\end{equation}
Here, $\mathcal{L}_k$ is deterministic and represents the asymptotic expectation of $\sqrt{N} (\lambda_k - \theta_k)$, while $\Theta_k$ captures the component contributing to the nonuniversality of the fluctuation. Upon removing these two components, the remaining component $\Phi_k$ exhibits an asymptotic Gaussian fluctuation with zero mean, and its variance depends on the entries $\{ x_{i \mu} \}$ only through their first four cumulants. We formulate this result in terms of the characteristic function of the joint distribution of the $\Phi_k$'s and $\Theta_k$'s. In essence, it asserts that there exist deterministic quantities $\mathcal{V}_{kj}$ and $\mathcal{W}_{kj}$ such that
\begin{equation*}
    \mathbb{E} \exp \Big [ {\mathrm{i} \sum_{k} (s_k \Phi_k + t_k \Theta_k)} \Big ]
    \approx \exp \Big [ {- \frac{1}{2} \sum_{k j} (s_k s_j \mathcal{V}_{kj}  + 2 s_k t_j \mathcal{W}_{kj})} \Big ]
    \cdot \mathbb{E} \exp \Big [ {\mathrm{i} \sum_{k} t_k \Theta_k} \Big ].
\end{equation*}
In particular, by letting $t_k = 0$, one readily deduce the asymptotic joint normality of the variables $\Phi_k$. It turns out that the quantity $\mathcal{V}_{kk}$ represents the asymptotic variance of $\Phi_k$, while $\mathcal{V}_{kj}$ characterizes the asymptotic covariance between $\Phi_k$ and $\Phi_j$. Another consequence of this result is that, the Gaussian component $\Phi_k$ is generally not independent of the nonuniversal component $\Theta_j$. Nevertheless, their mixed cumulants satisfy
\begin{equation*}
    \kappa_{p, q} (\Phi_k, \Theta_j) \to 
    \begin{cases}
        \mathcal{W}_{k j}, & \text{ if } p = q = 1, \\
        0, & \text{ if } p, q \geq 1 \text{ and } p + q \geq 3.
    \end{cases}
\end{equation*}
Here we recall that the multivariate cumulants of a bivariate vector $(\xi, \zeta)$ can be defined as the coefficients in Taylor expansion of its log-characteristic function,
\begin{equation*}
    \log \mathbb{E} \exp [ \mathrm{i} (s \xi + t \zeta) ]
    = \sum_{p, q = 0}^\infty \kappa_{p, q} (\xi, \zeta)
    \frac{(\mathrm{i} s)^p (\mathrm{i} t)^q}{p! q!}.
\end{equation*}
In other words, asymptotically, the dependence between $\Phi_k$ and $\Theta_j$ can be characterized by their covariance $\mathcal{W}_{k j}$. We highlight that establishment of this main result relies on a recursive estimate for the quantities $\mathbb{E} [ \Phi^{\ell} \mathrm{e}^{\mathrm{i} \Theta} ]$, which may have an independent interest. Essentially, we show that for each fixed $\ell \geq 0$, the following recurrence relation that bears a resemblance to the pattern of Hermite polynomials holds,
\begin{equation}
    \mathbb{E} [ {\Phi^{\ell + 2} \mathrm{e}^{\mathrm{i} \Theta}} ]
    \approx (\ell + 1) \mathcal{V} \mathbb{E} [ {\Phi^{\ell} \mathrm{e}^{\mathrm{i} \Theta}} ]
    + (\mathrm{i} \mathcal{W}) \mathbb{E} [ {\Phi^{\ell + 1} \mathrm{e}^{\mathrm{i} \Theta}} ].
    \label{eqn-approx-recurrence-moment}
\end{equation}

To summarise, this paper addresses the asymptotic joint distribution of the spiked eigenvalues associated with the signal-plus-noise model (\ref{model}), under fairly general assumptions on the structure of the pair $(S, \Sigma)$ and the distribution of the random entries $\{ x_{i \mu} \}$. Of particular significance, we discover the presence of nonuniversality in the fluctuations of these spiked eigenvalues around their almost sure limit, which does not emerge in the context of sample covariance matrices $\tilde{Q}_{\mathrm{m}}$ associated with the spiked population model. From a theoretical perspective, our contributions include: (i) a more adaptable approach for deriving the Green function representation of spiked eigenvalues in finite-rank deformation models, (ii) a general result regarding the asymptotic distribution of centralized sesquilinear forms of the resolvent $G(z)$, and (iii) an interesting recurrence relation for the joint moments of two interdependent variables, which proves valuable in characterizing the nonuniversality of spiked eigenvalues.

The subsequent sections of this paper are structured as follows. Section \ref{sec-main-results} is to present the principal findings of this study along with the necessary technical assumptions. To demonstrate the practical applicability of our results, Section \ref{sec-application} illustrates their use in solving a detection problem associated with mean heterogeneity. In section \ref{sec-proof}, we elucidate the two crucial steps leading to our main theorem. The majority of the technical proof is available in the Supplementary Material.

\begin{notation}
We denote the $k$-th largest eigenvalue of a symmetric or Hermitian matrix $A$ as $\lambda_k(A)$. For simplicity, we often abbreviate $\lambda_k \equiv \lambda_k (\tilde{Q}_{\mathrm{a}})$ whenever there is no ambiguity. Throughout the paper, we regard $N$ as the fundamental parameter and take $M \equiv M^{(N)}$. For simplicity, we frequently suppress the dependence on $N$ from the notations, bearing in mind that all quantities that are not explicitly constant may depend on $N$. Given the $N$-dependent quantities $a \equiv a^{(N)}$ and $b \equiv b^{(N)}$, we write $a \lesssim b$ if $a^{(N)} \leq C b^{(N)}$ for some $N$-independent constant $C > 0$. We also use the notation $a \asymp b$ if we have simultaneously $a \lesssim b$ and $a \gtrsim b$. Given positive integers $N_1 < N_2$, we denote $\llbracket N_1, N_2 \rrbracket := \{ N_1, N_1 + 1, \cdots, N_2 \}$ and abbreviate $\llbracket N_2 \rrbracket \equiv \llbracket 1, N_2 \rrbracket$ for simplicity. We reserve the letter $i$ to represent the indices within $\llbracket M \rrbracket$, while the symbol $\mu$ is reserved for indices within either $\llbracket N \rrbracket$ or $\llbracket M+1, M+N \rrbracket$, dependent on the specific context. Hence, there should be no ambiguity when we denote the canonical bases of $\mathbb{R}^M$ and $\mathbb{R}^N$ as $\{ \mathbf{e}_i \}_{i \in \llbracket M \rrbracket}$ and $\{ \mathbf{e}_\mu \}_{\mu \in \llbracket N \rrbracket}$, respectively. We use $(\mathbf{a})_k$ to represent the $k$th component of a vector $\mathbf{a}$. Moreover, we frequently identify vectors $\mathbf{u}, \mathbf{r} \in \mathbb{R}^M$ and $\mathbf{v}, \mathbf{w} \in \mathbb{R}^{N}$ and also use them to represent their respective natural embeddings in $\mathbb{R}^{M + N}$ when there is no confusion, i.e.
\begin{equation}
    \mathbf{u} \equiv \begin{bmatrix}
        \mathbf{u} \\ \mathbf{0}
    \end{bmatrix}, 
    \mathbf{r} \equiv \begin{bmatrix}
        \mathbf{r} \\ \mathbf{0}
    \end{bmatrix} \in \mathbb{R}^{M + N}
    \quad \text{ and } \quad
    \mathbf{v} \equiv \begin{bmatrix}
        \mathbf{0} \\ \mathbf{v}
    \end{bmatrix},
    \mathbf{w} \equiv \begin{bmatrix}
        \mathbf{0} \\ \mathbf{w}
    \end{bmatrix} \in \mathbb{R}^{M + N}.
    \label{def-embedding}
\end{equation}
Finally, throughout the paper, $\| \cdot \|$ denotes either the operator norm of a matrix or the $\ell_2$-norm of a vector. 
\end{notation}

\section{Main results} \label{sec-main-results}

\subsection{Assumptions}

Throughout this paper, we let $\tau > 0$ be some positive constant (independent of $N$) which can be chosen arbitrarily small. 

\begin{assumption}\label{assumption-bounded-moments}
We assume that $\{ x_{i \mu} \}$ are independent real random variables satisfying
\begin{equation}
    \mathbb{E} x_{i \mu} = 0, 
    \quad \mathbb{E} \big | \sqrt{N} x_{i \mu} \big |^2 = 1
    \qand
    \mathbb{E} \big | \sqrt{N} x_{i \mu} \big |^p \leq C_p < \infty,
    \quad \forall p \geq 3.
\end{equation}
\end{assumption}

Here the constants $C_{p}$ are assumed to be independent of $N$. We denote the $p$-th cumulant of the variables by $\kappa_p^{i \mu} := \kappa_p (\sqrt{N} x_{i \mu})$. It follows from Assumption \ref{assumption-bounded-moments} that
\begin{equation*}
    \kappa_1^{i \mu} = 0, 
    \quad \kappa_2^{i \mu} = 1
    \qand
    \kappa_p^{i \mu} \leq C_{p}^\prime < \infty,
    \quad \forall p \geq 3,
\end{equation*}
where the constants $C_{p}^\prime$ are independent of $N$. For the sake of clarity and presentation simplicity, in the following we state the theorems under the assumption that the variables share identical third and forth cumulants (or i.i.d). However, it is straightforward to relax this constraint by keeping track of the cumulants $\kappa_3^{i \mu}$ and $\kappa_4^{i \mu}$ when employing the cumulant expansion formula. For the necessary modifications to the main theorem when this restriction is lifted, we refer the readers to remark \ref{nonhomogenous} below. 

Let $\sigma_1 \geq \sigma_2 \geq \cdots \geq \sigma_M \geq 0$ be the eigenvalues of population covariance $\Sigma$ of the noise $Y$. We denote the corresponding empirical spectral distribution by
\begin{equation*}
    \nu \equiv \nu^{(N)} := \frac{1}{M} \sum_{i=1}^M \delta_{\sigma_i},
\end{equation*}

\begin{assumption}\label{assumption-dimension}
$\phi \equiv \phi_N := {M} / {N} \in [\tau, \tau^{-1}]$ for all sufficiently large $N$.
\end{assumption}

\begin{assumption}\label{assumption-spectrum}
$\| \Sigma \| = \sigma_1 \leq \tau^{-1}$ and $\nu([0, \tau]) \leq 1 - \tau$.
\end{assumption}

Assumption \ref{assumption-dimension} implies that we are operating within the high-dimensional regime. Assumption \ref{assumption-spectrum} ensures that the spectrum of $\Sigma$ remains uniformly bounded and avoids concentration of eigenvalues near zero. We now characterize the asymptotic spectral distribution of undeformed covariance matrix 
\begin{equation*}
    Q = YY^\top
\end{equation*}
via its Stieltjes transform. The following results are well-known and their proof can be found in \cite{silversteinAnalysisLimitingSpectral1995,knowlesAnisotropicLocalLaws2017}. For each $z \in \mathbb{C}_+$, there is a unique $m (z) \in \mathbb{C}_+$ satisfying
\begin{equation}
    z = - \frac{1}{m} + \phi \int \frac{s}{1 + s m} \nu (\mathrm{d} s).
    \label{eqn-self-consistent}
\end{equation}
Moreover, $m(z)$ is the Stieltjes transform of a probability measure $\varrho$ on $[0, \infty)$, i.e.
\begin{equation*}
    m(z) = \int \frac{1}{s - z} \varrho (\mathrm{d} s).
\end{equation*}
We call $\varrho$ the asymptotic spectral distribution, which is uniquely determined by $m$ using the inversion formula $\varrho (\mathrm{d} s) = {\pi}^{-1} \lim_{\eta \downarrow 0} \operatorname{Im} m(s + \mathrm{i} \eta)$. In particular, when $\nu = \delta_1$, i.e. $\Sigma = I$, the asymptotic spectral distribution $\varrho$ reduces to the prominent Marchenko-Pastur law. Given Assumptions \ref{assumption-dimension} and \ref{assumption-spectrum}, the supports of $\varrho$ remain uniformly bounded in $N$. Also note that $m (z)$ can be analytically extended to $\mathbb{C} \backslash \operatorname{supp} \varrho$. 

To describe the rightmost edge of $\varrho$, we introduce the function $f$, which serves as the inverse of $m$ outside the spectrum. Specifically, set
\begin{equation}
    f(w)
    := - \frac{1}{w} + \phi \int \frac{s}{1 + s w} \nu (\mathrm{d} s)
    = - \frac{1}{w} + \frac{1}{N} \operatorname{tr} \frac{\Sigma}{1 + w \Sigma},
\end{equation}
and it is analytic on $\{ w \in \mathbb{C}: -w^{-1} \notin \operatorname{supp} \nu \}$. The function $f(w)$ possesses a unique stationary point in $(-\sigma^{-1}_1, 0)$, which we denote as $w_{+}$, i.e. $w_{+}$ is the unique solution to
\begin{equation}
    f^\prime (w) 
    = \frac{1}{w^2} - \phi \int \left ( \frac{s}{1 + s w} \right )^2 \nu (\mathrm{d} s)
    = 0,
    \quad
    w \in (-\sigma^{-1}_1, 0).
    \label{def-critical-point}
\end{equation}
Now, the rightmost endpoint of $\operatorname{supp} \varrho$ can be characterized by $\lambda_{+} = f(w_{+})$. In addition, we have $w_{+} = m(\lambda_{+})$, where $m(\lambda_{+})$ is defined via continuous extension. Actually, when restricted to $[\lambda_+, \infty)$, the Stieltjes transform $m$ is negative, strictly increasing and invertible, with $f$ being its inverse. To prevent the emergence of large spiked eigenvalues in the spectrum of $\Sigma$, we impose the following regularity condition on the rightmost edge.

\begin{assumption}\label{assumption-regularity}
$w_{+} + \sigma_1^{-1} = m(\lambda_{+}) + \sigma_1^{-1} \geq \tau$.
\end{assumption}

As demonstrated by \cite{knowlesAnisotropicLocalLaws2017}, Assumptions \ref{assumption-bounded-moments}-\ref{assumption-regularity} ensure the asymptotic Tracy-Widom (TW) fluctuation for the largest eigenvalues of $Q$. Next, let us discuss the technical assumptions concerning the additive deformation $S$. 

\begin{assumption}\label{assumption-deformation}
The rank of the deformation matrix, $K$, is bounded. The strengths of deformation satisfy $d_k \in [\tau, \tau^{-1}]$ for all $k \in \llbracket K \rrbracket$. Here $\{ d_k \}_{k=1}^K$ are allowed to depend on $N$.
\end{assumption}

The population covariance of the deformed ensemble $\tilde{Y}$ is given by
\begin{equation*}
    \tilde{\Sigma} 
    := \mathbb{E} \tilde{Q}_{\mathrm{a}}
    = \mathbb{E} \tilde{Y} \tilde{Y}^\top
    = \Sigma + U D^2 U^\top.
\end{equation*}
The eigenvalues of $\tilde{\Sigma}$ play a crucial role in our theorem, which are denoted as 
\begin{equation*}
    \tilde{\sigma}_1 \geq \tilde{\sigma}_2 \geq \cdots \geq \tilde{\sigma}_M.
\end{equation*}
For each $k \in \llbracket K \rrbracket$, we denote the unit eigenvector of $\tilde{\Sigma}$ associated with $\tilde{\sigma}_k$ as $\boldsymbol{\psi}_k$, i.e.
\begin{equation*}
    \tilde{\Sigma} \boldsymbol{\psi}_k = \tilde{\sigma}_k \boldsymbol{\psi}_k, 
    \quad \text{ where } \quad
    \| \boldsymbol{\psi}_k \| = 1.
\end{equation*}
Here the selection of the sign of the eigenvectors $\boldsymbol{\psi}_k$ has no effect on our results. Note that by Cauchy's Interlace Theorem we have
\begin{equation*}
    \sigma_i \leq \tilde{\sigma}_i \leq \sigma_{i - K}, 
    \quad \forall K+1 \leq i \leq M.
\end{equation*}
In particular, the eigenvalues $\{ \tilde{\sigma}_i \}_{i = K+1}^M$ constitute the bulk of spectrum of $\tilde{\Sigma}$. On the other hand, $\{ \tilde{\sigma}_k \}_{k = 1}^K$ are the eigenvalues that can potentially induce spiked eigenvalues of $\tilde{Q}_{\mathrm{a}}$. Similar to other finite-rank deformation models, the largest eigenvalues of $\tilde{Q}_{\mathrm{a}}$ also undergo a BBP transition. The critical threshold for this transition is given by $-w_{+}^{-1}$. That is, to make the sample eigenvalues $\lambda_k$ jump out of the $\operatorname{supp} \varrho$, one needs to impose
\begin{equation}\label{cri1}
\tilde{\sigma}_k > -w_{+}^{-1}.
\end{equation}
The derivation of this threshold is presented in Section \ref{subsec-BBP}. For ease of presentation, in this paper we do not delve into the regime $\tilde{\sigma}_k + w_{+}^{-1} \ll 1$, although we believe that our arguments can be extended to the regime $\tilde{\sigma}_k + w_{+}^{-1} \gg N^{-1/3}$ by utilizing the corresponding local law. 

\begin{assumption}\label{assumption-spiked-spacing}
$\tilde{\sigma}_k - \tilde{\sigma}_{k+1} \geq \tau$ for all $k \in \llbracket K_0 \rrbracket$, where
 \begin{equation} \
    K_0 = \max \{ k : \tilde{\sigma}_k \geq -w_{+}^{-1} + 2\tau \}.
    \label{def-K0}
\end{equation}
We assume that $K_0\geq 1$.
\end{assumption}

The condition (\ref{def-K0}) ensures that there exists at least one spiked sample eigenvalue and Assumption \ref{assumption-spiked-spacing} also concerns the spacing between the largest eigenvalues of $\tilde{\Sigma}$. Again, our argument remains applicable when Assumption \ref{assumption-spiked-spacing} is relaxed to $\tilde{\sigma}_k - \tilde{\sigma}_{k+1} \gg N^{-1/2}$. However, in cases where some of the outliers are overlapping in the sense that $\tilde{\sigma}_k - \tilde{\sigma}_{k+1} \ll N^{-1/2}$, an intricate partitioning of the outliers is required. We leave the exploration of this scenario to future works.

\subsection{Asymptotic distribution of spiked eigenvalues}
\label{subsec-asymptotic-distribution}

To rigorously formulate our result, we adopt the following notion of high probability bounds from \cite{erdosLocalSemicircleLaw2013} to systematize statements of the form ``$\mathcal{X}$ is bounded with high probability by $\mathcal{Y}$ up to small powers of $N$''.
\begin{definition}[stochastic domination]
Consider two families of nonnegative random variables parameterized by $N \in \mathbb{N}$ and $t \in \mathbb{T}^{(N)}$,
\begin{equation*}
    \mathcal{X} \equiv \mathcal{X}^{(N)}(t)
    \qand
    \mathcal{Y} \equiv \mathcal{Y}^{(N)}(t),
\end{equation*}
where $\mathbb{T}^{(N)}$ is a possibly $N$-dependent parameter set. We say that $\mathcal{X}$ is \emph{stochastically dominated} by $\mathcal{Y}$, uniformly in $t$, if for all (small) $\varepsilon > 0$ and (large) $L > 0$ we have
\begin{equation*}
    \sup_{t \in \mathbb{T}^{(N)}} 
    \mathbb{P} \big \{ \mathcal{X}^{(N)}(t) > N^{\varepsilon} \mathcal{Y}^{(N)}(t) \big \} \leq N^{-L},
    \quad \forall N \geq N_0(\varepsilon, L).
\end{equation*}
\end{definition}
In this paper, $N_0$ may depend on parameters that are explicitly fixed (such as $\tau$ and $C_p$). We adopt the following notation conventions:
\begin{itemize}
    \item If $\mathcal{X}$ is stochastically dominated by $\mathcal{Y}$, uniformly in $t$, we denote $\mathcal{X} \prec \mathcal{Y}$.
    \item If $\mathcal{X}$ is complex and $|\mathcal{X}| \prec \mathcal{Y}$, we write $\mathcal{X} = \oprec(\mathcal{Y})$. 
    \item Let $A \equiv A^{(N)}(t)$ be a family of $n \times n$ random matrices where the dimension $n$ is independent of $N$. We write $A = \oprec(\mathcal{Y})$ if $| a_{ij} | \prec \mathcal{Y}$ uniformly in $t$ and $i,j \in [n]$.
\end{itemize}

The first result of this paper concerns the almost sure limit of the spiked eigenvalues of $\tilde{Q}_{\mathrm{a}}$. Let us begin with the following definition: as the Stieltjes transform $m$ is invertible on $\mathbb{C} \backslash \operatorname{supp} \varrho$ with its inverse given by $f$, we can introduce
\begin{equation}
    \theta (\tilde{\sigma}) := f(- \tilde{\sigma}^{-1})
    = \tilde{\sigma} + \frac{1}{N} \operatorname{tr} \frac{\tilde{\sigma} \Sigma}{\tilde{\sigma} - \Sigma},
    \quad \text{ for } \quad
    \tilde{\sigma} > - w_+^{-1}.
\end{equation}
Recall our definition of $K_0$ in (\ref{def-K0}). For each $k \in \llbracket K_0 \rrbracket$, we denote
\begin{equation}
    \theta_k := \theta (\tilde{\sigma}_k) 
    = \tilde{\sigma}_k + \frac{1}{N} \operatorname{tr} \frac{\tilde{\sigma}_k \Sigma}{\tilde{\sigma}_k - \Sigma}
    \qand
    \theta_k^\prime := \theta^\prime (\tilde{\sigma}_k) 
    = 1 - \frac{1}{N} \operatorname{tr} \frac{\Sigma^2}{(\tilde{\sigma}_k - \Sigma)^2},
    \label{def-theta-k}
\end{equation}
where $\theta^\prime$ represents the derivative of $\theta$ and the deterministic value $\theta_k$ represents the almost sure limit of the spiked eigenvalue $\lambda_k$. We direct interested readers to Section \ref{subsec-BBP} for the rationale behind the definition of $\theta_k$. The following result provides a large deviation bound for $\lambda_k - \theta_k$. Especially, this result extends the analysis conducted in \cite{dingHighDimensionalDeformed2020} to accommodate noises with a general covariance structure and the result in \cite{liuAsymptoticPropertiesSpiked2023} with a convergence rate. 

\begin{proposition}[large deviation bound] 
\label{prop-convergence-rate}
Under Assumptions \ref{assumption-bounded-moments}-\ref{assumption-spiked-spacing}, we have
\begin{equation}
    \lambda_k - \theta_k = \oprec(N^{-1/2}),
    \quad \forall k \in \llbracket K_0 \rrbracket.
\end{equation}
\end{proposition}

We note that the error bound provided in Proposition \ref{prop-convergence-rate} is optimal up to an $N^\varepsilon$ correction introduced by the definition of $\prec$. In particular, for a meaningful study of the asymptotic distribution of $\lambda_k - \theta_k$, one should scale it by $\sqrt{N}$. Before presenting our main result concerning the asymptotic joint distribution of these scaled fluctuations, let us introduce some necessary notations. We set $\mathbf{1}_N = (1, \cdots, 1) \in \mathbb{R}^N$. For each $k \in \llbracket K_0 \rrbracket$, we introduce a vector $\tilde{\boldsymbol{\pi}}_k \in \mathbb{R}^M$ such that its entries are given by
\begin{equation}
    (\tilde{\boldsymbol{\pi}}_k)_{i} = \left ( \frac{\Sigma}{\tilde{\sigma}_k - \Sigma} \right )_{i i}.
    \label{def-tilde-pi-vec}
\end{equation}
Let $\mathbf{a}_1, \cdots, \mathbf{a}_n \in \mathbb{R}^T$ be vectors of common length. Given positive integers $p_1, \cdots, p_n$, we use $\mathbb{M}$ to denote the mixed moment of $\mathbf{a}_1, \cdots, \mathbf{a}_n$ of order $p_1, \cdots, p_n$, i.e.
\begin{equation*}
    \mathbb{M}_{p_1, \cdots, p_n} (\mathbf{a}_1, \cdots, \mathbf{a}_n)
    = \sum_{t=1}^T [(\mathbf{a}_1)_t]^{p_1} \cdots [(\mathbf{a}_n)_t]^{p_n}.
\end{equation*}
We suppress the subscripts of $\mathbb{M}$ if $p_1 = \cdots = p_n = 1$. Note that $\mathbb{M} (\mathbf{a}_1, \mathbf{a}_2)$ is simply the inner product $\mathbf{a}_1^\top \mathbf{a}_2$. For ease of reference, in the following definition we summarize the quantities that are used to characterize the asymptotic distribution of the spiked eigenvalues.

\begin{definition}[asymptotic quantities for Theorem \ref{thm-spiked-distribution}] \label{def-spiked-quantities}
For each $k \in \llbracket K_0 \rrbracket$, define
\begin{equation*}
    \mathcal{L}_k = \frac{2 \kappa_{3} \theta^\prime_k}{N} 
    \tilde{\boldsymbol{\pi}}_k^\top \Sigma^{1/2} \boldsymbol{\psi}_k
    \cdot \mathbf{1}_N^\top S^\top \boldsymbol{\psi}_k.    
\end{equation*}
For each $k, j \in \llbracket K_0 \rrbracket$, we introduce
$$\mathcal{V}_{kj} = \mathcal{V}_{kj}^{(0, 1, 0)} + \mathcal{V}_{kj}^{(1, 2, 0)},$$ where
\begin{align*}
    \mathcal{V}_{kj}^{(0, 1, 0)} 
    & = \begin{cases}
        2 (\theta^\prime_k)^2 (\boldsymbol{\psi}_k ^\top \Sigma \boldsymbol{\psi}_k)^2
        + 2 \tilde{\sigma}_k^2 \theta^\prime_k - 2 \tilde{\sigma}_k^2 (\theta^\prime_k)^2,
        & \text{ if } k = j, \\
        2 \theta^\prime_k \theta^\prime_j (\boldsymbol{\psi}_k ^\top \Sigma \boldsymbol{\psi}_j)^2,
        & \text{ if } k \not= j,
    \end{cases} \\
    \mathcal{V}_{kj}^{(1, 2, 0)}
    & = \frac{\kappa_4 \theta^\prime_k \theta^\prime_j}{N} \Big [ 
    \mathbf{1}_N^\top \mathbf{1}_N \cdot \mathbb{M}_{2, 2} (\Sigma^{1/2} \boldsymbol{\psi}_k, \Sigma^{1/2} \boldsymbol{\psi}_j)
    + \tilde{\boldsymbol{\pi}}_k^\top \tilde{\boldsymbol{\pi}}_j \cdot \mathbb{M}_{2, 2} (S^\top \boldsymbol{\psi}_k, S^\top \boldsymbol{\psi}_j) \Big ].
\end{align*}
Finally, for each $k, j \in \llbracket K_0 \rrbracket$, let
\begin{equation*}
    \mathcal{W}_{kj}
    = \frac{2 \kappa_3 \theta^\prime_k \theta^\prime_j}{\sqrt{N}} \Big [ 
    \mathbf{1}_N^\top S^\top \boldsymbol{\psi}_j
    \cdot \mathbb{M}_{2, 1} (\Sigma^{1/2} \boldsymbol{\psi}_k, \Sigma^{1/2} \boldsymbol{\psi}_j )
    + \tilde{\boldsymbol{\pi}}_k^\top \Sigma^{1/2} \boldsymbol{\psi}_j
    \cdot \mathbb{M}_{2, 1} (S^\top \boldsymbol{\psi}_k, S^\top \boldsymbol{\psi}_j) \Big ].
\end{equation*}
\end{definition}

\begin{remark}\label{nonhomogenous}
Importantly, these quantities exhibit universality with respect to the entry distribution of $X$ as they only involve the first four cumulants of $\{ x_{i \mu} \}$. Moreover, each quantity involving $\kappa_3$ or $\kappa_4$ can be expressed as a product of two summations: one over $i \in \llbracket M \rrbracket$, and another over $\mu \in \llbracket N \rrbracket$, e.g.
\begin{equation*}
    \mathcal{L}_k = \frac{2 \kappa_{3} \theta^\prime_k}{N} 
    \sum_{i} (\tilde{\boldsymbol{\pi}}_k)_i (\Sigma^{1/2} \boldsymbol{\psi}_k)_i
    \sum_{\mu} (S^\top \boldsymbol{\psi}_k)_\mu.
\end{equation*}
Upon a meticulous examination of the proof, it becomes evident that to ease the constraint of the identical third and fourth cumulants, we just need to couple the two summations and treat $\kappa_3^{i \mu}$ or $\kappa_4^{i \mu}$ as a factor of the summand. For example, in that case $\mathcal{L}_k$ should be 
\begin{equation*}
    \mathcal{L}_k = \frac{2 \theta^\prime_k}{N} 
    \sum_{i \mu} \kappa_{3}^{i \mu} (\tilde{\boldsymbol{\pi}}_k)_i (\Sigma^{1/2} \boldsymbol{\psi}_k)_i (S^\top \boldsymbol{\psi}_k)_\mu.
\end{equation*}
\end{remark}

\begin{theorem}[fluctuation of spiked eigenvalues]
\label{thm-spiked-distribution}
Suppose that Assumptions \ref{assumption-bounded-moments}-\ref{assumption-spiked-spacing} are satisfied. Let $\mathcal{L}_{k}, \mathcal{V}_{kj}, \mathcal{W}_{kj}$ be given as in Definition \ref{def-spiked-quantities}. For each $k \in \llbracket K_0 \rrbracket$, define
\begin{equation*}
    \Theta_k := 2 \sqrt{N} \theta^\prime_k 
    \boldsymbol{\psi}_k^\top \Sigma^{1/2} X S^\top \boldsymbol{\psi}_k
    \quad \text{ and } \quad
    \Phi_k = \sqrt{N} ( \lambda_k - \theta_k ) - \Theta_k - \mathcal{L}_k.
\end{equation*}
Let 
$$\Phi := \sum_{k = 1}^{K_0} s_k \Phi_k,\quad \Theta := \sum_{k = 1}^{K_0} t_k \Theta_k$$
where $|s_k|, |t_k| \lesssim 1$ are deterministic. Then, given any $\varepsilon > 0$, there exists $N_0 \equiv N_0 (\varepsilon)$ such that for all $N \geq N_0$,
\begin{equation}
    \bigg | {\mathbb{E} \exp [{\mathrm{i} (\Phi + \Theta)}]
    - \exp \Big ( {- \frac{\mathcal{V} + 2 \mathcal{W}}{2}} \Big )
    \mathbb{E} \mathrm{e}^{\mathrm{i} \Theta}} \bigg | \leq \epsilon,
    \label{eqn-spiked-distribution}
\end{equation}
where 
$$\mathcal{V} = \sum_{k,j = 1}^{K_0} s_k s_j \mathcal{V}_{kj},\quad \mathcal{W} = \sum_{k,j = 1}^{K_0} s_k t_j \mathcal{W}_{kj}.$$
\end{theorem}


\begin{remark}\label{exclude}
We would like to point out that the sample spiked eigenvalues are due to the signal matrix $S$ rather than the noise matrix $Y$ by Assumptions \ref{assumption-regularity} and \ref{assumption-spiked-spacing}. Hence this excludes the case where $S=0$.
\end{remark}

\begin{remark}\label{asymptotic independ}
By setting $t_k = 0$ we readily establish the asymptotic joint normality of the $\Phi_k$'s, where $\mathcal{V}_{kj}$ represents the asymptotic covariance between $\Phi_k$ and $\Phi_j$. On the other hand, the dependence between $\Phi_k$ and $\Theta_j$ is explained by their asymptotic covariance $\mathcal{W}_{kj}$. We note that if $\kappa_3 = 0$, which holds especially when the distributions of the $x_{i \mu}$'s are symmetric, then $\mathcal{W}_{kj}$ automatically becomes zero and consequently, the families $\{ \Phi_k \}$ and $\{ \Theta_k \}$ are asymptotically independent.
\end{remark}

\begin{remark}\label{nonuniversal}
The nonuniversal part $\Theta_k$ in Theorem \ref{thm-spiked-distribution} is explicitly given by
    \begin{equation*}
        \Theta_k 
        = 2 \sqrt{N} \theta^\prime_k \sum_{i \mu}
        (\Sigma^{1/2} \boldsymbol{\psi}_k)_i (S^\top \boldsymbol{\psi}_k)_\mu
        x_{i \mu}. 
    \end{equation*}
    Hence, the characteristic function $\mathbb{E} \mathrm{e}^{\mathrm{i} \Theta}$ in (\ref{eqn-spiked-distribution}) actually admits an explicit form in terms of the characteristic functions of the $x_{i \mu}$'s. Moreover, according to the Lindeberg-L\'{e}vy CLT, if either $\Sigma^{1/2} \boldsymbol{\psi}_k$ or $S^\top \boldsymbol{\psi}_k$ is delocalized, i.e. 
    \begin{equation*}
        \| \Sigma^{1/2} \boldsymbol{\psi}_k \|_{\infty} \wedge \| S^\top \boldsymbol{\psi}_k \|_{\infty} = o(1),
    \end{equation*} 
    then asymptotically $\Theta_k$ exhibits Gaussian fluctuation. Here $\| \cdot \|_{\infty}$ is the supremum norm. On the other hand, if both $\Sigma^{1/2} \boldsymbol{\psi}_k$ and $S^\top \boldsymbol{\psi}_k$ are (nearly) localized, i.e. 
    \begin{equation*}
        \| \Sigma^{1/2} \boldsymbol{\psi}_k \|_{\infty} \wedge \| S^\top \boldsymbol{\psi}_k \|_{\infty} \gtrsim 1,
    \end{equation*}
    then the distribution of $\Theta_k$ is strongly influenced by the law of $x_{i \mu}$'s.
\end{remark}

\begin{remark}
As an corollary, we examine the special case where the noise part has an isotropic covariance $\Sigma = I$. In this scenario, the phase transition condition $\tilde{\sigma}_k > -w_{+}^{-1}$ reduces to $d_k^4 > \phi$. By a straightforward computation, we can also find that
\begin{equation*}
    \theta_k = 1 + d_k^2 + \phi \left ( 1 + \frac{1}{d_k^2} \right )
    \qand
    \theta_k^\prime = 1 - \frac{\phi}{d_k^4}.
\end{equation*}
If the signal intensities are ordered non-increasingly, i.e., $d_1 > d_2 > \cdots > d_K$, we can deduce that $\boldsymbol{\psi}_k = \bfk{u}_k$ by definition. Consequently, the nonuniversal components become
\begin{equation*}
    \Theta_k = 2 \sqrt{N} \theta^\prime_k d_k
    \bfk{u}_k^\top X \bfk{v}_k.
\end{equation*}
Furthermore, the asymptotic quantities given in Definition \ref{def-spiked-quantities} reduce to
\begin{align*}
    \mathcal{L}_k & = \frac{2 \kappa_{3} \theta^\prime_k}{N d_k} 
    (\mathbf{1}_M^\top \bfk{u}_k)
    (\mathbf{1}_N^\top \bfk{v}_k), \\
    \mathcal{V}_{kj}^{(0, 1, 0)} 
    & = 2 \theta^\prime_k (1 + \phi + 2\phi / d_k^2) \mathbbm{1}_{k = j}, \\
    \mathcal{V}_{kj}^{(1, 2, 0)}
    & = \frac{\kappa_4 \theta^\prime_k \theta^\prime_j}{N} \Big [ 
    (\mathbf{1}_N^\top \mathbf{1}_N) \mathbb{M}_{2, 2} (\bfk{u}_k, \bfk{u}_j)
    + (\mathbf{1}_M^\top \mathbf{1}_M) \mathbb{M}_{2, 2} (\bfk{v}_k, \bfk{v}_j) \Big ], \\
    \mathcal{W}_{kj}
    & = \frac{2 \kappa_3 \theta^\prime_k \theta^\prime_j d_j}{\sqrt{N}} \Big [ 
    (\mathbf{1}_N^\top \bfk{v}_j)
    \mathbb{M}_{2, 1} (\bfk{u}_k, \bfk{u}_j)
    + (\mathbf{1}_M^\top \bfk{u}_j)
    \mathbb{M}_{2, 1} (\bfk{v}_k, \bfk{v}_j) \Big ].
\end{align*}
\end{remark}

As highlighted in Section \ref{sec-introduction}, we can make comparisons between the behaviors of the spiked eigenvalues in our signal-plus-noise model (\ref{model}) and the spiked population model, particularly when the two models possess the same population covariance matrix. Recall $\tilde{\Sigma} = \Sigma + S S^\top$. We therefore compare the sample covariance matrices
\begin{equation*}
    \tilde{Q}_{\mathrm{a}} = (S + \Sigma^{1/2} X) (S + \Sigma^{1/2} X)^\top 
    \qand
    \tilde{Q}_{\mathrm{m}} = \tilde{\Sigma}^{1/2} X X^\top \tilde{\Sigma}^{1/2},
\end{equation*}
where $X$ satisfies Assumption \ref{assumption-bounded-moments}. Clearly, we have $\mathbb{E} \tilde{Q}_{\mathrm{a}} = \mathbb{E} \tilde{Q}_{\mathrm{m}} = \tilde{\Sigma}$. Leveraging Proposition \ref{prop-convergence-rate} alongside its counterpart as presented in \cite{baiSampleEigenvaluesGeneralized2012}, these two models share the same phase transition threshold and almost sure limit for supercritical eigenvalues. On the other hand, \cite{jiangGeneralizedFourMoment2021,zhangAsymptoticIndependenceSpiked2022} show that the asymptotic distributions of the sample spiked eigenvalues of $\tilde{Q}_{\mathrm{m}}$ are Gaussian and depends on the variables $\{ x_{i \mu} \}$ only via their first four cumulants (or moments). In other words, this asymptotic distribution is universal with respect to the distributions of entries in $X$.  In contrast, Theorem \ref{thm-spiked-distribution} shows that the asymptotic distributions of the sample spiked eigenvalues of $\tilde{Q}_{\mathrm{a}}$ may be nonuniversal, and dependent on the population eigenvectors and the distributions of entries in $X$. 

\begin{figure}
    \centering
    \includegraphics[width=\textwidth]{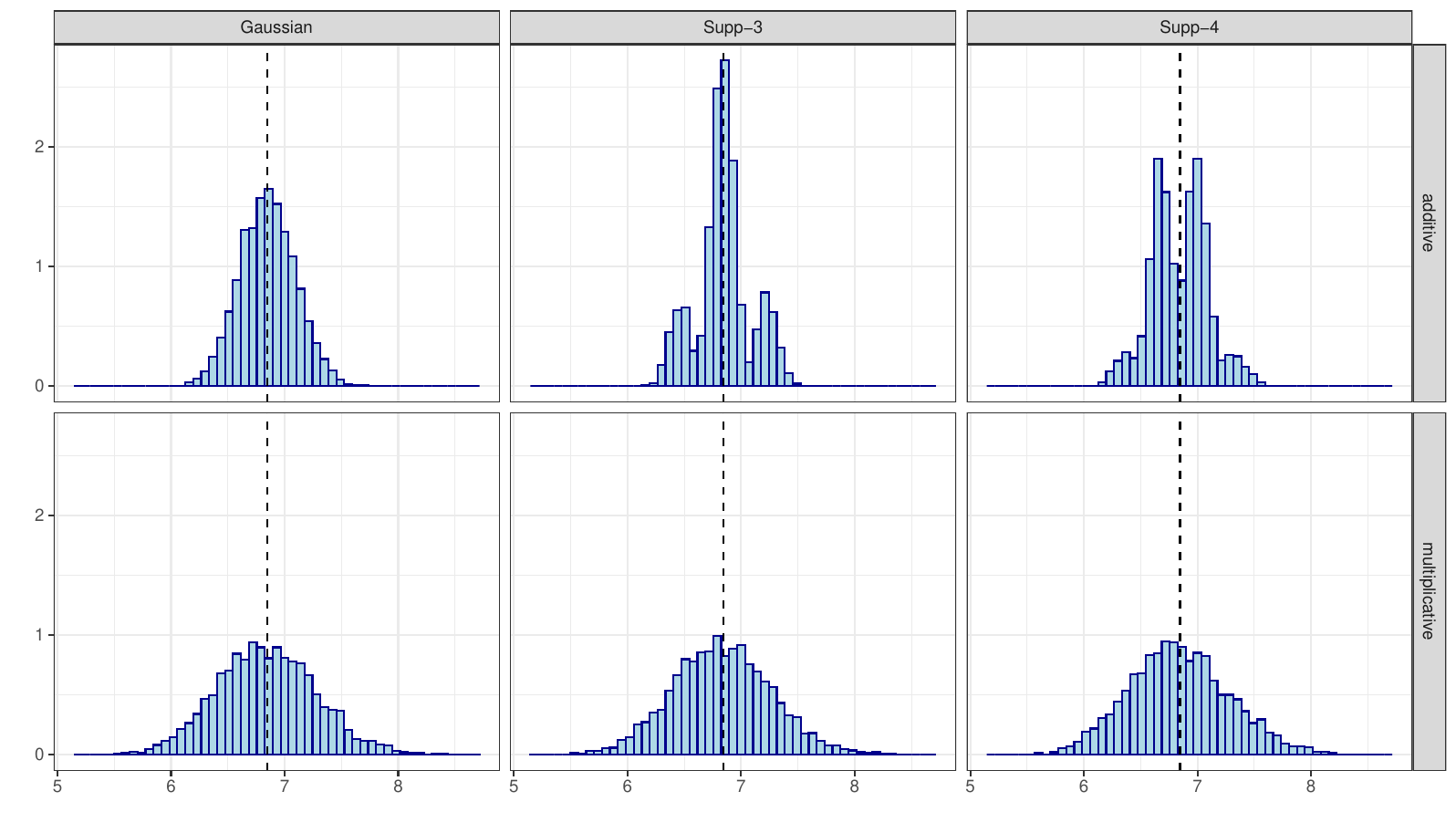}
    \caption{Comparison of spiked eigenvalues between signal-plus-noise model (Top) and spiked population model (Bottom). Parameters: 5,000 iterations, $M=200$, $N=400$, $\Sigma = I$ and $(S)_{i \mu} = \sqrt{5.25} \mathbbm{1} (i=1 \text{ and } \mu=1)$. Left: Standard Gaussian. Middle: $\mathbb{P}(x = \pm \sqrt{3}) = 1/6$ and $\mathbb{P}(x = 0) = 2/3$. Right: $\mathbb{P}(x = \pm 1/\sqrt{2}) = 4/9$ and $\mathbb{P}(x = \pm \sqrt{5}) = 1/18$. The latter two are tailored to match the first four moments of the standard Gaussian.}
    \label{fig-comparison}
\end{figure}

A visual comparison that illustrates this distinction is depicted in Figure \ref{fig-comparison}, where we conducted simulations with $\Sigma = I$ and $(S)_{i \mu} = \sqrt{5.25} \mathbbm{1} (i=1 \text{ and } \mu=1)$. In our simulations, we select three distributions for ${ x_{i \mu} }$, all of which share identical first four moments. As anticipated, the empirical distribution of the spiked eigenvalue in $\tilde{Q}_{\mathrm{m}}$ exhibits similarity across these three distributions. In contrast, the behavior of the spiked eigenvalue in $\tilde{Q}_{\mathrm{a}}$ is noticeably influenced by the specific characteristics of the distributions of ${ x_{i \mu} }$ beyond their first four moments. From a technical perspective, the key reason behind this distinction is that the Green function representation of $\lambda_k (\tilde{Q}_{\mathrm{a}})$ involves the off-diagonal blocks of the Green function, while $\lambda_k (\tilde{Q}_{\mathrm{m}})$ only involves the top-left block of the Green function. We refer interested readers to Section \ref{subsec-sesquilinear-forms} for more details.

\subsection{Application to mixture models}

As previously mentioned, a model that aligns seamlessly with our signal-plus-noise framework (\ref{model}) is the \emph{mixture model}. In this model, the distribution of the data under investigation takes the form of a mixture involving $K$ distributions, each characterized by a distinct mean. Specifically, assume that the $N$ i.i.d. observations are given by $\mathbf{a}_{\mu} + \Sigma^{1/2} \mathbf{w}_{\mu}$, where $\mathbf{a}_{\mu} \in \mathbb{R}^M$ are i.i.d. random vectors satisfying
\begin{equation*}
    \mathbb{P} (\mathbf{a}_{\mu} = \mathbf{c}_k) = \alpha_k,
    \quad \forall k \in \llbracket K \rrbracket.
\end{equation*}
Here $\sum_{k=1}^K \alpha_k = 1$. In other words, $\alpha_k$ signifies the likelihood that an observation pertains to the $k$-th cluster, and $\mathbf{c}_k$ represents the center of the $k$-th cluster. As for the noise vectors $\mathbf{x}_{\mu} = \mathbf{w}_{\mu} / \sqrt{N}$, we assume they are independent of the $\mathbf{a}_{\mu}$'s and have entries satisfying Assumption \ref{assumption-bounded-moments}. We note that this formulation encompasses the widely studied Gaussian mixture model.

Without loss of generality, we can condition on the randomness of the $\mathbf{a}_{\mu}$'s. Let $\mathcal{N}_k = \{ \mu \in \llbracket N \rrbracket: \mathbf{a}_{\mu} = \mathbf{c}_k \}$, and define $\mathbf{1}_{\mathcal{N}_k} \in \mathbb{R}^N$ as the vector with entries $(\mathbf{1}_{\mathcal{N}_k})_{\mu} = \mathbbm{1}(\mu \in \mathcal{N}_k)$. Now the scaled signal matrix can be expressed as
\begin{equation*}
    S = \frac{1}{\sqrt{N}} [\mathbf{c}_1, \cdots, \mathbf{c}_K] [\mathbf{1}_{\mathcal{N}_1}, \cdots, \mathbf{1}_{\mathcal{N}_K}]^\top.
\end{equation*}
Note that the low-rank assumption imposed on $S$ holds as long as the number of clusters remains bounded, and the boundness of $\| S \|$ is implied by the boundness of $\max_{k \in \llbracket K \rrbracket} \| \mathbf{c}_k \|$. Now, let us compute the right singular vectors of $S$ and see how they affect the behavior of the spiked eigenvalues. Here we set $\Sigma = I$ to simplify the discussion. Let $\hat{\alpha}_k = |\mathcal{N}_k| / N$, which represents the sample proportion of the $k$-th cluster. Assume that we have the singular value decomposition
\begin{equation*}
    [\sqrt{\hat{\alpha}_1} \mathbf{c}_1, \cdots, \sqrt{\hat{\alpha}_K} \mathbf{c}_K] 
    = U D R^\top
    = [\bfk{u}_1, \cdots, \bfk{u}_K] D R^\top,
\end{equation*}
where $D$ is diagonal and $R = (r_{jk})_{j,k=1}^K$ satisfies $R^\top R = I$. (For the moment, let us temporarily ignore the possibility that the rank of $S$ is strictly smaller than $K$ as it is not essential to the discussion here.) It follows that
\begin{equation*}
    S = U D V^\top,
    \quad \text{ where } \quad
    V = \left [ \frac{\mathbf{1}_{\mathcal{N}_1}}{\sqrt{N \hat{\alpha}_1}}, \cdots, \frac{\mathbf{1}_{\mathcal{N}_K}}{\sqrt{N \hat{\alpha}_K}} \right ] R.
\end{equation*}
Note that $V^\top V = I$. In other words, the right singular vectors of $S$ have the form
\begin{equation*}
    \bfk{v}_k = \sum_{j = 1}^K \frac{r_{j k} \mathbf{1}_{\mathcal{N}_j}}{\sqrt{N \hat{\alpha}_j}}.
\end{equation*}
In particular, when $N \hat{\alpha}_j \to \infty$ for all $j \in \llbracket K \rrbracket$, one has $\| \bfk{v}_k \|_{\infty} = o(1)$. According to the remark \ref{nonuniversal}, 
in this scenario, the spiked eigenvalues (if they exist) of the associated sample covariance matrix exhibit asymptotic Gaussian fluctuation. The requirement of $N \hat{\alpha}_j \to \infty$ is relatively weak, as it is commonly assumed in the literature that $\alpha_j \wedge (1 - \alpha_j) \asymp 1$ for all $j \in \llbracket K \rrbracket$, that is, the proportion of the $j$-th cluster is bounded away from both $0$ and $1$. In this case, according to the law of large numbers, we have $\| \bfk{v}_k \|_{\infty} \asymp N^{-1/2}$ almost surely. 

On the other hand, we would like to interpret the nonuniversality of spiked eigenvalues in the signal-plus-noise model (\ref{model}) from the perspective of the presence of abnormal or atypical observations that deviate significantly from the majority of observations. Specifically, 
one could show that that nonuniversality of spiked eigenvalues arises only when there are ``clusters'' within the dataset, that comprises only a small fraction of the total sample ($|\mathcal{N}_j| = N \hat{\alpha}_j \lesssim 1$). Furthermore, for these ``clusters'' to exert a substantial impact on the spiked eigenvalues, the corresponding observations must possess considerable magnitudes ($\| \mathbf{c}_j \| \gtrsim \sqrt{N}$). In such scenarios, it may be more appropriate to treat these observations as anomalies rather than a constituent of the mixture. It is worth mentioning that the simulation in Figure \ref{fig-comparison} serves as an illustrative example of this situation.

Summarising the discussion above we conclude the following:
\begin{itemize} 
\item The nonuniversality phenomenon does not arise for the sample spiked eigenvalues in the mixture model as long as the number of observations in each cluster diverges as $N \to \infty$. \item 
If there is at least a cluster with the number of observations $|\mathcal{N}_j| \lesssim 1$ out of $K$ clusters and $\| \mathbf{c}_j \| \gtrsim \sqrt{N}$, then the sample spiked eigenvalues exhibit nonuniversality.
\end{itemize}

\section{Detection of mean heterogeneity}
\label{sec-application}

In this section, we leverage the findings established in Section \ref{sec-main-results} to address a detection problem concerning mean heterogeneity. Assessing homogeneity holds significant importance as an initial stride before delving into subsequent data analysis endeavors, since numerous classical statistical methodologies hinge upon assumptions of homogeneity. In accordance with the mixture model introduced in Section \ref{sec-main-results}, our observations are defined as
\begin{equation*}
    \mathbf{a}_{\mu} + \Sigma^{1/2} \mathbf{w}_{\mu}
    \quad \text{ where } \quad
    \mathbb{P} (\mathbf{a}_{\mu} = \mathbf{c}_k) = \alpha_k.
\end{equation*}
In this context, the problem of detecting mean heterogeneity can be formulated as the following hypothesis test
\begin{equation}
    \mathrm{H}_0: K = 1
    \quad  \mathrm{ v.s. } \quad  
    \mathrm{H}_1: 1 < K \leq K_{\star},
    \label{hypothesis-test}
\end{equation}
where $K_{\star} \in \mathbb{N}$ represents an upper bound for the total number of clusters. This hypothesis testing problem has garnered considerable attention in literature, with noteworthy studies such as \cite{Chenjiahua1998,ChenChen2001,Sunjiayang2004,ChenLi2009,Chenjiahua2020}. These works primarily concentrate on scenarios where the data dimensionality $M$ remains fixed. However, the classical methodologies devised in these studies, like EM-test or likelihood ratio tests, will encounter formidable obstacles in high-dimensional settings where $M$ scales with the growth of $N$. In the high-dimensional regime, \cite{VERZELEN2017} introduced a set of eigenvalue-based and moment-based statistics designed specifically for Gaussian mixture models with $K=2$. Notably, their methodology was grounded in the assumption of a sparse difference between the centers of these two clusters. Consequently, our objective is to devise an eigenvalue-based statistic that transcends these limitations and is applicable to a broader array of mixture models.

For simplicity, we make the assumption that $\mathbb{E} \mathbf{a}_{\mu} = \sum_{k=1}^K \alpha_k \mathbf{c}_k = \mathbf{0}$. When $K=1$, this constraint implies that $\mathbf{c}_1 = \mathbf{0}$. Furthermore, under the alternative hypothesis, we implicitly assume that the intensity of signal cannot be excessively weak. In other words, the largest eigenvalue of $\tilde{\Sigma} = \Sigma + S S^\top$ must exceed the critical threshold $-\omega_{+}^{-1}$. Now, the scaled data matrix is given by
\begin{equation*}
    \tilde{Y} = [\tilde{\mathbf{y}}_1, \cdots, \tilde{\mathbf{y}}_N],
    \quad  \text{ where } \quad  
    \tilde{\mathbf{y}}_\mu = \frac{1}{\sqrt{N}} (\mathbf{a}_{\mu} + \Sigma^{1/2} \mathbf{w}_{\mu}).
\end{equation*}
Let $\lambda_i \equiv \lambda_i (\tilde{Y}^\top \tilde{Y})$. As indicated in \cite{knowlesAnisotropicLocalLaws2017}, under the null hypothesis, the first $K_{\star}$ largest eigenvalues (appropriately centred and normalized)
\begin{equation*}
    \frac{N^{2/3}}{\sigma_{\mathrm{TW}}} 
    (\lambda_1 - \lambda_{+}, \cdots, \lambda_{K^\star} - \lambda_{+})
\end{equation*}
converges weakly toward the Tracy-Widom-Airy random vector of type $1$. Here the scaled factor is given by $\sigma_{\mathrm{TW}} := (f^{\prime \prime} (\lambda_+)/2)^{1/3}$. In contrast, under the alternative hypothesis, the scaled fluctuation $N^{1/2} (\lambda_1 - \theta_1)$ converges to the distribution as elucidated in Theorem \ref{thm-spiked-distribution}. Furthermore, via a perturbation argument as demonstrated in \cite{benaych-georgesFluctuationsExtremeEigenvalues2011,bloemendalPrincipalComponentsSample2016}, $\lambda_{K^\star}$ remains adhered to the rightmost edge $\lambda_+$ and fluctuates at a level of $N^{-2/3}$. In view of this, we propose a new statistic for the hypothesis test (\ref{hypothesis-test}),
\begin{equation*}
    \mathrm{DS}_{K_{\star}}=\frac{\lambda_1-\lambda_{K_{\star}}}{\lambda_{K_{\star}}-\lambda_{2K_{\star}-1}}.
\end{equation*}
This test statistic is adapted from the one proposed by \cite{Onatski2006}, originally designed for testing the presence of factors,
\begin{equation*}
   \mathrm{RS}_{K_{\star}} = \frac{\lambda_1-\lambda_{K_{\star}}}{\lambda_{K_{\star}}-\lambda_{K_{\star}+1}}. 
\end{equation*}
These statistics circumvent the estimation of $\sigma_{\mathrm{TW}}$, which typically relies heavily on the spectrum of the unknown population covariance $\Sigma$. In accordance with Theorem \ref{thm-spiked-distribution}, both $\mathrm{DS}_{K_{\star}}$ and $\mathrm{RS}_{K_{\star}}$ exhibit a power that approaches $1$ as $N \to \infty$. However, as illustrated in the subsequent finite-sample analysis, the power of $\mathrm{RS}_{K_{\star}}$ diminishes significantly when the sample size is not sufficiently large and the signal strength is relatively weak.

To evaluate the finite-sample performance of the aforementioned statistics, we use $K_{\star}=4$ as an example for conducting simulation studies. To calculate the critical values of $\mathrm{DS}_{K_{\star}}$ and $\mathrm{RS}_{K_{\star}}$ through Monte Carlo simulations, we can construct Wishart matrices from $N_\star \times N_\star$ random matrices whose entries are i.i.d. Gaussian with zero mean and variance $1/N_\star$. Certainly, the choice of $N_\star$ should be tailored to the specific problem. If the data matrix is sufficiently large, for example, when $M \wedge N \geq 300$, we can use the critical values simulated using $N_\star = 1000$, as done in \cite{Onatski2006}. While for relatively small values of $M \wedge N$, it is recommended to adjust the value of $N_\star$ accordingly.

For our purpose, we examine two distinct pairs of $(N,M)$ values: $(N,M)=(200,100)$ and $(N,M)=(100,200)$. Given the relatively small values of $N$ and $M$, we set $N_\star = M \wedge N = 100$. The calculation of the corresponding critical values is conducted through the following steps:
\begin{enumerate}
    \item Generate a Wishart matrix from a $100 \times 100$ random matrix whose entries are i.i.d. Gaussian with zero mean and variance $1/100 = 0.01$, and compute its $(K_\star + 1)$ largest eigenvalues, denoted as $\nu_1, \ldots, \nu_{K_{\star}+1}$.
    \item Calculate the following two statistics
    \begin{equation*}
        t_1 = (\nu_1-\nu_{4})/(\nu_{4}-\nu_{7})
        \qand
        t_2 = (\nu_1-\nu_{4})/(\nu_{4}-\nu_{7}).
    \end{equation*}
    \item Repeat the above steps a total of 30,000 times independently. This yields two sequences: $T_1 = \{ t_{1,s} \}$ and $T_2 = \{ t_{2,s} \}$. We use the $95\%$ percentile of $T_1$ and $T_2$ as the approximate critical values for $\mathrm{DS}_4$ and $\mathrm{RS}_4$, respectively.
\end{enumerate}
Based on the results from our Monte Carlo simulations, the approximate critical values for $\mathrm{DS}_4$ and $\mathrm{RS}_4$ are given by $3.0251$ and $18.9920$, respectively.

\subsection{Size}

In this subsection we examine the finite-sample properties of the two test statistics with a focus on their sizes. Under the null hypothesis $K=1$, we have $\tilde{\mathbf{y}}_\mu = \Sigma^{1/2} \mathbf{w}_\mu / \sqrt{N}$. We consider three distinct forms of the covariance matrix $\Sigma$:
\begin{itemize}
  \item $\Sigma_1 = I$.
  \item $\Sigma_2 = \{\sigma_{ij}\}$, where $\sigma_{ij}=0.1^{|i-j|}$ for $1\leq i\leq j\leq M$.
  \item $\Sigma_3 = O \operatorname{diag} (\sigma_1, \ldots, \sigma_M) O^T$, where $\sigma_i \overset{\mathrm{i.i.d.}}{\sim} \operatorname{Unif} (1, 1.5)$ and $O \in \mathbb{R}^{M \times M}$ follows the Haar measure on the orthogonal group.
\end{itemize}
 In addition, we simultaneously consider two different distributions of $\mathbf{w}_\mu$, that is,
 \begin{itemize}
  \item The entries of $\mathbf{w}_\mu$ are i.i.d. from $N(0,1)$.
  \item The entries of $\mathbf{w}_\mu$ are i.i.d. from $\operatorname{Unif}(-\sqrt{3},\sqrt{3})$.
\end{itemize}
 
Conducting the experiment with 10,000 iterations, we record the ratio of rejections as an empirical measure of size. The outcomes of this investigation are consolidated in Table \ref{table2}. It is evident that all the statistics demonstrate strong performance in terms of size, with type I errors closely aligned with the predetermined value of $0.05$.
\begin{table}[ht]
\caption{Empirical size of different statistics.}
\label{table2} 
\centering
\begin{tabular}{ccc|cc|cc}
    \hline
    &		    &          &\multicolumn{2}{c|}{$N(0,1)$}  &\multicolumn{2}{c}{$U(-\sqrt{3},\sqrt{3})$} \\
    &		    &Method    &$(200,100)$     &$(100,200)$    &$(200,100)$            & $(100,200)$  \\
    \hline
    &	$\Sigma_1$	&$\mathrm{DS}_4$   &0.0475     &0.0460     &0.0465     &0.0493      \\
    &               &$\mathrm{RS}_4$   &0.0422     &0.0474     &0.0526     &0.0483         \\
    \hline
    &	$\Sigma_2$	&$\mathrm{DS}_4$   &0.0497     &0.0435     &0.0482     &0.0494      \\
    &               &$\mathrm{RS}_4$   &0.0521    &0.0511     &0.0531     &0.0543         \\
    \hline
    &	$\Sigma_3$	&$\mathrm{DS}_4$   &0.0499     &0.0475     &0.0463      &0.0495      \\
    &               &$\mathrm{RS}_4$   &0.0478     &0.0496     &0.0468     &0.0477         \\
    \hline
\end{tabular}
\end{table}

\subsection{Power}

In this subsection we turn to the power of these statistics. We explore three scenarios: $K = 2, 3, 4$, respectively. For simplicity we put $\alpha_k = 1/K$ and consider the three types of $\Sigma$ and two types of $\mathbf{w}_\mu$ introduced in the preceding subsection. Regarding the centers of the clusters, the configurations are outlined below for different values of $K$:
\begin{itemize}
    \item For $K = 2$, we assign $\mathbf{c}_1$ with i.i.d. entries drawn from $\operatorname{Unif} (0, 0.3)$, and set $\mathbf{c}_2 = -\mathbf{c}_1$.
    \item For $K = 3$, we assign $\mathbf{c}_1$ with i.i.d. entries drawn from $\operatorname{Unif} (0, 0.4)$, $\mathbf{c}_2$ with i.i.d. entries drawn from $\operatorname{Unif} (-0.3, 0)$, and set $\mathbf{c}_3=-(\mathbf{c}_1+\mathbf{c}_2)$.
    \item For $K = 4$, we assign $\mathbf{c}_1$ with i.i.d. entries drawn from $\operatorname{Unif} (0, 0.45)$, $\mathbf{c}_2$ with i.i.d. entries drawn from $\operatorname{Unif} (-0.3, 0)$, $\mathbf{c}_3$ with i.i.d. entries drawn from $\operatorname{Unif} (-0.1, 0.2)$, and set $\mathbf{c}_4=-(\mathbf{c}_1+\mathbf{c}_2+\mathbf{c}_3)$.
\end{itemize}

For each combination of $N$, $M$, and $\Sigma$ and $\mathbf{w}_\mu$, we conducted the experiments 10,000 times and recorded the rejection ratios to demonstrate the empirical power of different statistics. The summarized results are presented in Tables \ref{table3}, \ref{table4}, and \ref{table5} for $K = 2,3,4$ respectively. Notably, the statistics $\mathrm{DS}_{4}$ exhibit a significantly higher level of power compared to the statistics $\mathrm{RS}_{4}$. One plausible explanation for this phenomenon could be that the critical values of $\mathrm{RS}_{4}$ are relatively larger when contrasted with those of $\mathrm{DS}_{4}$. While this might not pose a significant issue when $N$ tends to infinity, it could considerably diminish the power of $\mathrm{RS}_{4}$ when dealing with finite values of $N$.

\begin{figure}
    \centering
    \includegraphics[width=\textwidth]{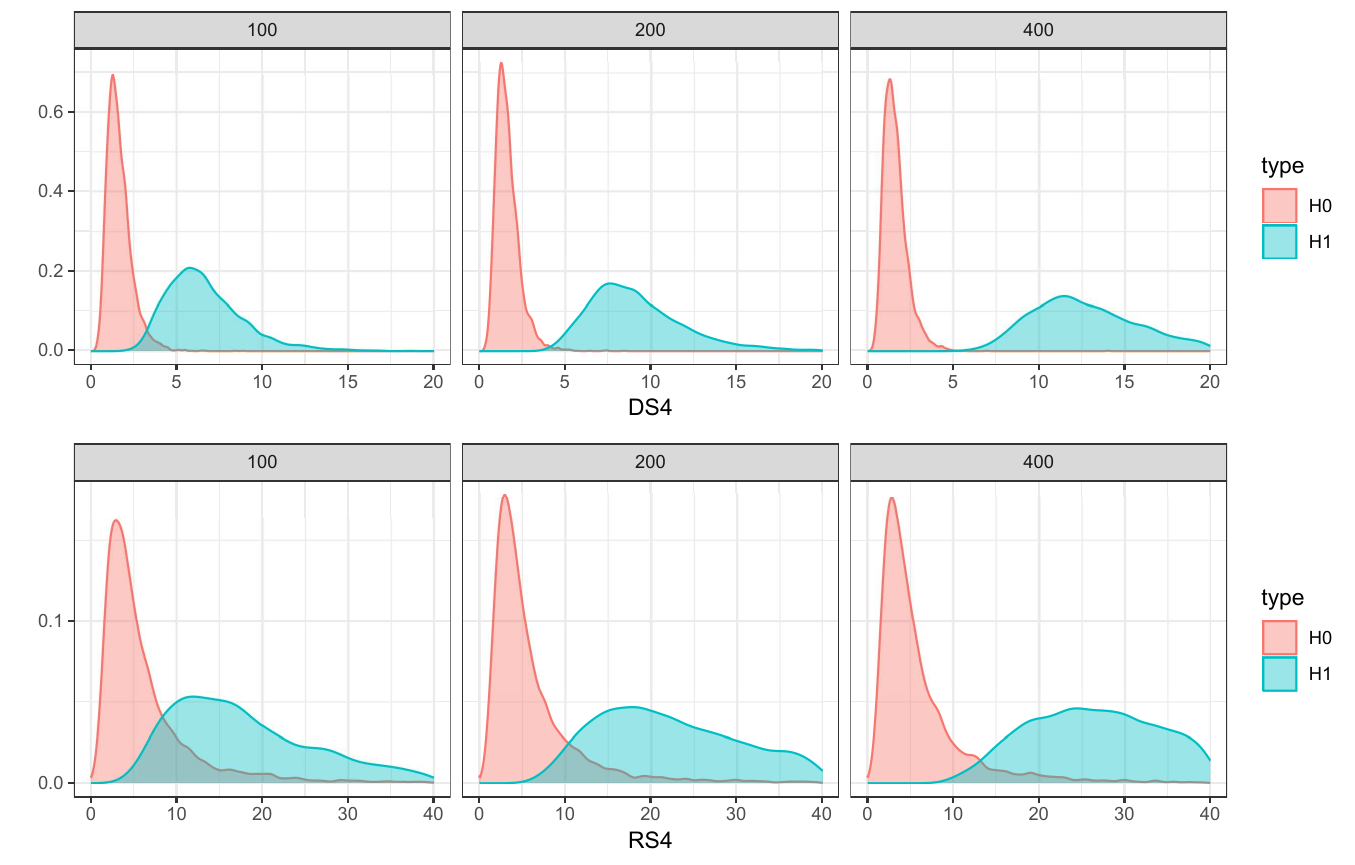}
    \caption{Distribution of statistics $\mathrm{DS}_4$ (Top) and $\mathrm{RS}_4$ (Bottom) under the null hypothesis (red) and alternative hypothesis (blue). Parameters: 5,000 iterations, $M=100, 200, 400$, $N=2M$, and $\Sigma = I$. The entries of $\mathbf{w}_\mu$ are i.i.d. from $N(0,1)$. In the simulation under $\mathrm{H}_1$, we specify $K = 2$ with $\mathbf{c}_1 = (1.5, 0, \cdots, 0)$ and $\mathbf{c}_2 = - \mathbf{c}_1$.}
    \label{fig-mod-stat}
\end{figure}

A visual representation that highlights the distinction in the tail behavior of the two statistics under $\mathrm{H}_0$ can be observed in Figure \ref{fig-mod-stat}. Notably, the tail of the distribution for $\mathrm{RS}_{4}$ under $H_0$ is substantially thicker than that of $\mathrm{DS}_{4}$ (please take note of the different $x$-axis scales for each statistic). This thick-tail property of $\mathrm{RS}_{4}$ under $\mathrm{H}_0$ leads to a significant overlap with the distribution under $\mathrm{H}_1$, particularly when the sample size is not sufficiently large. As $N \to \infty$, the distribution of $\mathrm{RS}_{4}$ under $H_1$ will gradually shift along the positive $x$-axis, causing the two distributions to separate. However, in finite-sample analysis, this overlap cannot be disregarded.

\begin{table}[ht]
\caption{Empirical power for $K=2$}
\label{table3} 
\centering
\begin{tabular}{ccc|cc|cc}
    \hline
    &		    &          &\multicolumn{2}{c|}{$N(0,1)$}  &\multicolumn{2}{c}{$U(-\sqrt{3},\sqrt{3})$} \\
    &		    &Method    &$(200,100)$     &$(100,200)$    &$(200,100)$         & $(100,200)$  \\
    \hline
    &	$\Sigma_1$	&$\mathrm{DS}_4$   &1            &1              &0.9998     &0.9999      \\
    &               &$\mathrm{RS}_4$   &0.7099       &0.7138         &0.7251     &0.7375         \\
    \hline
    &	$\Sigma_2$	&$\mathrm{DS}_4$   &0.9999     &1                &0.9998     &0.9999      \\
    &               &$\mathrm{RS}_4$   &0.7001     &0.7082            &0.7275     &0.7375         \\
    \hline
    &	$\Sigma_3$	&$\mathrm{DS}_4$   &0.9909     &0.9914           &0.9932      &0.9956      \\
    &               &$\mathrm{RS}_4$   &0.4089     &0.4907           &0.4789      &0.5044        \\
    \hline
\end{tabular}
\end{table}

\begin{table}[ht]
\caption{Empirical power for $K=3$}
\label{table4}
\centering
\begin{tabular}{ccc|cc|cc}
    \hline
    &		    &          &\multicolumn{2}{c|}{$N(0,1)$}  &\multicolumn{2}{c}{$U(-\sqrt{3},\sqrt{3})$} \\
    &		    &Method    &$(200,100)$     &$(100,200)$    &$(200,100)$         & $(100,200)$  \\
    \hline
    &	$\Sigma_1$	&$\mathrm{DS}_4$   &0.9984            &1              &0.9989     &1      \\
    &               &$\mathrm{RS}_4$   &0.5592            &0.9263         &0.5598     &0.9251         \\
    \hline
    &	$\Sigma_2$	&$\mathrm{DS}_4$   &0.9982     &1                &0.9998     &0.9987      \\
    &               &$\mathrm{RS}_4$   &0.5505     &0.9126            &0.6420     &0.5736         \\
    \hline
    &	$\Sigma_3$	&$\mathrm{DS}_4$   &0.9433     &1                &0.9605      &0.9999      \\
    &               &$\mathrm{RS}_4$   &0.3463     &0.7650           &0.3629      &0.6778        \\
    \hline
\end{tabular}
\end{table}

\begin{table}[ht]
\caption{Empirical power for $K=4$}
\label{table5} 
\centering
\begin{tabular}{ccc|cc|cc}
    \hline
    &		    &          &\multicolumn{2}{c|}{$N(0,1)$}  &\multicolumn{2}{c}{$U(-\sqrt{3},\sqrt{3})$} \\
    &		    &Method    &$(200,100)$     &$(100,200)$    &$(200,100)$         & $(100,200)$  \\
    \hline
    &	$\Sigma_1$	&$\mathrm{DS}_4$   &0.9999            &1              &0.9997     &1      \\
    &               &$\mathrm{RS}_4$   &0.6565            &0.9667        &0.6736     &0.9706         \\
    \hline
    &	$\Sigma_2$	&$\mathrm{DS}_4$   &1          &1                &0.9999     &1      \\
    &               &$\mathrm{RS}_4$   &0.6477     &0.9667            &0.6735     &0.9644        \\
    \hline
    &	$\Sigma_3$	&$\mathrm{DS}_4$   &0.9774     &0.9996           &0.8420      &1      \\
    &               &$\mathrm{RS}_4$   &0.4601     &0.6409           &0.2702      &0.6603        \\
    \hline
\end{tabular}
\end{table}

\section{Proof of the main results}
\label{sec-proof}

In this section we turn to the proof of the results provided in Section \ref{sec-main-results}. The proof of the large deviation bound, Proposition \ref{prop-convergence-rate}, is presented in the Supplementary Material. Here we focus on elucidating the two steps that lead to our main theorem regarding the asymptotic distribution of spiked eigenvalues. Throughout this section, we consistently present our results under Assumptions \ref{assumption-bounded-moments}-\ref{assumption-spiked-spacing}. 

\subsection{Local law}
\label{subsec-BBP}

Before proceeding to formally present the two aforementioned steps, let us start with some definitions and preliminary results that are necessary for our proof. We adopt the self-adjoint linearization technique, which is commonly employed in the literature of rectangular matrices. For $z \notin \operatorname{supp} \varrho$, let
\begin{equation*}
    \tilde{H}(z) := \begin{bmatrix}
    0 & \sqrt{z} \tilde{Y} \\
    \sqrt{z} \tilde{Y}^\top & 0
    \end{bmatrix}
    \qand
    H(z) := \begin{bmatrix}
    0 & \sqrt{z} Y \\
    \sqrt{z} Y^\top & 0
    \end{bmatrix}.
\end{equation*}
Below, we frequently abbreviate $\tilde{H} \equiv \tilde{H} (z)$, $H \equiv H(z)$ when there is no ambiguity. Note that we have
\begin{equation*}
    \tilde{H} = H + \sqrt{z} \mathfrak{U} \mathfrak{D} \mathfrak{U}^\top
    \quad \text{ where } \quad
    \mathfrak{U} := \begin{bmatrix}
    U & 0 \\
    0 & V
    \end{bmatrix}
    \qand
    \mathfrak{D} := \begin{bmatrix}
    0 & D \\
    D & 0
    \end{bmatrix}.
\end{equation*}
Recall $Y = \Sigma^{1/2} X$. Hence, we can express $H$ in the form
\begin{equation*}
    H = \sqrt{z} \underline{\Sigma}^{1/2} \begin{bmatrix}
    0 & X \\
    X^\top & 0
    \end{bmatrix} \underline{\Sigma}^{1/2},
    \quad \text{ where } \quad
    \underline{\Sigma} := \begin{bmatrix}
    \Sigma & 0 \\
    0 & I
    \end{bmatrix}.
\end{equation*}
By introducing the $(M+N) \times (M+N)$ rank-$2$ matrices
\begin{equation}
    \Delta^{i \mu} := \mathbf{e}_{i} \mathbf{e}_{\mu}^\top + \mathbf{e}_{\mu} \mathbf{e}_{i}^\top
    \qand
    \dimu := \underline{\Sigma}^{1/2} \Delta^{i \mu} \underline{\Sigma}^{1/2},
    \label{def-underline-delta}
\end{equation}
we can express $H$ as a sum of independent matrices
\begin{equation*}
    H 
    = \sqrt{z} \sum_{i \mu} x_{i \mu} \underline{\Sigma}^{1/2} \Delta^{i \mu} \underline{\Sigma}^{1/2}
    = \sqrt{z} \sum_{i \mu} x_{i \mu} \dimu.
\end{equation*}
Here we identify $\mathbf{e}_{i}, \mathbf{e}_{\mu}$ with their natural embedding in $\mathbb{R}^{M+N}$, as mentioned in (\ref{def-embedding}).

\begin{definition}
For any differentiable $g(H)$, we define its derivative w.r.t. $x_{i \mu}$ as
\begin{equation*}
    \partial_{i \mu} g(H)
    \equiv \frac{\partial}{\partial x_{i \mu}} g(H)
    := \left. \frac{\mathrm{d}}{\mathrm{d} t} g( H + t \sqrt{z} \dimu) \right|_{t=0}.
\end{equation*}
\end{definition}

Note that we have
\begin{equation*}
    \partial_{i \mu} H 
    = \sqrt{z} \dimu
    \qand
    \partial_{i \mu} G 
    = - G (\partial_{i \mu} H) G
    = - \sqrt{z} G \dimu G.
\end{equation*}

The main tool employed in this paper is the Green function (resolvent) 
\begin{equation*}
    G(z) := (H - z)^{-1} 
    = \begin{bmatrix}
    G_M & z^{-1/2} G_M Y \\
    z^{-1/2} Y^\top G_M & G_N
    \end{bmatrix}
    = \begin{bmatrix}
    G_M & z^{-1/2} Y G_N   \\
    z^{-1/2} G_N Y^\top & G_N
    \end{bmatrix},       
\end{equation*}
where
\begin{equation*}
    G_M (z) := (\Sigma^{1/2} XX^\top \Sigma^{1/2} - z)^{-1}
    \qand
    G_N (z) := (X^\top \Sigma X - z)^{-1}.
\end{equation*}
As in \cite{baoStatisticalInferencePrincipal2022,baoSingularVectorSingular2021,knowlesIsotropicSemicircleLaw2013}, our methodology heavily relies on the \emph{isotropic law} and \emph{average law} of the undeformed ensemble $Y$. The local law for undeformed sample covariance matrices 
$Q = YY^\top = \Sigma^{1/2} X X^\top \Sigma^{1/2}$ has been well established in \cite{baoUniversalityLargestEigenvalue2015,knowlesAnisotropicLocalLaws2017}. In essence, it states that with high probability, $G (z)$ can be well approximated by
\begin{equation}
    \Pi (z) := 
    \begin{bmatrix}
    \Pi_M(z) & 0 \\
    0 & \Pi_N(z)
    \end{bmatrix}
    = \begin{bmatrix}
    - z^{-1} (I + m(z)\Sigma)^{-1} & 0 \\
    0 & m(z) I
    \end{bmatrix},
    \label{def-Pi}
\end{equation}
where we recall that $m(z)$ represents the Stieltjes transform of the asymptotic spectral distribution of $Q$. We denote this approximation error by
\begin{equation}
    \Upsilon (z) := G (z) - \Pi (z).
    \label{def-Upsilon}
\end{equation}
For our purposes, we shall fix some small constants $\tau_1, \tau_2 > 0$ and focus on the following spectral domain outside the spectrum,
\begin{equation*}
    \mathbb{D} \equiv \mathbb{D}^{(N)} (\tau_1, \tau_2) := 
    \{ z = E + \mathrm{i} \eta \in \mathbb{C}: E \in [\lambda_{+} + \tau_1, \tau_1^{-1}],  
    \eta \in [- \tau_2, \tau_2] \}.
\end{equation*}
The projection of $\mathbb{D}$ onto the real line is denoted as 
\begin{equation*}
    \mathbb{D}_0 \equiv \mathbb{D}_0^{(N)} (\tau_1) := [\lambda_{+} + \tau_1, \tau_1^{-1}].
\end{equation*}
According to \cite{knowlesAnisotropicLocalLaws2017}, we have the following elementary estimates regarding $\varrho_N$ and its Stieltjes transform $m(z)$. In particular, the last estimate ensures the uniform boundness of $(1 + m \Sigma)^{-1}$, and thus $\Pi_M (z)$ for $z \in \mathbb{D}$.


\begin{lemma}
\label{lemma-basic-estimates-stieltjes}
Under Assumptions \ref{assumption-dimension}-\ref{assumption-regularity}, we have $\lambda_+ \asymp 1$. In addition, we have the following uniformly in $z \in \mathbb{D}$,
\begin{equation*}
    |z| \asymp 1,
    \quad 
    |m(z)| \asymp 1
    \qand
    \min_{i} | 1 + m (z) \sigma_i | \gtrsim 1.
\end{equation*}
\end{lemma}

In our proof, the following version of the Green function and its deterministic surrogate are frequently encountered,
\begin{equation*}
    \underline{G}(z) 
    := \underline{\Sigma}^{1/2} G(z) \underline{\Sigma}^{1/2}
    \qand
    \underline{\Pi}(z) 
    := \underline{\Sigma}^{1/2} \Pi(z) \underline{\Sigma}^{1/2}.
\end{equation*}
Here the upper-left block of $\underline{\Pi}$ is given by $\underline{\Pi}_M = \underline{\Sigma}^{1/2} \Pi_M \underline{\Sigma}^{1/2}$, while the bottom-right block is given by $\underline{\Pi}_N = \Pi_N$. Utilizing the self consistent equation (\ref{eqn-self-consistent}), it is not difficult to obtain
\begin{equation}
    \frac{1}{N} \operatorname{tr} \underline{\Pi}_M
    = \frac{1}{N} \operatorname{tr} \Pi_M \Sigma
    = - \frac{1 + zm}{zm}
    \qand
    \frac{1}{N} \operatorname{tr} \underline{\Pi}_N
    = \frac{1}{N} \operatorname{tr} \Pi_N
    = m.
\end{equation}

We employ the following local laws developed in \cite{knowlesAnisotropicLocalLaws2017}.

\begin{proposition}[local law]
\label{prop-local-law}
Under Assumptions \ref{assumption-bounded-moments}-\ref{assumption-regularity}, we have the following isotropic law uniformly in $z \in \mathbb{D}$ and deterministic $\bfk{u}, \bfk{v} \in \mathbb{R}^{M+N}$ with $\| \bfk{u} \|, \| \bfk{v} \| \lesssim 1$,
\begin{equation}
    \bfk{u}^\top \Upsilon(z) \bfk{v}
    = \oprec (N^{-1/2}).
    \label{eqn-isotropic-law}
\end{equation}
In addition, we have the following average law uniformly in $z \in \mathbb{D}$,
\begin{equation}
    N^{-1} \operatorname{tr} \big [ {\underline{G}_M(z) - \underline{\Pi}_M(z)} \big ], \ 
    N^{-1} \operatorname{tr} \big [ {\underline{G}_N(z) - \underline{\Pi}_N(z)} \big ] 
    = \oprec (N^{-1}).
    \label{eqn-average-law}
\end{equation}
\end{proposition}

Using Cauchy’s integral formula for derivatives, one can readily deduce from Proposition \ref{prop-local-law} the isotropic local law for $G^{\prime} (z) = \partial_z G(z)$. To be more specific, we have
\begin{equation}
    \bfk{u}^\top G^\prime(z) \bfk{v} 
    - \bfk{u}^\top \Pi^\prime(z) \bfk{v}
    = \oprec (N^{-1/2})
    \label{eqn-local-law-derivative}.
\end{equation}
Here the deterministic approximation $\Pi^\prime(z)$ is given by
\begin{equation}
    \Pi^{\prime}(z) \equiv \partial_z \Pi(z) = \begin{bmatrix}
    \partial_z \Pi_M(z) & 0 \\
    0 & \partial_z \Pi_N(z)
    \end{bmatrix}
    = \begin{bmatrix}
    z m^\prime \Pi_M \Sigma \Pi_M -z^{-1} \Pi_M & 0 \\
    0 & m^{\prime} I
    \end{bmatrix},
    \label{def-derivative-pi}
\end{equation}
where we abbreviate $m^\prime \equiv \partial_z m(z)$. On the other hand, by the definition of $H(z)$, we have
\begin{equation*}
    G^\prime(z) = - G(z) \left [ \frac{1}{2z} H(z) - I \right ] G(z)
    = - \frac{1}{2z} G(z) + \frac{1}{2} G(z)^2.
\end{equation*}
Hence, we can introduce
\begin{equation}
    \Pi_2(z) := 2 \Pi^\prime (z) + z^{-1} \Pi(z)
    = \begin{bmatrix}
    2z m^\prime \Pi_M \Sigma \Pi_M -z^{-1} \Pi_M & 0 \\
    0 & (2 m^\prime + z^{-1} m) I
    \end{bmatrix},
    \label{def-Pi-2}
\end{equation}
so that the following isotropic law for $G(z)^2$ also holds,
\begin{equation}
    \bfk{u}^\top G(z)^2 \bfk{v}
    - \bfk{u}^\top \Pi_2 (z) \bfk{v}
    = \oprec (N^{-1/2}).
\end{equation}

An important implication of the local laws is the rigidity of the eigenvalues of $Q$ around their classical locations. For our purposes, we introduce the the following rigidity estimate for the largest eigenvalue of $Q$, which allows us to control $\| H \|$ and $\| G \|$ simply by $\oprec (1)$ for $z \in \mathbb{D}$. For the proof and rigidity estimates for other eigenvalues of $Q$, we refer readers to \cite{knowlesAnisotropicLocalLaws2017} and the references therein.

\begin{lemma}[eigenvalue rigidity]
Under Assumptions \ref{assumption-bounded-moments}-\ref{assumption-regularity}, we have
\begin{equation}
   \lambda_1 (Q) - \lambda_+ = \oprec (N^{-2/3}).
   \label{eqn-eigenvalue-rigidity}
\end{equation}
\end{lemma}

We conclude this subsection by demonstrating that the spiked eigenvalues of $\tilde{Q}_{\mathrm{a}}$ are induced by those $\tilde{\sigma}_k$'s that surpass $- w_+^{-1}$, and through this analysis, we also obtain the formula of almost sure limit $\theta_k$. The following elementary lemma relates the spiked eigenvalues of $\tilde{Q}_{\mathrm{a}}$ with the singularity of a $2K \times 2K$ matrix formed by sesquilinear forms of $G(z)$. We note that this lemma has already been utilized in \cite{benaych-georgesSingularValuesVectors2012,dingHighDimensionalDeformed2020}.

\begin{lemma}
\label{lemma-determinant-equation}
For $\lambda > 0$ with $\det (Q - \lambda) \not = 0$, we have
\begin{equation*}
    \det ( \tilde{Q}_{\mathrm{a}} - \lambda ) = 0
    \quad \Longleftrightarrow \quad
    \det \big [ { \sqrt{\lambda} \mathfrak{U}^\top G(\lambda) \mathfrak{U} + \mathfrak{D}^{-1} } \big ] = 0.
\end{equation*}
\end{lemma}

As a consequence, we need to analyze the following two $2K \times 2K$ matrices,
\begin{equation}
    A_G (z) := \sqrt{z} \mathfrak{U}^\top G(z) \mathfrak{U} + \mathfrak{D}^{-1}
    \qand
    A_\Pi (z) := \sqrt{z} \mathfrak{U}^\top \Pi(z) \mathfrak{U} + \mathfrak{D}^{-1}.
    \label{def-mtr-AA-matrices}
\end{equation}
Note that $A_{\Pi} (z)$ is deterministic and according to the isotropic law (\ref{eqn-isotropic-law}), we have
\begin{equation}
    A_{\Upsilon} (z) 
    := A_G (z) - A_{\Pi} (z) = \sqrt{z} \mathfrak{U}^\top \Upsilon(z) \mathfrak{U}
    = \oprec(N^{-1/2}).
\end{equation}
In other words, one can consider solutions to $\det A_{\Pi} (\theta) = 0$ as deterministic surrogates for solutions to $\det A_{G} (\lambda) = 0$. We shall now proceed by determining the solutions to $\det A_{\Pi} (\theta) = 0$. Note that for $\theta > \lambda_+$, the matrix $I + m(\theta) \Sigma$ is invertible and thus
\begin{align*}
    \det A_{\Pi} (\theta) = 0
    & \quad \Longleftrightarrow \quad
    \det \begin{bmatrix}
    - \theta^{-1/2} U^\top (I + m(\theta) \Sigma)^{-1} U & D^{-1} \\
    D^{-1} & \theta^{1/2} m(\theta)
    \end{bmatrix} = 0 \\
    & \quad \Longleftrightarrow \quad
    \det \big [ m(\theta) U^\top (I + m(\theta) \Sigma)^{-1} U + D^{-2} \big ] = 0 \\
    & \quad \Longleftrightarrow \quad
    \det \big [ {{m(\theta)}^{-1} + \Sigma + U D^2 U^\top} \big ] = 0,
\end{align*}
where we utilized the commutativity between $D^{-1}$ and $z^{1/2} m(z)$ in the second step and employed the Weinstein-Aronszajn identity in the last step. In particular, if $\theta > \lambda_+$ is a solution to $\det A_{\Pi} (\theta) = 0$, then $-{m(\theta)}^{-1}$ should be an eigenvalue of $\tilde{\Sigma} = \Sigma + U D^2 U^\top$, the population covariance of the deformed ensemble. Recall that $m$ is negative and strictly increasing on $[\lambda_+, \infty)$, i.e.
\begin{equation*}
    - m(\theta)^{-1} > - m(\lambda_+)^{-1} = - w_+^{-1},
    \quad \forall \theta > \lambda_+.
\end{equation*}
Now it becomes evident that the spiked eigenvalues of $\tilde{Q}_{\mathrm{a}}$ are induced by those $\tilde{\sigma}_k$'s that exceed $- w_+^{-1}$. Each of these $\tilde{\sigma}_k$'s corresponds to a solution of $\det A_{\Pi} (\theta) = 0$, characterized by $m(\theta)^{-1} + \tilde{\sigma}_k = 0$. This readily leads to the definition of $\theta_k$ as given in (\ref{def-theta-k}).

\subsection{Green function representation}

The first step toward Theorem \ref{thm-spiked-distribution} involves deriving the Green function representation of the scaled fluctuations $\sqrt{N} (\lambda_k - \theta_k)$. The derivation here builds upon the equations $\operatorname{det} A_G (\lambda_k) = 0$ and $\operatorname{det} A_\Pi (\theta_k) = 0$, while incorporating our large deviation bound for $\lambda_k - \theta_k$. We highlight that the derivation in this paper requires more careful consideration compared to previous works such as \cite{knowlesIsotropicSemicircleLaw2013,zhangAsymptoticIndependenceSpiked2022}, primarily due to the absence of a block-diagonal structure in the reference matrix $A_{\Pi} (\theta_k)$. Moreover, the situation is further complicated by the non-diagonal nature of the top-left $K \times K$ block of $A_{\Pi} (\theta_k)$, which typically arises when $\Sigma \not= I$. Consequently, we cannot analyze the singularity of $A_{G} (\lambda_k)$ solely based on its diagonal entries. We overcome this difficulty by employing the Hadamard first variation formula, where the null space of $A_\Pi (\theta_k)$ plays a crucial role. The argument we present can be adapted to many other finite-rank deformation models. We also note that in the establishment of the asymptotic distribution in \cite{benaych-georgesSingularValuesVectors2012}, this procedure was bypassed as $A_G (\lambda_k)$ and $A_{\Pi} (\theta_k)$ reduce to $2 \times 2$ matrices in the case of rank-one deformation.


Let us first introduce some necessary definitions before presenting the result. Note that $A_{\Pi} (\theta_k)$ is singular as $\det A_{\Pi} (\theta_k) = 0$. Our argument below relies on the null space of $A_{\Pi} (\theta_k)$. Recall that for $k \in \llbracket K_0 \rrbracket$, we have $\tilde{\Sigma} \boldsymbol{\psi}_k = \tilde{\sigma}_k \boldsymbol{\psi}_k$ and $m(\theta_k) = -\tilde{\sigma}_k^{-1}$ by definition of $\theta_k$ in (\ref{def-theta-k}). In other words, we can write
\begin{equation*}
    \big [ {I + m(\theta_k) \tilde{\Sigma}} \big ] \boldsymbol{\psi}_k 
    = \mathbf{0}.
\end{equation*}
It follows that
\begin{equation*}
    \left [ U^\top \frac{m(\theta_k)}{I + m (\theta_k) \Sigma} U + D^{-2} \right ] (D^2 U^\top \boldsymbol{\psi}_k)
    = U^\top \big [ I + m (\theta_k) \Sigma \big ]^{-1} \big [ I + m (\theta_k) \tilde{\Sigma} \big ] \boldsymbol{\psi}_k
    = \mathbf{0}.
\end{equation*}
By definition of $A_\Pi (\theta_k)$, we therefore deduce that
\begin{equation}
    A_{\Pi} (\theta_k) \boldsymbol{\xi}_k = \mathbf{0}
    \quad \text{ where } \quad
    \boldsymbol{\xi}_k 
    := \begin{bmatrix} 
    - \sqrt{\theta_k} m(\theta_k) D^2 U^\top \boldsymbol{\psi}_k \\ 
    D U^\top \boldsymbol{\psi}_k 
    \end{bmatrix}.
    \label{def-eigenvector-A-Pi}
\end{equation}
Furthermore, due to Assumption \ref{assumption-spiked-spacing} on eigenvalue spacing, it can be readily verified that $0$ is a simple eigenvalue of $I + m(\theta_k) \tilde{\Sigma}$. Consequently, the null space of $A_{\Pi} (\theta_k)$ is spanned by $\boldsymbol{\xi}_k$. Next, for each $k \in \llbracket K_0 \rrbracket$, let us introduce the vectors $\mathbf{u}_k \in \mathbb{R}^M$ and $\mathbf{v}_k \in \mathbb{R}^N$ such that
\begin{equation*}
    \begin{bmatrix} \mathbf{u}_k \\ \mathbf{v}_k \end{bmatrix}
    := \mathfrak{U} \boldsymbol{\xi}_k
    = \begin{bmatrix} 
    - \sqrt{\theta_k} m(\theta_k) U D^2 U^\top \boldsymbol{\psi}_k \\ 
    V D U^\top \boldsymbol{\psi}_k 
    \end{bmatrix}.
\end{equation*}
Note that using $\big [ {I + m(\theta_k) \tilde{\Sigma}} \big ] \boldsymbol{\psi}_k = \mathbf{0}$ and $S = U D V^\top$, one may also write
\begin{equation}
    \mathbf{u}_k 
    = \sqrt{\theta_k} \big [ {I + m(\theta_k) \Sigma} \big ] \boldsymbol{\psi}_k
    \qand
    \mathbf{v}_k
    = S^\top \boldsymbol{\psi}_k .
    \label{def-uv-vectors}
\end{equation}
Importantly, unless $\Sigma = I$, the vectors $\mathbf{u}_k$ and $\mathbf{v}_k$ in general cannot be solely determined by $\bfk{u}_k$ and $\bfk{v}_k$, the singular vectors associated with the $k$-th largest singular value of the signal. Consequently, these vectors, and therefore the fluctuation of $\lambda_k$, are generally influenced by the other spectral components of $S$. Now, we are in a position to present the Green function representation of the spiked eigenvalues.

\begin{proposition}[Green function representation] 
\label{prop-green-function-representation}
For each $k \in \llbracket K_0 \rrbracket$, we have
\begin{equation}
    \sqrt{N} ( \sqrt{\lambda_k} - \sqrt{\theta_k} )
    = -\sqrt{N} \cdot \frac{\boldsymbol{\xi}_k^\top A_{\Upsilon} (\theta_k) \boldsymbol{\xi}_k}
    {\boldsymbol{\xi}_k^\top B_\Pi (\theta_k) \boldsymbol{\xi}_k} 
    + \oprec (N^{-1/2}),
    \label{eqn-green-function-representation}
\end{equation}
where $B_\Pi (z) := z \mathfrak{U}^\top \Pi_2 (z) \mathfrak{U}$ and $\Pi_2$ is defined in (\ref{def-Pi-2}). Consequently,
\begin{equation}
    \sqrt{N} ( \lambda_k - \theta_k )
    = -\sqrt{N} \tilde{\sigma}_k \theta^\prime_k
    \begin{bmatrix} \mathbf{u}_k \\ \mathbf{v}_k \end{bmatrix}^\top 
    \Upsilon (\theta_k)
    \begin{bmatrix} \mathbf{u}_k \\ \mathbf{v}_k \end{bmatrix}
    + \oprec (N^{-1/2}).
    \label{eqn-simplified-representation}
\end{equation}
\end{proposition}

In order to present the proof of Proposition \ref{prop-green-function-representation}, we need to revisit some elementary results regarding the first-order perturbation of determinants. Fix some $n > 0$. Let $A$ and $F$ be $n \times n$ symmetric matrices with columns
\begin{equation*}
    A = [\mathbf{a}_1, \cdots, \mathbf{a}_n]
    \quad \text{ and } \quad
    F = [\mathbf{f}_1, \cdots, \mathbf{f}_n].
\end{equation*}
As the determinant is a multilinear function of the matrix columns, we have for any $t \in \mathbb{R}$,
\begin{equation}
    \det (A + tF) 
    = \sum_{\mathcal{J} \subset [n]} t^{|\mathcal{J}|} \det A^{(\mathcal{J})}
    = \sum_{p = 0}^n t^{p} \sum_{| \mathcal{J} | = p} \det A^{(\mathcal{J})},
    \label{eqn-determinant-expansion}
\end{equation}
where the the $n \times n$ matrix $A^{(\mathcal{J})}$ has columns $\mathbf{a}_{j}^{(\mathcal{J})} = \mathbbm{1} (j \notin \mathcal{J}) \mathbf{a}_{j} + \mathbbm{1} (j \in \mathcal{J}) \mathbf{f}_{j}$. In particular, for each $0 \leq p \leq n$, we have
\begin{equation*}
    \sum_{| \mathcal{J} | = p} \det A^{(\mathcal{J})} 
    = \frac{1}{p!} \big [ \partial_t^p \det (A + tF) \big ]_{t = 0}.
\end{equation*}
The following expression of $\big [ \partial_t \det (A + tF) \big ]_{t = 0}$ for singular $A$ played a crucial role in our proof of Proposition \ref{prop-green-function-representation}.

\begin{lemma}
\label{lemma-perturbation-singular}
Let $A$ and $F$ be $n \times n$ symmetric matrices such that $0$ is a simple eigenvalue of $A$. Denote the eigenvalues of $A$ by $s_1, \cdots, s_n$ and suppose $\boldsymbol{\zeta}$ is the unit eigenvector of $A$ associated with $0$. Then, we have
\begin{equation*}
    \big [ \partial_t \det (A + tF) \big ]_{t = 0}
    = \boldsymbol{\zeta}^\top F \boldsymbol{\zeta} \cdot \prod_{k: s_k \not= 0} s_k.
\end{equation*}
\end{lemma}

\begin{proof}[Proof of Lemma \ref{lemma-perturbation-singular}]
Denote the spectral decomposition of $(A + tF)$ by
\begin{equation*}
    A + tF = \sum_{j = 1}^n s_j (t) \boldsymbol{\zeta}_j (t) \boldsymbol{\zeta}_j (t)^\top,
\end{equation*}
where $s_j (0) = s_j$. Let $k_0 \in [n]$ be the index with $s_{k_0} (0) = 0$ and $\boldsymbol{\zeta}_{k_0} (0) = \boldsymbol{\zeta}$. Using the Hadamard first variation formula (see, e.g. \cite[Equation (1.73)]{tao2012topics}), we have for each $j \in [n]$,
\begin{equation*}
    \partial_t s_j (t) = \boldsymbol{\zeta}_j (t)^\top F \boldsymbol{\zeta}_j (t).
\end{equation*}
It follows that,
\begin{equation*}
    \partial_t \det (A + tF) 
    = \partial_t \left [ \prod_{j = 1}^n s_j (t) \right ] 
    = \sum_{j = 1}^n \left [ \boldsymbol{\zeta}_j (t)^\top F \boldsymbol{\zeta}_j (t) 
    \cdot \prod_{k \not= j} s_k (t) \right ]. 
\end{equation*}
Now, setting $t = 0$, we get
\begin{equation*}
    \big [ \partial_t \det (A + tF) \big ]_{t = 0} 
    = \boldsymbol{\zeta}_{k_0} (0)^\top F \boldsymbol{\zeta}_{k_0} (0) 
    \cdot \prod_{k \not= k_0} s_k (0)
    = \boldsymbol{\zeta}^\top F \boldsymbol{\zeta} \cdot \prod_{k \not= k_0} s_k,
\end{equation*}
which concludes the proof.
\end{proof}

\begin{proof}[Proof of Proposition \ref{prop-green-function-representation}]
For simplicity, we omit the subscript $k$ as it applies to all $k \in \llbracket K_0 \rrbracket$. Firstly, we have the following resolvent expansion formula, which follows directly from the definition of the resolvent,
\begin{equation}
    \sqrt{\lambda} G(\lambda) - \sqrt{\theta} G(\theta) 
    = \sqrt{\lambda \theta} (\sqrt{\lambda} - \sqrt{\theta}) G(\lambda) G(\theta).
    \label{eqn-resolvent-expansion}
\end{equation}
The following two matrices will arise as a consequence of the resolvent expansion,
\begin{equation}
    B_G (z) := z \mathfrak{U}^\top G(z)^{2} \mathfrak{U}
    \qand
    B_{\Upsilon} (z)
    = B_G (z) - B_{\Pi} (z)
    = z \mathfrak{U}^\top \Upsilon_2 (z) \mathfrak{U}.      
    \label{def-mtr-B-matrices}
\end{equation}
Recalling (\ref{def-mtr-AA-matrices}) and applying (\ref{eqn-resolvent-expansion}) twice, we have
\begin{align*}
    A_{G} (\lambda)
    & = A_{G} (\theta) + \sqrt{\lambda \theta} (\sqrt{\lambda} - \sqrt{\theta}) \mathfrak{U}^\top G(\lambda) G(\theta) \mathfrak{U}^\top \\
    & = A_{G} (\theta) + (\sqrt{\lambda} - \sqrt{\theta}) B_{G} (\theta) + \sqrt{\lambda} \theta (\sqrt{\lambda} - \sqrt{\theta})^2 \mathfrak{U}^\top G(\lambda) G(\theta)^2 \mathfrak{U} \\
    & = A_{\Pi} (\theta) + A_{\Upsilon} (\theta) + (\sqrt{\lambda} - \sqrt{\theta}) B_{\Pi} (\theta) + \oprec (N^{-1}),
\end{align*}
where in the last step we absorbed the third term of the second line and $(\sqrt{\lambda} - \sqrt{\theta}) B_{\Upsilon} (z)$ into $\oprec (N^{-1})$ using the isotropic law (\ref{eqn-isotropic-law}) and the preliminary estimate $\sqrt{\lambda} - \sqrt{\theta} = \oprec (N^{-1/2})$ from Proposition \ref{prop-convergence-rate}. Recall $\det A_{G} (\lambda) = 0$. Hence, by (\ref{eqn-determinant-expansion}) we can get
\begin{equation*}
    \det \big [ { A_{\Pi} (\theta) + A_{\Upsilon} (\theta) + (\sqrt{\lambda} - \sqrt{\theta}) B_{\Pi} (\theta) } \big ] 
    = \oprec (N^{-1}).
\end{equation*}
Now we apply the determinant expansion formula (\ref{eqn-determinant-expansion}) to l.h.s. of the above equation. By noting $A_{\Upsilon} (\theta)$ and $\sqrt{\lambda} - \sqrt{\theta} = \oprec (N^{-1/2})$ again, we get
\begin{align*}
    & \ \det \big [ { A_{\Pi} (\theta) + A_{\Upsilon} (\theta) + (\sqrt{\lambda} - \sqrt{\theta}) B_{\Pi} (\theta) } \big ] \\
    = & \ \partial_t \big [ { \det ( A_{\Pi} (\theta) + t A_{\Upsilon} (\theta) ) } \big ]_{t = 0} 
    + (\sqrt{\lambda} - \sqrt{\theta}) \partial_t \big [ { \det ( A_{\Pi} (\theta) + t B_{\Pi} (\theta) ) } \big ]_{t = 0} 
    + \oprec (N^{-1}),
\end{align*}
where we also used $\det A_{\Pi} (\theta) = 0$. Recalling (\ref{def-eigenvector-A-Pi}), the unit eigenvector of $A_{\Pi} (\theta)$ associated with $0$ is given by $\boldsymbol{\xi}$. Hence, by applying Lemma \ref{lemma-perturbation-singular} we can obtain,
\begin{equation*}
    \sqrt{\lambda} - \sqrt{\theta}
    = -\frac{\partial_t \big [ { \det ( A_{\Pi} (\theta) + t A_{\Upsilon} (\theta) ) } \big ]_{t = 0}}
    {\partial_t \big [ { \det ( A_{\Pi} (\theta) + t B_{\Pi} (\theta) ) } \big ]_{t = 0}} 
    + \oprec (N^{-1})
    = -\frac{\boldsymbol{\xi}^\top A_{\Upsilon} (\theta) \boldsymbol{\xi}}{\boldsymbol{\xi}^\top B_{\Pi} (\theta) \boldsymbol{\xi}} + \oprec (N^{-1}).
\end{equation*}
Multiplying both sides with $\sqrt{N}$ concludes the proof of (\ref{eqn-green-function-representation}). As for (\ref{eqn-simplified-representation}), let us recall the definition of the vectors $\mathbf{u}_k$ and $\mathbf{v}_k$ in (\ref{def-uv-vectors}) and $\Pi_2$ in (\ref{def-Pi-2}), which yields
\begin{align*}
    \boldsymbol{\xi}^\top B_\Pi (\theta) \boldsymbol{\xi} 
    & = 2\theta^2 m^\prime(\theta) \mathbf{u}^\top \Pi_M(\theta) \Sigma \Pi_M(\theta) \mathbf{u}
    - \mathbf{u}^\top \Pi_M(\theta) \mathbf{u}
    + \big [ {2 \theta m^\prime(\theta) + m(\theta)} \big ] \mathbf{v}^\top \mathbf{v} \\
    & = 2 \theta m^\prime(\theta) \boldsymbol{\psi}^\top \Sigma \boldsymbol{\psi}
    + 2 \theta m^\prime(\theta) \boldsymbol{\psi}^\top S S^\top \boldsymbol{\psi} \\
    & = 2 \theta m^\prime(\theta) \boldsymbol{\psi}^\top \tilde{\Sigma} \boldsymbol{\psi}
    = \frac{2 \theta }{\tilde{\sigma} \theta^\prime},
\end{align*}
where in the second step we used
\begin{equation*}
    \mathbf{u}^\top \Pi_M(\theta) \mathbf{u} 
    = - \boldsymbol{\psi}^\top \big [ {I + m(\theta) \Sigma} \big ] \boldsymbol{\psi}
    = m(\theta) \boldsymbol{\psi}^\top S S^\top \boldsymbol{\psi}
    = m(\theta) \mathbf{v}^\top \mathbf{v} .
\end{equation*}
On the other hand, by the large deviation bound on $\sqrt{\lambda} - \sqrt{\theta}$, we have
\begin{align*}
    \sqrt{N} ( \lambda - \theta )
    = 2 \sqrt{N \theta} ( \sqrt{\lambda} - \sqrt{\theta} ) + \oprec (N^{-1/2}).
\end{align*}
Summarizing these two results concludes the proof of (\ref{eqn-simplified-representation}).
\end{proof}

\subsection{Joint distribution of sesquilinear forms}
\label{subsec-sesquilinear-forms}

With the Green function representation in place, our subsequent step towards proving Theorem \ref{thm-spiked-distribution} entails developing the asymptotic joint distribution of centralized sesquilinear forms of the resolvent. This constitutes the primary technical achievement of this paper, and we formally state it in Theorem \ref{thm-characteristic-function}. It is worth noting that the analysis of sesquilinear forms of the resolvents has been a subject of extensive research across various contexts. In our specific notation, \cite{zhangAsymptoticIndependenceSpiked2022,baoStatisticalInferencePrincipal2022} examined the distribution of $\mathbf{u}^\top \Upsilon (z) \mathbf{u}$ where $\mathbf{u} \in \mathbb{R}^M$. In these contexts, the resulting asymptotic distributions do not manifest any nonuniversality. On the other hand, the nonuniversality of sesquilinear forms was discussed in the works \cite{knowlesIsotropicSemicircleLaw2013,baoSingularVectorSingular2021}, where the underlying deterministic surrogate $\Pi$ is isotropic. An evident advantage of our approach in characterizing nonuniversality, compared to these works, is the circumvention of an additional decomposition of $\Theta$ into a Gaussian and a non-Gaussian part. This decomposition, based on the magnitudes of the coefficients in front of the $x_{i \mu}$'s within $\Theta$, would significantly aggravate the computation in our case due to the presence of the anisotropic $\Sigma$.


Let $\{ z_k \}_{k \in \mathcal{J}} \subset \mathbb{D}_0$ be a finite set of spectral parameters (not necessarily distinct). Assume that we are given the deterministic vectors $\mathbf{u}_k, \mathbf{r}_k \in \mathbb{R}^M$, $\mathbf{v}_k, \mathbf{w}_k \in \mathbb{R}^N$, as well as the coefficients $\alpha_k, \beta_k, \gamma_k, \delta_k \in \mathbb{R}$, such that
\begin{equation}
    \| \mathbf{u}_k \|, \| \mathbf{r}_k \|, \| \mathbf{v}_k \|, \| \mathbf{w}_k \| \asymp 1
    \qand
    \max \big \{ |\alpha_k|, |\beta_k|, |\gamma_k|, |\delta_k| \big \} \asymp 1.
    \label{assumption-coefficient-size}
\end{equation}
We are interested in the asymptotic joint distribution of the sesquilinear forms
\begin{equation}
    \mathcal{Q}_k 
    = \sqrt{N} \big [ {\alpha_k \mathbf{u}_k^\top \Upsilon (z_k) \mathbf{r}_k
    + \beta_k \mathbf{v}_k^\top \Upsilon (z_k) \mathbf{w}_k
    + \gamma_k \mathbf{u}_k^\top \Upsilon (z_k) \mathbf{w}_k
    + \delta_k \mathbf{v}_k^\top \Upsilon (z_k) \mathbf{r}_k} \big ].
    \label{def-quantites-Q}
\end{equation}
Especially, when $\alpha_k, \beta_k, \gamma_k, \delta_k = 1$, we can write
\begin{equation*}
    \mathcal{Q}_k = \sqrt{N}
    \begin{bmatrix} \mathbf{u}_k \\ \mathbf{v}_k \end{bmatrix}^\top 
    \Upsilon (z_k)
    \begin{bmatrix} \mathbf{r}_k \\ \mathbf{w}_k \end{bmatrix}.
\end{equation*}
Although the derivation of Theorem \ref{thm-spiked-distribution} only involves sesquilinear forms with $\mathbf{u}_k = \mathbf{r}_k$ and $\mathbf{v}_k = \mathbf{w}_k$, as demonstrated by the Green function representation (\ref{eqn-simplified-representation}), we choose to allow flexibility in the choice of vectors to maintain generality and facilitate future investigations. 

Before describing the asymptotic distribution of the $\mathcal{Q}_k$'s, let us introduce some additional notations to simplify the presentation. We will use the abbreviation $m_k \equiv m(z_k)$ and denote the diagonals of $\underline{\Pi}_{M} (z_k)$ and $\underline{\Pi}_{N} (z_k)$ as $\boldsymbol{\pi}_k$ and $\mathbf{m}_k$, respectively. In other words,
\begin{equation}
    (\boldsymbol{\pi}_k)_i = \underline{\Pi}_{i i} (z_k) = \left [-{\frac{\Sigma}{z_k (I + m_k \Sigma)}} \right ]_{i i}
    \qand
    (\mathbf{m}_k)_\mu = \underline{\Pi}_{\mu \mu} (z_k) = m_k.
    \label{def-pi-vec}
\end{equation}
Given a vector in $\mathbb{R}^M$ or $\mathbb{R}^N$ and indexed by $k$, denote the result of multiplying it by the matrix $\underline{\Sigma}^{1/2} \Pi (z_k)$ by placing an underline. For instance, 
\begin{equation}
    \underline{\mathbf{u}}_k 
    = \Sigma^{1/2} \Pi_M (z_k) \mathbf{u}_k
    = - \frac{\Sigma^{1/2}}{z_k (1 + m_k \Sigma)} \mathbf{u}_k
    \qand
    \underline{\mathbf{v}}_k 
    = \Pi_N (z_k) \mathbf{v}_k
    \equiv m_k \mathbf{v}_k.
    \label{def-underline-vec}
\end{equation}
We also adopt the notation $m [z_k, z_j]$ for divided differences of the Stieltjes transform. That is, $m [z_k, z_j] = m^\prime (z_k)$ if $z_k = z_j$. Otherwise,
\begin{equation*}
    m [z_k, z_j] := \frac{m(z_k) - m(z_j)}{z_k - z_j}.
\end{equation*}
Similar to Definition \ref{def-spiked-quantities}, in the following definition we summarize the quantities that are necessary for the asymptotic distribution of the $\mathcal{Q}_k$'s. Recall the notion of mixed moments that we introduced before Definition \ref{def-spiked-quantities}. Also note that by Assumption \ref{assumption-coefficient-size} and Lemma \ref{lemma-basic-estimates-stieltjes}, the vectors $\underline{\mathbf{u}}_k, \underline{\mathbf{v}}_k, \underline{\mathbf{r}}_k, \underline{\mathbf{w}}_k$, and thus the quantities defined below, remain bounded as $N$ grows.

\begin{definition}[asymptotic quantities for Theorem \ref{thm-characteristic-function}]
\label{def-sesquilinear-quantities}
For each $k \in \mathcal{J}$, let
\begin{equation*}
    \mathcal{L}_k := -\frac{\kappa_3 z_k \sqrt{z_k}}{N}
    \Big [ \gamma_k \mathbb{M} (\boldsymbol{\pi}_k, \underline{\mathbf{u}}_k )
    \mathbb{M} (\mathbf{m}_k, \underline{\mathbf{w}}_k ) +
    \delta_k \mathbb{M} (\mathbf{m}_k, \underline{\mathbf{v}}_k ) 
    \mathbb{M} (\boldsymbol{\pi}_k, \underline{\mathbf{r}}_k ) \Big ].
\end{equation*}
For each $k, j \in \mathcal{J}$, we introduce $\mathcal{V}_{kj} := \mathcal{V}_{kj}^{(0, 1, 0)} + \mathcal{V}_{kj}^{(1, 2, 0)}$. Here
\begin{align*}
    \mathcal{V}_{kj}^{(0, 1, 0)} :=
    z_k z_j m [z_k, z_j]
    & \Big [ {\alpha_k \alpha_j \mathbb{M} (\underline{\mathbf{u}}_k, \underline{\mathbf{u}}_j )
    \mathbb{M} (\underline{\mathbf{r}}_k, \underline{\mathbf{r}}_j ) 
    + \alpha_k \alpha_j \mathbb{M} (\underline{\mathbf{u}}_k, \underline{\mathbf{r}}_j )
    \mathbb{M} (\underline{\mathbf{r}}_k, \underline{\mathbf{u}}_j )} \Big ] \\
    + \frac{m [z_k, z_j] - m_k m_j}{m_k^2 m_j^2} 
    & \Big [ {\beta_k \beta_j \mathbb{M} (\underline{\mathbf{v}}_k, \underline{\mathbf{v}}_j )
    \mathbb{M} (\underline{\mathbf{w}}_k, \underline{\mathbf{w}}_j ) 
    + \beta_k \beta_j \mathbb{M} (\underline{\mathbf{v}}_k, \underline{\mathbf{w}}_j )
    \mathbb{M} (\underline{\mathbf{w}}_k, \underline{\mathbf{v}}_j )} \Big ] \\
    + \frac{\sqrt{z_k z_j} (m [z_k, z_j] - m_k m_j)}{m_k m_j}
    & \Big [ \gamma_k \gamma_j \mathbb{M} (\underline{\mathbf{u}}_k, \underline{\mathbf{u}}_j )
    \mathbb{M} (\underline{\mathbf{w}}_k, \underline{\mathbf{w}}_j ) 
    + \gamma_k \delta_j \mathbb{M} (\underline{\mathbf{u}}_k, \underline{\mathbf{r}}_j )
    \mathbb{M} (\underline{\mathbf{w}}_k, \underline{\mathbf{v}}_j ) \\
    & + \delta_k \gamma_j \mathbb{M} (\underline{\mathbf{v}}_k, \underline{\mathbf{w}}_j )
    \mathbb{M} (\underline{\mathbf{r}}_k, \underline{\mathbf{u}}_j )
    + \delta_k \delta_j \mathbb{M} (\underline{\mathbf{v}}_k, \underline{\mathbf{v}}_j )
    \mathbb{M} (\underline{\mathbf{r}}_k, \underline{\mathbf{r}}_j ) \Big ],
\end{align*}
\begin{equation*}
    \mathcal{V}_{kj}^{(1, 2, 0)} := \frac{\kappa_{4} z_k z_j}{N} 
    \Big [ {\alpha_k \alpha_j 
    \mathbb{M} (\underline{\mathbf{u}}_k, \underline{\mathbf{r}}_k, \underline{\mathbf{u}}_j, \underline{\mathbf{r}}_j )
    \mathbb{M} (\mathbf{m}_k, \mathbf{m}_j )
    + \beta_k \beta_j 
    \mathbb{M} (\underline{\mathbf{v}}_k, \underline{\mathbf{w}}_k, \underline{\mathbf{v}}_j, \underline{\mathbf{w}}_j ) 
    \mathbb{M} (\boldsymbol{\pi}_k, \boldsymbol{\pi}_j )} \Big ] .
\end{equation*}
Finally, for each $k, j \in \mathcal{J}$, we let
\begin{align*}
    \mathcal{W}_{kj} 
    := - \frac{\kappa_{3} z_k \sqrt{z_j}}{N^{1/2}} 
    & \Big [ \alpha_k \gamma_j \mathbb{M} (\underline{\mathbf{u}}_k, \underline{\mathbf{r}}_k, \underline{\mathbf{u}}_j )
    \mathbb{M} (\mathbf{m}_k, \underline{\mathbf{w}}_j )
    + \alpha_k \delta_j \mathbb{M} (\underline{\mathbf{u}}_k, \underline{\mathbf{r}}_k, \underline{\mathbf{r}}_j )
    \mathbb{M} (\mathbf{m}_k, \underline{\mathbf{v}}_j ) \\
    & + \beta_k \gamma_j \mathbb{M} (\underline{\mathbf{v}}_k, \underline{\mathbf{w}}_k, \underline{\mathbf{w}}_j )
    \mathbb{M} (\boldsymbol{\pi}_k, \underline{\mathbf{u}}_j )
    + \beta_k \delta_j \mathbb{M} (\underline{\mathbf{v}}_k, \underline{\mathbf{w}}_k, \underline{\mathbf{v}}_j )
    \mathbb{M} (\boldsymbol{\pi}_k, \underline{\mathbf{r}}_j ) \Big ].
\end{align*}
\end{definition}

\begin{theorem}[distribution of sesquilinear forms]
\label{thm-characteristic-function}
Consider the sesquilinear forms $\{ \mathcal{Q}_k \}_{k \in \mathcal{J}}$ defined in (\ref{def-quantites-Q}). Let $\mathcal{L}_{k}, \mathcal{V}_{kj}, \mathcal{W}_{kj}$ be given as in Definition \ref{def-sesquilinear-quantities}. For each $k \in \mathcal{J}$, set
\begin{align*}
    \Theta_k & := - \sqrt{N} \big [ {
    \gamma_k \mathbf{u}_k^\top \Pi (z_k) H (z_k) \Pi (z_k) \mathbf{w}_k 
    + \delta_k \mathbf{v}_k^\top \Pi (z_k) H (z_k) \Pi (z_k) \mathbf{r}_k} \big ], \\
    \Phi_k & := \mathcal{Q}_k - \Theta_k- \mathcal{L}_k.
\end{align*}
Let $\Phi := \sum_{k} s_k \Phi_k$ and $\Theta := \sum_{k} t_k \Theta_k$ where $|s_k|, |t_k| \lesssim 1$ are deterministic. Then, given any $\varepsilon > 0$, there exists $N_0 \equiv N_0 (\varepsilon)$ such that for all $N \geq N_0$,
\begin{equation}
    \bigg | {\mathbb{E} \exp [{\mathrm{i} (\Phi + \Theta)}]
    - \exp \left ( - \frac{\mathcal{V} + 2 \mathcal{W}}{2} \right )
    \mathbb{E} \exp ({\mathrm{i} \Theta})} \bigg | \leq \epsilon,
    \label{eqn-characteristic-function}
\end{equation}
where $\mathcal{V} = \sum_{k,j} s_k s_j \mathcal{V}_{kj}$ and $\mathcal{W} = \sum_{k,j} s_k t_j \mathcal{W}_{kj}$.
\end{theorem}

Let us provide some comments on Theorem \ref{thm-characteristic-function}. For any $k, j$, the two sesquilinear forms 
\begin{equation*}
    \mathbf{u}_k^\top \Upsilon (z_k) \mathbf{r}_k
    = \begin{bmatrix} \mathbf{u}_k \\ \mathbf{0} \end{bmatrix}^\top 
    \Upsilon (z_k)
    \begin{bmatrix} \mathbf{r}_k \\ \mathbf{0} \end{bmatrix}
    \qand
    \mathbf{v}_j^\top \Upsilon (z_j) \mathbf{w}_j
    = \begin{bmatrix} \mathbf{0} \\ \mathbf{v}_j \end{bmatrix}^\top 
    \Upsilon (z_j)
    \begin{bmatrix} \mathbf{0} \\ \mathbf{w}_j \end{bmatrix}
\end{equation*}
do not contribute to the nonuniversality of $\mathcal{Q}_k$. Instead, they exhibit an asymptotic joint Gaussian fluctuation that does not depend on the cumulants $\kappa_p$ with $p > 4$. In addition, they are asymptotically independent, as the coefficient $\alpha_k \beta_j$ is absent in definition of $\mathcal{V}_{kj}$. We highlight that for the spiked population model, the Green function representation for the fluctuation of spiked eigenvalues only involves sesquilinear forms of the type $\mathbf{u}^\top \Upsilon (z) \mathbf{u}$, where $\mathbf{u} \in \mathbb{R}^M$ (see e.g. \cite{zhangAsymptoticIndependenceSpiked2022}). In particular, this elucidates the absence of nonuniversality in the context of spiked covariance matrices.


Here, we would also like to briefly explain how we arrive at the specific form of $\Theta_k$ in Theorem \ref{thm-characteristic-function}. Let us temporarily set aside the subscript $k$. Recall the derivative formulas $\partial_{i \mu} H (z) = \sqrt{z} \dimu$ and $\partial_{i \mu} G (z) = - \sqrt{z} G(z) \dimu G(z)$, Therefore,
\begin{equation*}
    \partial_{i \mu} \mathcal{Q}
    = - \sqrt{N z} \big [
    \alpha \mathbf{u}^\top G \dimu G \mathbf{r}
    + \beta \mathbf{v}^\top G \dimu G \mathbf{w}
    + \gamma \mathbf{u}^\top G \dimu G \mathbf{w}
    + \delta \mathbf{v}^\top G \dimu G \mathbf{r} \big ]    .
\end{equation*}
With the isotropic local law (\ref{eqn-isotropic-law}), it is natural to consider substituting $G$ with $\Pi$ in the above expression. The first two terms after substitution, $\mathbf{u}^\top \Pi \dimu \Pi \mathbf{r}$ and $\mathbf{v}^\top \Pi \dimu \Pi \mathbf{w}$, vanish due to the block anti-diagonal structure of $\dimu$. As a result, this substitution yields 
\begin{equation*}
    - \sqrt{N z} \big [
    \gamma \mathbf{u}^\top \Pi \dimu \Pi \mathbf{w}
    + \delta \mathbf{v}^\top \Pi \dimu \Pi \mathbf{r} \big ],
\end{equation*}
which precisely captures the derivative of $\Theta$ w.r.t. $x_{i \mu}$. That is, the definition of $\Theta$ enables effective control over $\partial_{i \mu} \Phi = \partial_{i \mu} \mathcal{Q} - \partial_{i \mu} \Theta$, which plays a crucial role in derivation of the asymptotic Gaussianity of $\Phi$. This heuristic may be valuable in extracting the nonuniversal component of other functionals of resolvents, especially in models where the deterministic surrogate $\Pi (z)$ of the resolvent is anisotropic.

Now, the asymptotic joint distribution of the spiked eigenvalues can be obtained via a straightforward combination of Proposition \ref{prop-green-function-representation} and Theorem \ref{thm-characteristic-function}.

\begin{proof}[Proof of Theorem \ref{thm-spiked-distribution}]
According to Proposition \ref{prop-green-function-representation}, we need to specify $z_k = \theta_k$, and letting $\mathbf{u}_k = \mathbf{r}_k$ and $\mathbf{v}_k = \mathbf{w}_k$ be as given in (\ref{def-uv-vectors}). Recall (\ref{def-underline-vec}). In this case, we have
\begin{equation*}
    \underline{\mathbf{u}}_k = \underline{\mathbf{r}}_k
    = - \Sigma^{1/2} \boldsymbol{\psi}_k / \sqrt{\theta_k}
    \qand
    \underline{\mathbf{v}}_k = \underline{\mathbf{w}}_k
    = - S^\top \boldsymbol{\psi}_k / \tilde{\sigma}_k.
\end{equation*}
In addition, by $m(\theta_k) = - \tilde{\sigma}_k^{-1}$ along with the definitions of $\boldsymbol{\pi}_k$ and $\tilde{\boldsymbol{\pi}}_k$ provided in (\ref{def-pi-vec}) and (\ref{def-tilde-pi-vec}), we have
\begin{equation*}
    \boldsymbol{\pi}_k= -\tilde{\sigma}_k \tilde{\boldsymbol{\pi}}_k/\theta_k
    \qand
    \mathbf{m}_k = -\mathbf{1}_N/\tilde{\sigma}_k.
\end{equation*}
We also need to specify $\alpha_k, \beta_k, \gamma_k, \delta_k = 1$. Now, the expressions for $\mathcal{L}_{k}, \mathcal{V}_{kj}^{(1,2,0)}, \mathcal{W}_{kj}$ in Definition \ref{def-spiked-quantities} can be obtained by combining Definition \ref{def-sesquilinear-quantities} with the prefactor $\tilde{\sigma}_k \theta^\prime_k$ in (\ref{eqn-simplified-representation}). As for $\mathcal{V}_{kj}^{(0,1,0)}$, we note that
\begin{equation*}
    m [\theta_k, \theta_k] = m^\prime (\theta_k) = \frac{1}{\tilde{\sigma}_k^2 \theta_k^\prime}
    \qand
    \boldsymbol{\psi}_k^\top S S^\top \boldsymbol{\psi}_j
    = \tilde{\sigma}_k \mathbbm{1}_{k = j} - \boldsymbol{\psi}_k^\top \Sigma \boldsymbol{\psi}_j,
\end{equation*}
where the second one follows from $\tilde{\sigma}_j \boldsymbol{\psi}_j = (\Sigma + SS^\top) \boldsymbol{\psi}_j$. Plugging these expressions into Definition \ref{def-sesquilinear-quantities} yields the desired result via Theorem \ref{thm-characteristic-function}.
\end{proof}

The proof of Theorem \ref{thm-characteristic-function} is provided in detail in the Supplementary Material. The main idea is to establish the recursive estimate (\ref{eqn-approx-recurrence-moment}), which is derived by incorporating the cumulant expansion technique and the local laws for the undeformed ensemble. Let us conclude this section by highlighting how we establish the convergence of the characteristic function (\ref{eqn-spiked-distribution}) or (\ref{eqn-characteristic-function}) through the recurrence relation of moments (\ref{eqn-approx-recurrence-moment}). While the overall idea is straightforward, we are not aware of its previous appearance in the field of RMT. Recall that the Hermite polynomials can be defined as
\begin{equation*}
    \mathrm{He}_{\ell} (\xi) 
    = (-1)^\ell \mathrm{e}^{\xi^2/2} \frac{\mathrm{d}^\ell (\mathrm{e}^{-\xi^2/2})}{\mathrm{d} \xi^\ell},
    \quad \ell \geq 0.
\end{equation*}
In addition, for each $\ell \geq 0$, we have the recurrence relation
\begin{equation*}
    \mathrm{He}_{\ell + 2} (\xi) = -(\ell + 1) \mathrm{He}_{\ell} (\xi) + \xi \mathrm{He}_{\ell + 1} (\xi).
\end{equation*}
Hence, once we establish the recurrence relation (\ref{eqn-approx-recurrence-moment}), we can proceed with a straightforward induction to derive 
\begin{equation*}
    \mathbb{E} [ {\Phi^{\ell} \mathrm{e}^{\mathrm{i} \Theta}} ]
    \approx (\mathrm{i} \sqrt{\mathcal{V}})^{\ell} \mathrm{He}_{\ell} (\mathcal{W} / \sqrt{\mathcal{V}}) 
    \mathbb{E} \mathrm{e}^{\mathrm{i} \Theta}.
\end{equation*}
Now, by utilizing Taylor's expansion, we are able to deduce
\begin{align*}
    \mathbb{E} \exp [ \mathrm{i} (\Phi + \Theta) ]
    & = \sum_{\ell = 0}^{\infty} \frac{1}{\ell !} \mathbb{E} [ {(\mathrm{i} \Phi)^{\ell} \mathrm{e}^{\mathrm{i} \Theta}} ]
    \approx \sum_{\ell = 0}^{\infty} \frac{(- \sqrt{\mathcal{V}})^\ell}{\ell !} \mathrm{He}_{\ell} (\mathcal{W} / \sqrt{\mathcal{V}}) 
    \mathbb{E} \mathrm{e}^{\mathrm{i} \Theta} \\
    & = \sum_{\ell = 0}^{\infty} \frac{\mathcal{V}^{\ell/2}}{\ell !}  
    \left . \frac{\mathrm{d}^\ell (\mathrm{e}^{-\xi^2/2})}{\mathrm{d} \xi^\ell} \right |_{\xi = \mathcal{W} / \sqrt{\mathcal{V}}}
    \mathrm{e}^{\mathcal{W}^2/(2 \mathcal{V})} 
    \mathbb{E} \mathrm{e}^{\mathrm{i} \Theta}
    = \exp \left ( - \frac{\mathcal{V} + 2 \mathcal{W}}{2} \right ) 
    \mathbb{E} \mathrm{e}^{\mathrm{i} \Theta}.
\end{align*}
A rigorous justification of this argument is presented in the Supplementary Material.


\begin{funding}
G.M.Pan was partially supported by MOE Tier 1 grant of RG76/21 and MOE Tier 2 grant of T2EP20123-0046.

\end{funding}


\begin{supplement}
Supplement to ``Asymptotic distribution of spiked eigenvalues in large signal-plus-noise model''. This document comprises the technical proofs for Proposition \ref{prop-convergence-rate} and Theorem \ref{thm-characteristic-function}. 
\end{supplement}


\bibliographystyle{imsart-number} 
\bibliography{bibliography}       


\newpage
\setcounter{equation}{0}
\setcounter{section}{0}
\setcounter{page}{1}
\renewcommand{\theequation}{S.\arabic{equation}}
\renewcommand{\thesection}{S.\arabic{section}}
\renewcommand{\thepage}{S\arabic{page}}

\begin{frontmatter}

\title{Supplementary Material for ``Asymptotic distribution of spiked eigenvalues in large signal-plus-noise model''}
\runtitle{Supplementary Material}

\begin{aug}
\author[A]{\fnms{Zeqin}~\snm{Lin}},
\author[A]{\fnms{Guangming}~\snm{Pan}}
\author[B]{\fnms{Peng}~\snm{Zhao}}
\and
\author[C]{\fnms{Jia}~\snm{Zhou}}
\address[A]{School of Physical and Mathematical Sciences,
Nanyang Technological University\printead[presep={,\ }]{e1,e2}}

\address[B]{Department,
University or Company Name\printead[presep={,\ }]{e3}}

\address[C]{Department,
University or Company Name\printead[presep={,\ }]{e4}}
\end{aug}

\begin{abstract}
In Section \ref{sec-large-deviation}, we present the proof for the large deviation bound, Proposition \ref{prop-convergence-rate}. Section \ref{sec-recurrence-equation} is dedicated to establishing the asymptotic joint distribution of sesquilinear forms of resolvent, as presented in Theorem \ref{thm-characteristic-function}. Finally, in Section \ref{sec-two-resolvent}, we provide the derivation of a double-resolvent isotropic law for sample covariance matrices, which plays a crucial role in Section \ref{sec-recurrence-equation}.
\end{abstract}

\end{frontmatter}

\section{Large deviation bound}
\label{sec-large-deviation}

In this section we focus on the proof of Proposition \ref{prop-convergence-rate}. Recall the definition of $w_{+}$ in (\ref{def-critical-point}) and set $K_1 := \max \{ k : \tilde{\sigma}_k > -w_{+}^{-1} \}$. Then $\{ \tilde{\sigma}_k \}_{k \in \llbracket K_1 \rrbracket}$ are the outlier eigenvalues in spectrum of $\tilde{\Sigma}$ that go beyond the critical threshold. In this paper we do not discuss the behavior of $\lambda_k$ for $k \in \llbracket K_1 \rrbracket \backslash \llbracket K_0 \rrbracket$. Nevertheless, we can extend the definition $\theta_k = f(-\tilde{\sigma}_k^{-1})$ to all $k \in \llbracket K_1 \rrbracket$. Our arguments heavily rely on local laws at the $\theta_k$'s. Hence, we need to fix some small $\tau_1 > 0$ such that $\mathbb{D}_0 (\tau_1)$ accommodates all the $\theta_k$'s.

\begin{lemma}
\label{lemma-theta-location}
There exists some small $\tau_1 > 0$, depending only on $\tau$, such that 
\begin{equation}
    \lambda_+ + 2 \tau_1 
    \leq \theta_k 
    \leq \tau_1^{-1},
    \quad \forall k \in \llbracket K_0 \rrbracket.
    \label{eqn-location-bound}
\end{equation}
By choosing $\tau_1 > 0$ small enough, we can also ensure
\begin{equation}
    \theta_{k} - \theta_{k + 1}
    \geq \tau_1,
    \quad \forall k \in \llbracket K_0 \rrbracket,
    \label{eqn-location-separation}
\end{equation}
where we adopt the convention $\theta_{K_1 + 1} = \lambda_+$.
\end{lemma}

In the proof of Lemma \ref{lemma-theta-location} we need to use the following elementary lemma taken from \cite[Lemma A.3]{knowlesAnisotropicLocalLaws2017}, which concerns the behavior of $f$ near its critical point $w_+$. 

\begin{lemma}
\label{lemma-derivative-f-23}
We have
\begin{equation*}
    - w_+ \asymp 1 
    \quad \text{ and } \quad
    \lambda_+ = f(w_+) \asymp 1.
\end{equation*}
Furthermore, there exist constants $c, C, \tau_0 > 0$, depending only on $\tau$, such that
\begin{equation*}
    f^{\prime \prime} (w_+) \geq 2c
    \quad \text{ and } \quad
    | f^{\prime \prime \prime} (w) | \leq C, 
    \quad \text{ if } |w - w_+| \leq \tau_0.
\end{equation*}
\end{lemma}

\begin{proof}[Proof of Lemma \ref{lemma-theta-location}]
Note that $\sigma_1 < - w_+^{-1}$ since $w_+ \in (-\sigma_1^{-1}, 0)$. Recall the definition of $K_0$ in (\ref{def-K0}). We have
\begin{equation*}
    \tilde{\sigma}_k - \sigma_1 
    \geq \tilde{\sigma}_k + w_+^{-1}
    \geq 2 \tau \gtrsim 1,
    \quad \forall k \in \llbracket K_0 \rrbracket.
\end{equation*}
On the other hand, we have $\tilde{\sigma}_1 \lesssim 1$ as we assume $\| \Sigma \| \leq \tau^{-1}$ and $d_1 \leq \tau^{-1}$. Consequently,
\begin{equation*}
    \theta_k = f(- \tilde{\sigma}_k^{-1})
    = \tilde{\sigma}_k \left [ 1 + \phi \int \frac{s}{\tilde{\sigma}_k - s} \nu (\mathrm{d} s) \right ]
    \lesssim 1,
\end{equation*}
which concludes the proof of the upper bound in (\ref{eqn-location-bound}). While for the lower bound, we note that by $\tilde{\sigma}_k + w_+^{-1} \geq 2 \tau$,
\begin{equation*}
    - \tilde{\sigma}_k^{-1} - w_+ 
    \geq \frac{w_{+}}{1 - 2 \tau w_{+}} - w_+ 
    \geq \frac{2 \tau w_{+}^2}{1 - 2 \tau w_{+}}
    \gtrsim 1,
\end{equation*}
where we also used $-w_+ \asymp 1$ and that $\tau > 0$ is chosen to be sufficiently small. Let $c, C, \tau_0 > 0$ be constants in Lemma \ref{lemma-derivative-f-23}. Without loss of generality, we may assume $\tau_0$ is chosen sufficiently small such that $\tau_0 < c/ C$. Let $k \in \llbracket K_0 \rrbracket$. Utilizing $\theta_k = f(- \tilde{\sigma}_k^{-1})$ and the monotone increasing property of $f$ on $(w_+, 0)$, it suffices to consider the $k$'s with $- \tilde{\sigma}_k^{-1} < w_+ + \tau_0$. In this case, we have
\begin{equation*}
    \theta_k - \lambda_+
    = f(- \tilde{\sigma}_k^{-1}) - f(w_+) \\
    = f^{\prime \prime} (w_+) (- \tilde{\sigma}_k^{-1} - w_+)^2
    + f^{\prime \prime \prime} (w) (- \tilde{\sigma}_k^{-1} - w_+)^3,
\end{equation*}
where $w \in (- \tilde{\sigma}_k^{-1}, w_+ + \tau_0)$. Using $f^{\prime \prime} (w_+) \geq 2c$ and $| f^{\prime \prime \prime} (w) | \leq C \leq c/\tau_0$, we can obtain
\begin{align*}
    \theta_k - \lambda_+
    \geq c (- \tilde{\sigma}_k^{-1} - w_+)^2
    \gtrsim 1,
\end{align*}
which concludes the proof of the lower bound in (\ref{eqn-location-bound}). 

We now turn to the lower bound (\ref{eqn-location-separation}) for gaps between the $\theta_k$'s. Note that $m(\theta)$ is increasing for $\theta > \lambda_{+}$ with
\begin{equation}
    m^\prime (\theta) 
    = \int \frac{1}{(\lambda - \theta)^2} \varrho (\mathrm{d} \lambda) \asymp 1,
    \quad \text{ if } \theta - \lambda_{+} \gtrsim 1 \text{ and } \theta \lesssim 1.
    \label{eqn-estimate-mprime}
\end{equation}
Combined with $m(\theta_{k}) = -\tilde{\sigma}_k^{-1}$ and our Assumption \ref{assumption-spiked-spacing}, the estimate (\ref{eqn-estimate-mprime}) readily yields the lower bound (\ref{eqn-location-separation}) for $k \leq K_0 - 1$,
\begin{equation*}
    | \theta_{k} - \theta_{k+1} |
    \asymp | m(\theta_{k}) - m(\theta_{k+1}) |
    = | \tilde{\sigma}_k^{-1} - \tilde{\sigma}_{k+1}^{-1} |
    \asymp 1,
\end{equation*}
where we also used $\tilde{\sigma}_k \asymp 1$ since $\nu([0, \tau]) \leq 1 - \tau$ and $K$ is fixed. For $k = K_0$, it suffices to consider the case $K_1 > K_0$ since if $K_1 = K_0$ then by our convention $\theta_{K_0 + 1} = \theta_{K_1 + 1} = \lambda_+$. In this case, we may let $\tau_1 > 0$ be a constant such that (\ref{eqn-location-bound}) holds. If $\theta_{K_0 + 1} - \lambda_+ \leq \tau_1$, then we have automatically $\theta_{K_0} - \theta_{K_0 + 1} \geq \tau_1$. Otherwise, we have $\theta_{K_0 + 1} - \lambda_+ > \tau_1$, for which we may apply the estimate (\ref{eqn-estimate-mprime}) again to deduce
\begin{equation*}
    | \theta_{K_0} - \theta_{K_0+1} | 
    \asymp | m(\theta_{K_0}) - m(\theta_{K_0+1}) |
    = | \tilde{\sigma}_{K_0}^{-1} - \tilde{\sigma}_{K_0+1}^{-1} | \asymp 1.
\end{equation*}
Now (\ref{eqn-location-separation}) follows after some appropriate renaming of $\tau_1 > 0$.
\end{proof}

In what follows, we fix some $\tau_1 > 0$ such that (\ref{eqn-location-bound}) and (\ref{eqn-location-separation}) hold. To prove Proposition \ref{prop-convergence-rate}, we need to show that for arbitrary small $\varepsilon > 0$ and large $L > 0$,
\begin{equation*}
    | \lambda_k - \theta_k | \leq N^{-1/2 + 2 \varepsilon}, 
    \quad \forall k \in \llbracket K_0 \rrbracket
\end{equation*}
holds with probability at least $1 - N^{-L}$. To this end, we define 
\begin{equation}
    \mathcal{I} := \bigcup_{k=1}^{K_0} \mathcal{I}_k,
    \quad \text{ where } \quad \mathcal{I}_k = [\theta_{k-}, \theta_{k+}]
    \text{ with } \theta_{k \pm} := \theta_k \pm N^{-1/2 + 2 \varepsilon}.
    \label{def-interval}
\end{equation}
We need to prove, with probability at least $1 - N^{-L}$, each interval $\mathcal{I}_k$ contains exactly one spiked eigenvalue $\lambda_k$. 

Applying the Weinstein–Aronszajn identity, we may write
\begin{align}
\begin{split}
    \det A_{\Pi} & (z) 
    = (-1)^K \det \left [ U^\top \frac{m(z)}{I + m (z) \Sigma} U + D^{-2} \right ] \\
    & = (-1)^K (\det D)^{-2}
    \frac{\det \big [ {I + m (z) \Sigma + m(z) U D^2 U^\top} \big ]}{\det \big [{I + m(z) \Sigma} \big ]} \\
    & = (-1)^K (\det D)^{-2}
    \frac{\det \big [ {I + m (z) \tilde{\Sigma}} \big ]}{\det \big [ {I + m (z) \Sigma} \big ]}
    = (-1)^K (\det D)^{-2} \prod_{i = 1}^M
    \frac{1 + m(z) \tilde{\sigma}_i}{1 + m(z) \sigma_i}.    
\end{split} \label{eqn-determinant-expression}
\end{align}

The next two lemmas provide some useful estimates regarding the quantities $|1 + m(z) \tilde{\sigma}_i|$. For $z = E + \mathrm{i} \eta \in \mathbb{C} \backslash \operatorname{supp} \varrho$, we write $m(z) = m_1(z) + \mathrm{i} m_2(z)$. Let us first introduce some estimates that are useful when extending estimates on the real line to general complex spectral parameters,
\begin{align}
    0 < m_1(z) - m(E) 
    = \int \frac{\eta^2}{(E - \lambda) [ (\lambda - E)^2 + \eta^2 ]} \varrho (\mathrm{d} \lambda)
    \lesssim \eta^2,
    \quad \text{ if } E \in \mathbb{D}_0 
    \label{eqn-extension-real} \\
    | m_2(z) | = \int \frac{| \eta |}{(\lambda - E)^2 + \eta^2} \varrho (\mathrm{d} \lambda)
    \lesssim | \eta |,
    \quad \text{ if } E \in \mathbb{D}_0
    \label{eqn-extension-imaginary}
\end{align}
where we recall $\mathbb{D}_0 = [\lambda_+ + \tau_1, \tau_1^{-1}]$. We also use the following estimate of $m^\prime(z)$,
\begin{equation}
    m^\prime (z) = \int \frac{1}{(\lambda - z)^2} \varrho (\mathrm{d} \lambda) \gtrsim 1,
    \quad \forall z \in (\lambda_{+}, \tau_1^{-1}].
    \label{eqn-extension-derivative}
\end{equation}

\begin{lemma}
\label{lemma-1+ms-non-spiked}
There exists some small $\tau_2 > 0$, which depends only on $\tau$ and $\tau_1$, such that the following holds uniformly in $z = E + \mathrm{i} \eta$ with $E \in \mathbb{D}_0$ and $| \eta | \leq \tau_2$,
\begin{align*}
    |1 + m(z) \sigma_{i-K}| 
    \leq |1 + m(z) \tilde{\sigma}_i|
    \leq | 1 + m(z) \sigma_i|,
    & \quad \forall K+1 \leq i \leq M, \\
    |1 + m(z) \tilde{\sigma}_{k}| \asymp 1,
    & \quad \forall K_1+1 \leq k \leq K.
\end{align*}
\end{lemma}

\begin{proof}[Proof of Lemma \ref{lemma-1+ms-non-spiked}]
Consider the inequalities in the first line. By monotonicity of $m(z)$ on $(\lambda_{+}, \tau_1^{-1}]$ and Assumption \ref{assumption-regularity}, we have
\begin{equation*}
    1 + \sigma_1 m(E) 
    > 1 + \sigma_1 m(\lambda_+)
    = \sigma_1 ( w_+ + \sigma_1^{-1} )
    \geq \tau^2,
    \quad \forall E \in (\lambda_{+}, \tau_1^{-1}].    
\end{equation*}
Using the estimates (\ref{eqn-extension-real}), (\ref{eqn-extension-imaginary}) and 
\begin{equation*}
    | m_1(z) | = \int \frac{E - \lambda}{(\lambda - E)^2 + \eta^2} \varrho (\mathrm{d} \lambda)
    \gtrsim E - \lambda_+ \gtrsim 1,
    \quad \text{ if } E \in \mathbb{D}_0,
\end{equation*}
we find that for small enough $\tau_2 > 0$, the following holds as long as $E \in \mathbb{D}_0$ and $| \eta | \leq \tau_2$,
\begin{equation*}
    | m_2(z) |^2 / | m_1(z) | \leq \tau^3 / 4
    \quad \text{ and } \quad
    0 < m_1(z) - m(E) \leq \tau^3 / 4.
\end{equation*}
Then, using $\sigma_1 < \tau^{-1}$, we can deduce
\begin{equation*}
    1 + \sigma_1 \left [ m_1(z) + \frac{m_2(z)^2}{m_1(z)} \right ]
    \geq 1 + \sigma_1 m(E) - \sigma_1 \left | m_1(z) - m(E) + \frac{m_2(z)^2}{m_1(z)} \right |
    \geq \tau^2 - \frac{\tau^2}{2}
    = \frac{\tau^2}{2}.    
\end{equation*}
Note that $m_1 (z) < 0$ when $E \in \mathbb{D}_0$. We therefore conclude
\begin{equation}
    \sigma_1
    < - \frac{1}{m_1(z) + m_2(z)^2 / m_1(z)} = -\frac{m_1 (z)}{|m(z)|^2},
    \quad \text{ if } E \in \mathbb{D}_0
    \text{ and } | \eta | \leq \tau_2.
    \label{eqn-sigma1-bound-by-m}
\end{equation}
Since the mapping
\begin{equation*}
    \sigma \quad \mapsto \quad |1 + m(z) \sigma|^2 
    = 1 + 2 m_1(z) \sigma + | m(z) |^2 \sigma^2    
\end{equation*}
is decreasing for $\sigma \in [0, -m_1(z)/|m(z)|^2 ]$, we can conclude the proof of the first line using (\ref{eqn-sigma1-bound-by-m}) and the interlacing inequalities $\sigma_i \leq \tilde{\sigma}_i \leq \sigma_{i - K}$ for $K+1 \leq i \leq M$. 

Now let us turn to the proof of the second line. Using the derivative estimate (\ref{eqn-extension-derivative}) and that $\tilde{\sigma}_k \leq - w_+^{-1}$ for $K_1 + 1 \leq k \leq K$, we have for $E \in \mathbb{D}_0$,
\begin{align*}
    1 + m(E) \tilde{\sigma}_k 
    & \geq [1 + m(E) \tilde{\sigma}_k] - [1 + m(\lambda_+) \tilde{\sigma}_k] \\
    & = [m(E) - m(\lambda_+)] \tilde{\sigma}_k
    \gtrsim (E - \lambda_+) \tilde{\sigma}_k
    \gtrsim 1,
\end{align*}
where we used $\tilde{\sigma}_k \asymp 1$ again. To extend this lower bound to complex spectral parameters, we may activate (\ref{eqn-extension-real}) and (\ref{eqn-extension-imaginary}) again to get
\begin{equation*}
    \big | [1 + m(z) \tilde{\sigma}_k] - [1 + m(E) \tilde{\sigma}_k] \big |
    \lesssim |m_1(z) - m(E)| + |m_2(z)|
    \lesssim \eta^2 + |\eta|.
\end{equation*}
Therefore, by letting $\tau_2 > 0$ small enough, one can ensure 
\begin{equation*}
    |1 + m(z) \tilde{\sigma}_k| \gtrsim 1,
    \quad \forall K_1 + 1 \leq k \leq K
\end{equation*}
as long as $E \in \mathbb{D}_0$ and $| \eta | \leq \tau_2$. This concludes our proof of the second line.
\end{proof}

In what follows, we shall fix some $\tau_2 > 0$ such that the inequalities in Lemma \ref{lemma-1+ms-non-spiked} hold. Without loss of generality, we assume $\tau_2 \leq \tau_1/2$. Consider the positive-oriented contours
\begin{equation*}
    \mathcal{C}_k := \{ z \in \mathbb{C}: |z - \theta_k| = \tau_2 \},
    \quad \forall k \in \llbracket K_0 \rrbracket,
\end{equation*}
and set $\mathcal{C} := \bigcup_{k = 1}^{K_0} \mathcal{C}_k$. By Lemma \ref{lemma-theta-location}, we have $\mathcal{C} \subset \mathbb{D}$ and each contour $\mathcal{C}_k$ encloses $\theta_k$ but no other $\theta_j$, $j \not= k$. 

\begin{lemma}
\label{lemma-1+ms-spiked}
The following holds uniformly in $z \in \mathbb{D}_0$,
\begin{equation*}
    |1 + m(z) \tilde{\sigma}_k| 
    \gtrsim |z - \theta_k|,
    \quad \forall k \in \llbracket K_1 \rrbracket.
\end{equation*}
On the other hand, the following holds uniformly in $z \in \mathcal{C}$,
\begin{equation*}
    |1 + m(z) \tilde{\sigma}_k| 
    \gtrsim 1,
    \quad \forall k \in \llbracket K_1 \rrbracket.
\end{equation*}
\end{lemma}

\begin{proof}[Proof of Lemma \ref{lemma-1+ms-spiked}]
The first estimate follows directly from $1 + m(\theta_k) \tilde{\sigma}_k = 0$ and the derivative estimate (\ref{eqn-extension-derivative}). For the second estimate, we note that by (\ref{eqn-extension-real}-\ref{eqn-extension-derivative}), there exist constants $c \leq 1$ and $C \geq 1$ such that we have the following uniformly in $k \in \llbracket K_1 \rrbracket$,
\begin{equation}
    | m(z) - m(E) | \tilde{\sigma}_k \leq C | \eta |,
    \quad \forall z \in \mathbb{D}
    \quad \text{ and } \quad 
    m^\prime(z) \tilde{\sigma}_k \geq c,
    \quad \forall z \in (\lambda_{+}, \tau_1^{-1}],
    \label{eqn-contour-estimate}
\end{equation}
where we used $\tilde{\sigma}_k \asymp 1$ for $k \in \llbracket K_1 \rrbracket$.
Consider $z \in \mathcal{C} \subset \mathbb{D}$. If $| \eta | \geq c\tau_2 / (2C)$, then
\begin{equation*}
    |1 + m(z) \tilde{\sigma}_k| 
    \geq |m_2(z)| \tilde{\sigma}_k
    \asymp | \eta | \gtrsim 1,
\end{equation*}
where we used $|m_2(z)| \asymp |\eta|$ for $z \in \mathbb{D}$. On the other hand, if $\eta < c\tau_2 / (2C) \leq \tau_2/2$, then by Lemma \ref{lemma-theta-location}, $\tau_2 \leq \tau_1/2$ and definition of the contours $\mathcal{C}_k$, we must have 
\begin{equation*}
    |E - \theta_k| \geq \tau_2/2,
    \quad \forall k \in \llbracket K_1 \rrbracket.
\end{equation*}
Note that $E \in \mathbb{D}_0 \subset (\lambda_{+}, \tau_1^{-1}]$ and $\{ \theta_k \}_{k \in \llbracket K_1 \rrbracket} \subset (\lambda_{+}, \tau_1^{-1}]$. Hence, in this case, by the first estimate of (\ref{eqn-contour-estimate}), we have
\begin{equation*}
    | 1 + m(E) \tilde{\sigma}_k | 
    = |m(E) - m(\theta_k)| \tilde{\sigma}_k 
    \geq c |E - \theta_k|
    \geq c\tau_2/2.
\end{equation*}
Together with the second estimate of (\ref{eqn-contour-estimate}), we get
\begin{equation*}
    | 1 + m(z) \tilde{\sigma}_k | 
    \geq | 1 + m(E) \tilde{\sigma}_k | - | m(z) - m(E) | \tilde{\sigma}_k
    \geq c\tau_2/2 - C \eta
    \geq c\tau_2/4
    \gtrsim 1,
\end{equation*}
which concludes the proof.
\end{proof}

We are now in a position to present the proof of Proposition \ref{prop-convergence-rate}.

\begin{proof}[Proof of Proposition \ref{prop-convergence-rate}]
According to Proposition \ref{prop-local-law} and the definitions of $A_{G} (z)$ and $A_{\Pi} (z)$ in (\ref{def-mtr-AA-matrices}), we have $A_{G} (z) = A_{\Pi} (z) + \oprec (N^{-1/2})$. Consequently,
\begin{equation*}
    \det A_{G} (z) = \det A_{\Pi} (z) + \oprec (N^{-1/2})
\end{equation*}
uniformly in $z \in \mathbb{D}$. Using the estimates in Lemma \ref{lemma-basic-estimates-stieltjes} and $\| G \| \prec 1$ for $z \in \mathbb{D}$, both $\det A_{G} (z)$ and $\det A_{\Pi} (z)$ are $N^{C_0}$-Lipschitz in $\mathbb{D}$ with high probability for some constant $C_0 > 0$. Therefore, via a standard $\varepsilon$-net argument, we can extend the above probabilistic bound such that it holds simultaneously for all $z \in \mathbb{D}$, i.e. for any fixed $\varepsilon > 0$ and $L > 0$,
\begin{equation}
    | {\det A_{G} (z) - \det A_{\Pi} (z)} | \leq N^{-1/2 + \varepsilon},
    \quad \forall z \in \mathbb{D}
    \label{eqn-AG-APi-bound}
\end{equation}
holds with probability at least $1 - N^{-2L}$. We refer to Remark 2.7 of \cite{benaych-georgesLecturesLocalSemicircle2018} for more details. Furthermore, by Assumption \ref{assumption-deformation} and the rigidity estimate (\ref{eqn-eigenvalue-rigidity}), the following event also holds with probability at least $1 - N^{-2L}$,
\begin{equation*}
    \lambda_1 \leq \tau_1^{-1}.
\end{equation*}
Let $\Omega$ be the intersection of the two high probability events mentioned above. Our analysis will be conditioned on $\Omega$. Recall that by (\ref{eqn-determinant-expression}),
\begin{equation*}
    \det A_{\Pi} (z) 
    = (-1)^K (\det D)^{-2} \prod_{i = 1}^M
    \frac{1 + m(z) \tilde{\sigma}_i}{1 + m(z) \sigma_i},
\end{equation*}
where $\det D \asymp 1$ due to Assumption \ref{assumption-deformation}. By Lemma \ref{lemma-1+ms-non-spiked}, for $z \in \mathbb{D}$ we have
\begin{equation*}
    \frac{\prod_{i = K+1}^M | 1 + m(z) \tilde{\sigma}_i |}{\prod_{i = 1}^M | 1 + m(z) \sigma_i |} \asymp 1
    \quad \text{ and } \quad
    \prod_{k = K_1+1}^{K} | 1 + m(z) \tilde{\sigma}_k | \gtrsim 1,
\end{equation*}
where in the first estimate we used
\begin{equation*}
    \prod_{i = M-K+1}^M |1 + m(z) \sigma_i|^{-1}
    \leq \frac{\prod_{i = K+1}^M | 1 + m(z) \tilde{\sigma}_i |}{\prod_{i = 1}^M | 1 + m(z) \sigma_i |}
    \leq \prod_{i = 1}^K |1 + m(z) \sigma_i|^{-1}
\end{equation*}
and $| 1 + m(z) \sigma_i | \asymp 1$ from Lemma \ref{lemma-basic-estimates-stieltjes}. We firstly claim that, on the event $\Omega$, we have 
\begin{equation}
    \det A_{G} (z) \not= 0,
    \quad \forall z \in [\theta_{K_0} - \tau_1/2, \infty) \backslash \mathcal{I}.
    \label{eqn-singularity-determinant}
\end{equation}
Indeed, as $\lambda_1 < \tau_1^{-1}$, it suffices to consider $z \in [\theta_{K_0} - \tau_1/2, \tau_1^{-1}] \backslash \mathcal{I}$. Note that $z \in \mathbb{D}_0$. Set
\begin{equation*}
    k_0 \equiv k_0 (z) := \underset{j \in \llbracket K_0 \rrbracket}{\arg \min} |z - \theta_j|.
\end{equation*}
Then, by the definition of $\mathcal{I}$ in (\ref{def-interval}) and Lemma \ref{lemma-theta-location}, we have
\begin{equation*}
    |z - \theta_{k_0}| \geq N^{-1/2 + 2\varepsilon} 
    \quad \text{ and } \quad
    |z - \theta_j| \geq \tau_1 / 2,
    \quad \forall j \in \llbracket K_1 \rrbracket \backslash \{ k \}.
\end{equation*}
Applying Lemma \ref{lemma-1+ms-spiked}, we therefore deduce
\begin{equation*}
    | {\det A_{\Pi} (z)} |
    \gtrsim \prod_{j = 1}^{K_1} | 1 + m(z) \tilde{\sigma}_j | 
    \gtrsim N^{-1/2 + 2\varepsilon}.
\end{equation*}
Combined with (\ref{eqn-AG-APi-bound}), we conclude that for all sufficiently large $N$,
\begin{equation*}
    | {\det A_{\Pi} (z)} | > | {\det A_{\Pi} (z) - \det A_{G} (z)} |,
    \quad \forall z \in [\theta_{K_0} - \tau_1/2, \infty) \backslash \mathcal{I}.
\end{equation*}
from which then claim (\ref{eqn-singularity-determinant}) follows. Next, we claim that on $\Omega$, each interval $\mathcal{I}_k$ contains exactly one spiked eigenvalue $\lambda_k$. Indeed, let $k \in \llbracket K_0 \rrbracket$ and consider any $z \in \mathcal{C}_k$. Note that in this case we have $|z - \theta_j| \geq \tau_2$ for all $j \in \llbracket K_1 \rrbracket$. Applying Lemma \ref{lemma-1+ms-spiked} again, we have $| {\det A_{\Pi} (z)} | \gtrsim 1$ and therefore
\begin{equation*}
    | {\det A_{\Pi} (z)} | > | {\det A_{\Pi} (z) - \det A_{G} (z)} |,
    \quad \forall z \in \mathcal{C}_k.
\end{equation*}
Now, applying Roch\'{e}'s theorem, $\det A_{G} (z)$ and $\det A_{\Pi} (z)$ have the same number of zeros inside $\mathcal{C}_k$. Since $\det A_{\Pi} (z)$ contains exactly one zero $\theta_k$ inside $\mathcal{C}_k$, we know that $\det A_{G} (z)$ also have exactly one zero inside $\mathcal{C}_k$. Combined with (\ref{eqn-singularity-determinant}), this zero must be contained in $\mathcal{I}_k$. In conclusion, we have proved that on the high probability event $\Omega$,
\begin{equation*}
    \lambda_k \in \mathcal{I}_k
    \quad \text{ i.e. } \quad
    |\lambda_k - \theta_k| \leq N^{-1/2 + 2\varepsilon},
    \quad \forall k \in \llbracket K_0 \rrbracket.
\end{equation*}
This completes the proof of Proposition \ref{prop-convergence-rate}.
\end{proof}

\section{Recursive estimate}
\label{sec-recurrence-equation}

This section is dedicated to the proof of our main technical result, Theorem \ref{thm-characteristic-function}. As previously mentioned, the proof involves incorporating the cumulant expansion technique and the local laws for the undeformed ensemble. Here we utilize the following version of cumulant expansion formula, whose proof can be obtained through a minor modification of the proof for \cite[Proposition 3.1]{lytovaCentralLimitTheorem2009} and, therefore, we omit it for brevity.

\begin{lemma}[cumulant expansion formula]
\label{lemma-cumulant-expansion}
Let $n \in \mathbb{N}$ be fixed and $g \in C^{n+1}(\mathbb{R})$. Supposed $\xi$ is a centered random variable with finite moments to order $n + 2$. Denote the $r$-th cumulant of $\xi$ by $\kappa_r (\xi)$. Then, we have
\begin{equation*}
    \mathbb{E} [ \xi g(\xi) ]
    = \sum_{r=0}^{n} \frac{\kappa_{r+1}(\xi)}{r !} 
    \mathbb{E} g^{(r)}(\xi) + \mathcal{R}_n,
\end{equation*}
assuming that all expectations on the r.h.s. exist, where the remainder term $\mathcal{R}_n$ depends on $g$ and $\xi$ and satisfies, for any $s > 0$,
\begin{equation*}
    | \mathcal{R}_n | 
    \leq \hat{C}_{n} \left [ \mathbb{E} |\xi|^{n+2} 
    \cdot \sup_{|t| \leq s} \big | { g^{(n+1)}(t) } \big |
    + \mathbb{E} \big [ |\xi|^{n+2} \mathbbm{1}(|\xi|>s) \big ]
    \cdot \sup_{t \in \mathbb{R}} \big | { g^{(n+1)}(t) } \big | \right ].
\end{equation*}
Here $\hat{C}_{n} \leq (Cn)^n / n!$ for some absolute constant $C > 0$.
\end{lemma}

We also employ the following elementary lemma, which summarizes the key property of stochastic domination when combined with expectation. It can be proved by a straightforward application of the Cauchy-Schwartz inequality.

\begin{lemma}
\label{lemma-expectation-domination}
Let $\xi \prec \zeta$ where $\zeta$ is deterministic. If there exists some large constant $C > 0$ such that $\mathbb{E} \xi^2 \leq N^C$ and $\zeta \geq N^{-C}$, then $\mathbb{E} |\xi| \prec \zeta$.
\end{lemma}

In particular, this lemma introduces a technical challenge that we must address before delving into the proof. Specifically, to translate the high probability bounds encoded by $\prec$ into bounds for expectations, a crude deterministic bound for $\| G \|$ is usually required if we intend to apply Lemma \ref{lemma-expectation-domination}. Unfortunately, such a bound is typically unavailable for real $z$. However, as explained in \cite{baoSingularVectorSingular2021}, we can overcome this obstacle by fixing some sufficiently large constant $C_1 > 0$ and considering spectral parameters from 
\begin{equation*}
    \mathbb{D}_1 \equiv \mathbb{D}_1^{N} (\tau_1, \tau_2, C_1) := 
    \{ E + \mathrm{i} N^{-C_1}: E \in [\lambda_+ + \tau_1, \tau_1^{-1}] \}.
\end{equation*}
We will first prove Theorem \ref{thm-characteristic-function} by assuming $\{ z_k \} \subset \mathbb{D}_1$, where we can freely employ the elementary bound $\| G \| \leq N^{C_1}$. The result can then be extended to $\{ z_k \} \subset \mathbb{D}_0$ via a simple continuity argument. Note that the asymptotic quantities in Definition \ref{def-sesquilinear-quantities} are still well defined when $\{ z_k \} \subset \mathbb{D}_1$. Now, we are prepared to formalize the recursive estimate (\ref{eqn-approx-recurrence-moment}) and provide a rigorous justification for the argument mentioned at the end of Section \ref{subsec-sesquilinear-forms}.

\begin{proposition}
\label{prop-recurrence-moment}
For each fixed $\ell \in \mathbb{N}$, the following holds uniformly in $\{ z_k \} \subset \mathbb{D}_1$,
\begin{equation}
    \mathbb{E} [ {\Phi^{\ell + 2} \mathrm{e}^{\mathrm{i} \Theta}} ]
    = (\ell + 1) \mathcal{V} \mathbb{E} [ {\Phi^{\ell} \mathrm{e}^{\mathrm{i} \Theta}} ]
    + (\mathrm{i} \mathcal{W}) \mathbb{E} [ {\Phi^{\ell + 1} \mathrm{e}^{\mathrm{i} \Theta}} ]
    + \oprec (N^{-1/2}).
    \label{eqn-recurrence-moment-ell}
\end{equation}
In addition,
\begin{equation}
    \mathbb{E} [ {\Phi \mathrm{e}^{\mathrm{i} \Theta}} ]
    = (\mathrm{i} \mathcal{W}) \mathbb{E} [\mathrm{e}^{\mathrm{i} \Theta}]
    + \oprec (N^{-1/2}).
    \label{eqn-recurrence-moment-one}
\end{equation}
\end{proposition}

\begin{remark}
Letting $t_k = 0$ for all $k$ in Proposition \ref{prop-recurrence-moment} yields
\begin{equation*}
    \mathbb{E} \Phi^{2 \ell + 2}
    = (\ell + 1) \mathcal{V} \mathbb{E} \Phi^{2 \ell} 
    + \oprec (N^{-1/2}) .
\end{equation*}
When $\{ z_k \} \subset \mathbb{D}_1$, both $\Phi$ and the asymptotic quantity $\mathcal{V}$ are complex, so that the above recursive estimate is insufficient for concluding the asymptotic Gaussianity of $\Phi$. However, the argument we present in this section can be easily adapted to develop a recursive estimate for $\mathbb{E} [\Phi^{\ell_1} (\Phi^*)^{\ell_2}]$ where $\ell_1, \ell_2 \in \mathbb{N}$, which can be used to justify the the asymptotic Gaussianity of $\Phi$. In particular, for $\mathbb{E} |\Phi|^{2 \ell} = \mathbb{E} [\Phi^{\ell} (\Phi^*)^{\ell}]$, one can deduce
\begin{equation}
    \mathbb{E} |\Phi|^{2 \ell + 2}
    = (\ell + 1) \mathcal{V}_{\mathrm{c}} \mathbb{E} |\Phi|^{2 \ell} 
    + \oprec (N^{-1/2}),
    \label{eqn-recurrence-moment-module}
\end{equation}
where $\mathcal{V}_{\mathrm{c}} \lesssim 1$ is real and agrees with $\mathcal{V}$ when $\{ z_k \} \subset \mathbb{D}_0$. We will utilize (\ref{eqn-recurrence-moment-module}) to derive a preliminary estimate for $\mathbb{E} |\Phi|^{2 \ell}$ with $\{ z_k \} \subset \mathbb{D}_1$. However, for the sake of simplicity, we will not present the proof of (\ref{eqn-recurrence-moment-module}) and the exact expression of $\mathcal{V}_{\mathrm{c}}$.
\end{remark}

\begin{proof}[Proof of Theorem \ref{thm-characteristic-function}]
As highlighted, we firstly present the proof of (\ref{eqn-characteristic-function}) for $\{ z_k \} \subset \mathbb{D}_1$. Without loss of generality, we may assume $\max_k |s_k| \asymp 1$. Then, by (\ref{assumption-coefficient-size}) and Lemma \ref{lemma-basic-estimates-stieltjes}, one can easily verify that $| \mathcal{V} | \asymp 1$. On the other hand, we have $\mathcal{V}_{\mathrm{c}}, |\mathcal{W}| \lesssim 1$, where $\mathcal{V}_{\mathrm{c}}$ is introduced in the remark following Proposition \ref{prop-recurrence-moment}.

Fix arbitrary $\varepsilon > 0$. By Taylor's expansion,
\begin{equation*}
    \bigg | \mathrm{e}^{\mathrm{i} \xi}
    - \sum_{\ell = 0}^{2n-1} \frac{(\mathrm{i} \xi)^\ell}{\ell !} \bigg | \leq \frac{|\xi|^{2n}}{(2n) !},
    \quad \forall \xi \in \mathbb{R}.
\end{equation*}
Let $n \in \mathbb{N}$ be sufficiently large such that $\mathcal{V}_{\mathrm{c}}^n / n! \leq \varepsilon$. Note that $n$ is $N$-independent. Therefore, a recursive application of (\ref{eqn-recurrence-moment-module}) implies
\begin{equation*}
    \mathbb{E} |\Phi|^{2n}
    \leq (2n - 1) !! (2 \mathcal{V}_{\mathrm{c}})^n.
\end{equation*}
As a result, for large enough $n \in \mathbb{N}$,
\begin{equation}
    \left | \mathbb{E} \exp [ {\mathrm{i} (\Phi + \Theta)} ]
    - \sum_{\ell = 0}^{2n-1} \frac{1}{\ell !} \mathbb{E} \big [ {(\mathrm{i} \Phi)^{\ell} \mathrm{e}^{\mathrm{i} \Theta}} \big ] \right | \\
    \leq \frac{\mathbb{E} |\Phi|^{2 n}}{(2n)!} \leq \frac{\mathcal{V}_{\mathrm{c}}^n}{n!} 
    \leq \varepsilon.
    \label{eqn-taylor-expansion-1}
\end{equation}
By the definition of Hermite polynomials, we have $\mathrm{He}_{0} (\xi) = 1, \mathrm{He}_{1} (\xi) = \xi$. Therefore, using (\ref{eqn-recurrence-moment-one}), we have the following for $\ell = 0, 1$,
\begin{equation}
    \mathbb{E} [ {\Phi^{\ell} \mathrm{e}^{\mathrm{i} \Theta}} ]
    = (\mathrm{i} \sqrt{\mathcal{V}})^{\ell} \mathrm{He}_{\ell} (\mathcal{W} / \sqrt{\mathcal{V}}) 
    \mathbb{E} [{\mathrm{e}^{\mathrm{i} \Theta}}]
    + \mathcal{O}_{\prec} (N^{-1/2}).
    \label{eqn-hermite-estimate}
\end{equation}
Applying (\ref{eqn-recurrence-moment-ell}), we can generalize (\ref{eqn-hermite-estimate}) to arbitrary fixed $\ell \geq 0$ via a simple induction,
\begin{align*}
    & \ \mathbb{E} [ {\Phi^{\ell + 2} \mathrm{e}^{\mathrm{i} \Theta}} ] \\
    = & \ (\ell + 1) \mathcal{V} \mathbb{E} [ {\Phi^{\ell} \mathrm{e}^{\mathrm{i} \Theta}} ] 
    + (\mathrm{i} \mathcal{W}) \mathbb{E} [ {\Phi^{\ell + 1} \mathrm{e}^{\mathrm{i} \Theta}} ]
    + \mathcal{O}_{\prec} (N^{-1/2}) \\
    = & \ (\mathrm{i} \sqrt{\mathcal{V}})^{\ell + 2} \big [ {
    - (\ell + 1) \mathrm{He}_{\ell} (\mathcal{W} / \sqrt{\mathcal{V}}) 
    + (\mathcal{W} / \sqrt{\mathcal{V}}) \mathrm{He}_{\ell + 1} (\mathcal{W} / \sqrt{\mathcal{V}}) 
    } \big ] \mathbb{E} [{\mathrm{e}^{\mathrm{i} \Theta}}] 
    + \mathcal{O}_{\prec} (N^{-1/2}) \\
    = & \ (\mathrm{i} \sqrt{\mathcal{V}})^{\ell + 2} \mathrm{He}_{\ell + 2} (\mathcal{W} / \sqrt{\mathcal{V}}) 
    \mathbb{E} [{\mathrm{e}^{\mathrm{i} \Theta}}] + \mathcal{O}_{\prec} (N^{-1/2}),
\end{align*}
where in the last step we used the recurrence relation of Hermite polynomials. Let us denote $g(\xi) = \mathrm{e}^{-\xi^2/2}$ and recall $\mathrm{He}_{\ell} (\xi) = (-1)^\ell \mathrm{e}^{\xi^2/2} g^{(\ell)} (\xi)$. Therefore, (\ref{eqn-hermite-estimate}) implies that, for all $0 \leq \ell \leq 2n-1$,
\begin{equation}
    \frac{1}{\ell !} \mathbb{E} [ {(\mathrm{i}\Phi)^{\ell} \mathrm{e}^{\mathrm{i} \Theta}} ] 
    = \frac{g^{(\ell)} (\mathcal{W} / \sqrt{\mathcal{V}})}{\ell !} 
    \mathcal{V}^{\ell/2} 
    \cdot \mathrm{e}^{\mathcal{W}^2/(2 \mathcal{V})} 
    \mathbb{E} [{\mathrm{e}^{\mathrm{i} \Theta}}]
    + \mathcal{O}_{\prec} (N^{-1/2}).
    \label{eqn-taylor-expansion-2}
\end{equation}
By Taylor's expansion, we also have
\begin{equation*}
    g( \mathcal{W} / \sqrt{\mathcal{V}} + \sqrt{\mathcal{V}} ) =
    \sum_{\ell = 0}^{2 n - 1} \frac{g^{(\ell)} (\mathcal{W} / \sqrt{\mathcal{V}})}{\ell !} 
    \mathcal{V}^{\ell/2} 
    + \frac{g^{(2n)} (\zeta)}{(2n) !} \mathcal{V}^{n},
\end{equation*}
where $\zeta$ lies between $\mathcal{W} / \sqrt{\mathcal{V}}$ and $\mathcal{W} / \sqrt{\mathcal{V}} + \sqrt{\mathcal{V}}$. In particular, $\zeta$ belongs to some $N$-independent compact set since $|\mathcal{V}| \asymp 1$ and $|\mathcal{W}| \lesssim 1$. To control $g^{(2n)} (\zeta) = \mathrm{e}^{-\zeta^2/2} \mathrm{He}_{2n} (\zeta)$, we can apply the following asymptotic expansion of Hermite polynomials deduced from \cite[Theorem 8.22.1]{szeg1939orthogonal},
\begin{equation*}
    \mathrm{He}_{2n}(\zeta)
    = (-1)^n 2^{n} n ! \left [ \pi^{-1/2} \mathrm{e}^{\zeta^2/4} n^{-1/2}
    \cos (\sqrt{2 n \zeta^2} ) + O(n^{-1}) \right ],
\end{equation*}
where the reminder term is uniform in compact sets. Hence, by letting $n \in \mathbb{N}$ large enough, 
\begin{equation}
    \left | {g( \mathcal{W} / \sqrt{\mathcal{V}} + \sqrt{\mathcal{V}} ) - \sum_{\ell = 0}^{2 n - 1} \frac{g^{(\ell)} (\mathcal{W} / \sqrt{\mathcal{V}})}{\ell !} \mathcal{V}^{\ell/2} } \right |
    \mathrm{e}^{\mathcal{W}^2/(2 \mathcal{V})} 
    \leq \varepsilon.
    \label{eqn-taylor-expansion-3}
\end{equation}
Summarising the inequalities (\ref{eqn-taylor-expansion-1}) to (\ref{eqn-taylor-expansion-3}), we get
\begin{equation*}
    \left | \mathbb{E} \exp [ \mathrm{i} (\Phi + \Theta) ]
    - g( \mathcal{W} / \sqrt{\mathcal{V}} + \sqrt{\mathcal{V}} ) \cdot
    \mathrm{e}^{\mathcal{W}^2/(2 \mathcal{V})} 
    \mathbb{E} \mathrm{e}^{\mathrm{i} \Theta} \right |
    \leq 2 \varepsilon + \mathcal{O}_{\prec} (N^{-1/2}).
\end{equation*}
Since $\varepsilon > 0$ is arbitrary, by plugging $g(\xi) = \mathrm{e}^{-\xi^2/2}$ we therefore conclude the proof of (\ref{eqn-characteristic-function}) when $\{ z_k \} \subset \mathbb{D}_1$. 

Now we turn to the case where $\{ z_k \} \subset \mathbb{D}_0$. We can introduce the perturbed spectral parameters by
\begin{equation*}
    \hat{z}_k = z_k + \mathrm{i} N^{-C_1} \in \mathbb{D}_1.
\end{equation*}
Actually, using the estimates in Lemma \ref{lemma-basic-estimates-stieltjes} and $\| H \|, \| G \| \prec 1$ for $z \in \mathbb{D}$, we can fix some $C_0 > 0$ such that for arbitrary large $L > 1$, there exists some high-probability event $\Omega$ with $\mathbb{P}(\Omega) \geq 1 - N^{-L}$ on which $\Phi (\mathbf{z})$ and $\Theta ((\mathbf{z}))$ are $N^{C_0}$-Lipschitz in $\mathbb{D}$. In other words,
\begin{equation*}
    | \Phi(\mathbf{z}) - \Phi(\hat{\mathbf{z}}) |,
    | \Theta(\mathbf{z}) - \Theta(\hat{\mathbf{z}}) | \leq N^{C_0} \| \mathbf{z} - \hat{\mathbf{z}} \|,
    \quad \text{ where } \quad
    \| \mathbf{z} - \hat{\mathbf{z}} \| = \sum_{k \in \mathcal{J}} |z_k - \hat{z}_k|.
\end{equation*}
On the other hand, by choosing $C_0$ large enough, we can also ensure the $N^{C_0}$-Lipschitz continuity of the mapping
\begin{equation*}
    \mathbf{z} = (z_k)_{k \in \mathcal{J}} \quad \mapsto \quad
    \exp \left(-\frac{\mathcal{V}(\mathbf{z}) + 2 \mathcal{W}(\mathbf{z})}{2}\right).
\end{equation*}
Now, we let the constant $C_1$ in definition of $\mathbb{D}_1$ to be larger than $C_0 + 1$. On the ``good event'' $\Omega$, we have by $N^{C_0}$-Lipschitz continuity of $\Phi$ and $\Theta$, 
\begin{equation*}
    \big | { \exp [ {\mathrm{i}(\Phi(\mathbf{z}) +\Theta(\mathbf{z}))} ]
    - \exp [ {\mathrm{i}(\Phi(\hat{\mathbf{z}}) +\Theta(\hat{\mathbf{z}}))} ] } \big | \mathbbm{1}_{\Omega}
    \leq 2 N^{C_0} \| \mathbf{z} - \hat{\mathbf{z}} \| 
    \asymp N^{-1}.
\end{equation*}
While on the ``bad event'' $\Omega^c$, we may use the trivial bound $|\mathrm{e}^{\mathrm{i} \zeta}| \leq 1$. Combined with $\mathbb{P}(\Omega^c) \leq N^{-L}$, we therefore obtain
\begin{equation*}
    \big | { \mathbb{E} \exp [ {\mathrm{i}(\Phi(\mathbf{z}) +\Theta(\mathbf{z}))} ]
    - \mathbb{E} \exp [ {\mathrm{i}(\Phi(\hat{\mathbf{z}}) +\Theta(\hat{\mathbf{z}}))} ] } \big |
    \leq N^{-1/2}.
\end{equation*}
Similarly, we have 
\begin{equation*}
    \big | { \mathbb{E} \exp [ {\mathrm{i} \Theta(\mathbf{z})} ]
    - \mathbb{E} \exp [ {\mathrm{i} \Theta(\hat{\mathbf{z}})} ] } \big |
    \leq N^{-1/2}.
\end{equation*}
In conclusion, substituting the spectral parameters $\{ \hat{z}_k \}$ in (\ref{eqn-characteristic-function}) with $\{ z_k \}$ will only result in a small error that can be controlled by $N^{-1/2}$. This concludes our proof of Theorem \ref{thm-characteristic-function}.
\end{proof}

For the remaining parts of this document, we focus on the proof of Proposition \ref{prop-recurrence-moment}. We only present the derivation of (\ref{eqn-recurrence-moment-ell}), as the argument also applies to (\ref{eqn-recurrence-moment-one}). By employing cumulant expansion formula, we reshape the problem into the estimation of certain functionals of the resolvent, which can be attained by leveraging local laws.

\subsection{Cumulant expansion}

To simplify the notation, let us introduce the following linear functional on $(M+N) \times (M+N)$ matrices,
\begin{equation}
    \varphi_k (A) := \alpha_k \mathbf{u}_k^\top A \mathbf{r}_k
    + \beta_k \mathbf{v}_k^\top A \mathbf{w}_k
    + \gamma_k \mathbf{u}_k^\top A \mathbf{w}_k
    + \delta_k \mathbf{v}_k^\top A \mathbf{r}_k.
    \label{def-functional}
\end{equation}
Note that the definition of $\varphi_k$ depends on the coefficients $\alpha_k, \beta_k, \gamma_k, \delta_k$, and the deterministic vectors $\mathbf{u}_k, \mathbf{r}_k \in \mathbb{R}^M$ and $\mathbf{v}_k, \mathbf{w}_k \in \mathbb{R}^N$. Recall the definitions of $\mathcal{Q}_k, \Theta_k$ and $\Phi_k$ in Theorem \ref{thm-characteristic-function}. Using $\varphi_k$ we can express them as
\begin{align}
    \mathcal{Q}_k & = \sqrt{N} \varphi_k (\Upsilon_k) = \sqrt{N} \varphi_k (G_k - \Pi_k), \\
    \Theta_k & = - \sqrt{N} \varphi_k (\Pi_k H_k \Pi_k), 
    \label{def-Theta-k} \\
    \Phi_k & = \mathcal{Q}_k - \Theta_k - \mathcal{L}_k
    = \sqrt{N} \varphi_k (G_k - \Pi_k + \Pi_k H_k \Pi_k) - \mathcal{L}_k, 
    \label{def-Phi-k}
\end{align}
where we used the subscript $k$ for matrices $H, G, \Pi, \Upsilon$ to indicate that the corresponding spectral parameter is given by $z_k$, e.g. $G_k \equiv G(z_k)$. Also recall that we have defined $\Phi = \sum_{k} s_k \Phi_k$ and $\Theta = \sum_{k} t_k \Theta_k$. Using the isotropic law, $x_{i \mu} = \oprec(N^{-1/2})$ that follows from Assumption (\ref{assumption-bounded-moments}), along with the boundedness of the vectors $\underline{\mathbf{u}}_k, \underline{\mathbf{v}}_k, \underline{\mathbf{r}}_k, \underline{\mathbf{w}}_k$, we have the elementary bounds 
\begin{equation*}
    \Phi = \oprec(1)
    \qand
    \Theta = \oprec(1).
\end{equation*}
We will make use of these bounds in our proof without making explicit references to them. As mentioned, our proof of Proposition \ref{prop-recurrence-moment} relies on combining the local laws with the cumulant expansion formula. Specifically, our proof of Proposition \ref{prop-recurrence-moment} begins by applying the cumulant expansion formula to
\begin{align}
\begin{split}
    & - \sqrt{N} \mathbb{E} \bigg [ {\sum_{k} s_k \varphi_k (\Pi_k H_k \Upsilon_k) 
    \Phi^{\ell + 1} \mathrm{e}^{\mathrm{i} \Theta}} \bigg ] \label{eqn-before-expansion} \\
    = & - \sqrt{N} \sum_{i \mu} \mathbb{E} \bigg [ {x_{i \mu} \sum_{k} s_k \sqrt{z_k}
    \varphi_k (\Pi_k \dimu \Upsilon_k)
    \Phi^{\ell + 1} \mathrm{e}^{\mathrm{i} \Theta}} \bigg ] \\
    =: & \ \sqrt{N} \sum_{i \mu} \mathbb{E} \big [ {x_{i \mu} \mathcal{Z}^{i \mu} \Phi^{\ell + 1} \mathrm{e}^{\mathrm{i} \Theta}} \big ],
\end{split}
\end{align}
where we used $H = \sqrt{z} \sum_{i \mu} x_{i \mu} \dimu$ and defined
\begin{equation}
    \mathcal{Z}^{i \mu} := - \sum_{k} s_k \sqrt{z_k} \varphi_k (\Pi_k \dimu \Upsilon_k).
    \label{def-Z-imu}
\end{equation}
Expanding each summand in (\ref{eqn-before-expansion}) using Lemma \ref{lemma-cumulant-expansion} up to order $4$, we obtain
\begin{align*}
    & \ \sqrt{N} \sum_{i \mu} \mathbb{E} \big [ {x_{i \mu} \mathcal{Z}^{i \mu} \Phi^{\ell + 1} \mathrm{e}^{\mathrm{i} \Theta}} \big ] \\
    = & \ \sum_{i \mu} \sum_{1 \leq p + q + r \leq 3} \frac{\kappa_{p + q + r + 1}}{p! q! r!} \frac{1}{N^{(p + q + r)/2}}
    \mathbb{E} \big [ {\partial_{i \mu}^p \mathcal{Z}_{i \mu}
    \cdot \partial_{i \mu}^q \Phi^{\ell + 1} \cdot \partial_{i \mu}^r \mathrm{e}^{\mathrm{i} \Theta}} \big ] + \mathcal{R}_{i \mu} \\
    =: & \ \sum_{1 \leq p + q + r \leq 3} \mathbb{E} \mathcal{T}_{p, q, r} + \mathcal{R},
\end{align*}
where we used $\kappa_1 = 0$ and introduced the notation
\begin{align}
\begin{split}
    \mathcal{T}_{p, q, r}
    & = \frac{\kappa_{p + q + r + 1}}{p! q! r!} \frac{1}{N^{(p + q + r)/2}}
    \sum_{i \mu} \partial_{i \mu}^p \mathcal{Z}_{i \mu}
    \cdot \partial_{i \mu}^q \Phi^{\ell + 1} \cdot \partial_{i \mu}^r \mathrm{e}^{\mathrm{i} \Theta} \\
    & =: \frac{\kappa_{p + q + r + 1}}{p! q! r!} \frac{1}{N^{(p + q + r)/2}}
    \sum_{i \mu} \mathcal{D}_{p, q, r}^{i \mu}.    
\end{split} \label{def-T-pqr}
\end{align}
Regarding the remainder term $\mathcal{R} = \sum_{i \mu} \mathcal{R}_{i \mu}$, we can utilize the estimate $x_{i \mu} = \oprec (N^{-1/2})$ and Lemma \ref{lemma-cumulant-expansion} to deduce that, for arbitrary small $\varepsilon > 0$ and large $L > 0$,
\begin{equation}
    | \mathcal{R}^{i \mu} | 
    \leq \sum_{p + q + r = 4} \left [ 
    N^{-2} \sup_{|t| \leq N^{-1/2 + \varepsilon}} 
    \big |\mathcal{D}_{p, q, r}^{i \mu} (t) \big |
    + N^{-L} \sup_{t \in \mathbb{R}} 
    \big |\mathcal{D}_{p, q, r}^{i \mu} (t) \big | \right ],
    \label{eqn-error-imu-bound}
\end{equation}
where $\mathcal{D}_{p, q, r}^{i \mu} (t)$ is constructed based on the matrix $X^{i \mu}(t)$, which is obtained from $X$ by replacing the $(i, \mu)$-th entry with the deterministic quantity $t$. Now, to associate the l.h.s. of (\ref{eqn-before-expansion}) with $\mathbb{E} \big [ {\Phi^{\ell + 2} \mathrm{e}^{\mathrm{i} \Theta}} \big ]$, we utilize the following estimate concerning $\mathcal{T}_{1, 0, 0}$.

\begin{proposition}
\label{prop-self-correction}
The following holds uniformly in $\{ z_k \} \subset \mathbb{D}$,
\begin{equation}
    \mathcal{T}_{1, 0, 0} = - \sqrt{N} \sum_k s_k \varphi_k (G_k + z_k \Pi_k G_k)  
    \Phi^{\ell + 1} \mathrm{e}^{\mathrm{i} \Theta} + \oprec (N^{-1/2}).
    \label{eqn-term-100}
\end{equation}
\end{proposition}

\begin{remark}
As mentioned in Section \ref{subsec-sesquilinear-forms}, for $z \in \mathbb{D}_1$, we have the trivial deterministic bound $\| G \| \leq N^{C_1}$. Hence, by Lemma \ref{lemma-expectation-domination}, we can take the expectation of all explicit terms in (\ref{eqn-term-100}) without affecting its validity when $\{ z_k \} \subset \mathbb{D}_1$. In conclusion, from Proposition \ref{prop-self-correction} we can deduce that, when $\{ z_k \} \subset \mathbb{D}_1$,
\begin{equation*}
    \mathbb{E} \mathcal{T}_{1, 0, 0} 
    = - \sqrt{N} \mathbb{E} \bigg [ {\sum_k s_k \varphi_k (G_k + z_k \Pi_k G_k) \Phi^{\ell + 1} \mathrm{e}^{\mathrm{i} \Theta} } \bigg ] 
    + \oprec (N^{-1/2}).
\end{equation*}
\end{remark}

Now, we can add this estimation of $\mathbb{E} \mathcal{T}_{1, 0, 0}$ to the l.h.s. of (\ref{eqn-before-expansion}). This, together with the resolvent identity $HG = I + zG$, yields
\begin{align*}
    & \ \sqrt{N} \mathbb{E} \bigg [ {\sum_{k} s_k \big [ 
    {-\varphi_k (\Pi_k + z_k \Pi_k G_k - \Pi_k H_k \Pi_k)
    + \varphi_k (G_k + z_k \Pi_k G_k)} \big ]
    \Phi^{\ell + 1} \mathrm{e}^{\mathrm{i} \Theta}} \bigg ] \\
    = & \ \sqrt{N} \mathbb{E} \bigg [ {\sum_{k} s_k 
    \varphi_k (G_k - \Pi_k + \Pi_k H_k \Pi_k)
    \Phi^{\ell + 1} \mathrm{e}^{\mathrm{i} \Theta}} \bigg ] \\
    = & \ \mathbb{E} \bigg [ {\sum_{k} s_k (\Phi_k + \mathcal{L}_k)
    \Phi^{\ell + 1} \mathrm{e}^{\mathrm{i} \Theta}} \bigg ],
\end{align*}
where we recall the definition of the $\Phi_k$'s in (\ref{def-Phi-k}) for the last step. In other words, for $\{ z_k \} \subset \mathbb{D}_1$, we have obtained
\begin{align}
\begin{split}
    \mathbb{E} \big [ {\Phi^{\ell + 2} \mathrm{e}^{\mathrm{i} \Theta}} \big ]
    & = \mathbb{E} \bigg [ {\sum_{k} s_k \Phi_k
    \Phi^{\ell + 1} \mathrm{e}^{\mathrm{i} \Theta}} \bigg ] \\
    & = \sum_{\substack{1 \leq p + q + r \leq 3 \\ (p, q, r) \not= (1, 0, 0)}} \mathbb{E} \mathcal{T}_{p, q, r}
    - \sum_{k} s_k \mathcal{L}_k \mathbb{E} \big [ {\Phi^{\ell + 1} \mathrm{e}^{\mathrm{i} \Theta}} \big ]
    + \mathcal{R}
    + \oprec (N^{-1/2}).
\end{split} \label{eqn-after-expansion}
\end{align}
Our next step is to estimate the the r.h.s. of (\ref{eqn-after-expansion}) using local laws. Specifically, we have the following proposition regarding the estimation for the $\mathcal{T}_{p, q, r}$'s.

\begin{proposition} 
\label{prop-estimation-terms}
The following holds uniformly in $\{ z_k \} \subset \mathbb{D}$.
\begin{enumerate}
    \item For the term of index $(2, 0, 0)$, we have
    \begin{equation}
        \mathcal{T}_{2, 0, 0} 
        = \sum_{k} s_k \mathcal{L}_k \Phi^{\ell + 1} \mathrm{e}^{\mathrm{i} \Theta}
        + \oprec (N^{-1/2}).
        \label{eqn-estimate-200}
    \end{equation} 
    \item For the term of index $(0, 1, 0)$, we have
    \begin{equation}
        \mathcal{T}_{0, 1, 0}
        = (\ell + 1) \sum_{k j} s_k s_j \mathcal{V}_{kj}^{(0, 1, 0)} \Phi^{\ell} \mathrm{e}^{\mathrm{i} \Theta} 
        + \oprec (N^{-1/2}) 
        \label{eqn-estimate-010}
    \end{equation}
    \item For the term of index $(1, 2, 0)$, we have
    \begin{equation}
        \mathcal{T}_{1, 2, 0}
        = (\ell + 1) \sum_{k j} s_k s_j \mathcal{V}_{kj}^{(1, 2, 0)} \Phi^{\ell} \mathrm{e}^{\mathrm{i} \Theta}
        + \oprec (N^{-1/2}) 
        \label{eqn-estimate-120}
    \end{equation}
    \item For the term of index $(1, 0, 1)$, we have
    \begin{equation}
        \mathcal{T}_{1, 0, 1}
        = \mathrm{i} \sum_{k j} s_k t_j \mathcal{W}_{kj} \Phi^{\ell + 1} \mathrm{e}^{\mathrm{i} \Theta} + \oprec (N^{-1/2}).
        \label{eqn-estimate-101}
    \end{equation}
    \item If the index $(p, q, r)$ satisfies $1 \leq p + q + r \leq 3$ but
    \begin{equation*}
        (p, q, r) \notin \{ (1, 0, 0), (2, 0, 0), (0, 1, 0), (1, 2, 0), (1, 0, 1) \},
    \end{equation*}
    then the corresponding term is negligible in the sense that
    \begin{equation*}
        \mathcal{T}_{p, q, r} = \oprec (N^{-1/2}).
    \end{equation*}
\end{enumerate}
\end{proposition}

Again, by the remark following Proposition \ref{prop-self-correction}, one deduce the estimates of $\mathbb{E} \mathcal{T}_{p, q, r}$ from Proposition \ref{prop-estimation-terms} when $\{ z_k \} \subset \mathbb{D}_1$. Finally, for the remainder terms $\mathcal{R}$, we have the following proposition.

\begin{proposition}
\label{prop-estimation-error}
We have the following uniformly in $\{ z_k \} \subset \mathbb{D}_1$.
\begin{equation*}
    \mathcal{R} = \oprec (N^{-1/2}).
\end{equation*}
\end{proposition}

We can now conclude the proof of Proposition \ref{prop-recurrence-moment} by applying Propositions \ref{prop-estimation-terms} and \ref{prop-estimation-error} to r.h.s. of (\ref{eqn-after-expansion}). Consequently, our remaining objective is to prove Propositions \ref{prop-self-correction} to \ref{prop-estimation-error}. To initiate this process, let us start by introducing some technical lemmas that will be recurrently employed throughout the proof.

\subsection{Technical lemmas}
\label{subsec-tech-lemmas}

In accordance with the definition of the $\mathcal{T}_{p,q,r}$'s in (\ref{def-T-pqr}), the following derivative formulas are necessary for the proof. The first two derivatives of $\mathcal{Z}^{i \mu}$ w.r.t. $x_{i \mu}$ are given by
\begin{align*}
    \partial_{i \mu} \mathcal{Z}^{i \mu} 
    & = \sum_{k} s_k z_k \varphi_k (\Pi_k \dimu G_k \dimu G_k), \\
    \partial_{i \mu}^2 \mathcal{Z}^{i \mu} 
    & = -2 \sum_{k} s_k z_k^{3/2} \varphi_k (\Pi_k \dimu G_k \dimu G_k \dimu G_k).
\end{align*}
Recalling (\ref{def-Phi-k}), the first two derivatives of $\Phi = \sum_k s_k \Phi_k$ w.r.t. $x_{i \mu}$ are given by
\begin{align*}
    \partial_{i \mu} \Phi 
    & = - \sqrt{N} \sum_{k} s_k \sqrt{z_k} \varphi_k (G_k \dimu G_k - \Pi_k \dimu \Pi_k), \\
    \partial_{i \mu}^2 \Phi 
    & = 2 \sqrt{N} \sum_{k} s_k z_k \varphi_k (G_k \dimu G_k \dimu G_k).
\end{align*}
As for the nonuniversal component $\Theta = \sum_k t_k \Theta_k$ where the $\Theta_k$'s are given in (\ref{def-Theta-k}), we have $\partial_{i \mu}^p \Theta = 0$ for $p \geq 2$ and
\begin{equation*}
    \partial_{i \mu} \Theta 
    = - \sqrt{N} \sum_{k} t_k \sqrt{z_k} \varphi_k (\Pi_k \dimu \Pi_k).
\end{equation*}
The following two matrices will recurrently emerge in the proof,
\begin{equation}
    \underline{\Sigma}_M := \begin{bmatrix} \Sigma & 0 \\ 0 & 0 \end{bmatrix}
    = \sum_{i} \underline{\Sigma}^{1/2} \mathbf{e}_i 
    \mathbf{e}_i^\top \underline{\Sigma}^{1/2} 
    \qand
    \underline{I}_N := \begin{bmatrix} 0 & 0 \\ 0 & I \end{bmatrix}
    = \sum_{\mu} \underline{\Sigma}^{1/2} \mathbf{e}_\mu 
    \mathbf{e}_\mu^\top \underline{\Sigma}^{1/2}.
    \label{def-udl-sigma-I}
\end{equation}
As mentioned, by the rigidity estimate (\ref{eqn-eigenvalue-rigidity}), we have $\| H \|, \| G \| \prec 1$ for $z \in \mathbb{D}$. Consequently, we have following elementary but useful bounds. Let $z \in \mathbb{D}$ and $\bfk{u}, \bfk{v} \in \mathbb{R}^{M + N}$ be deterministic with $\| \bfk{u} \|, \| \bfk{v} \| \lesssim 1$. Then,
\begin{equation}
    \sum_{i} (\underline{\Sigma}^{1/2} G \bfk{u})_i^2 
    = \sum_{i} \bfk{u}^\top G \underline{\Sigma}^{1/2} \mathbf{e}_i 
    \mathbf{e}_i^\top \underline{\Sigma}^{1/2} G \bfk{u}
    = \bfk{u}^\top G \underline{\Sigma}_M G \bfk{u} = \oprec (1).
\end{equation}
Consequently, by the Cauchy–Schwarz inequality, we have
\begin{equation}
    \sum_{i} \big | {(\underline{\Sigma}^{1/2} G \bfk{u})_i} \big |
    \leq \sqrt{N} \bigg [ {\sum_{i} (\underline{\Sigma}^{1/2} G \bfk{u})_i^2} \bigg ]^{1/2}
    = \oprec (N^{1/2}).
    \label{eqn-cauchy-2}
\end{equation}
Using $\bfk{v}^\top \Pi \underline{\Sigma}_M \Pi \bfk{v} \lesssim 1$ and the Cauchy–Schwarz inequality again, we have
\begin{equation}
    \sum_{i} \big | {(\underline{\Sigma}^{1/2} G \bfk{u})_i (\underline{\Sigma}^{1/2} G \bfk{v})_i} \big | = \oprec (1),
    \quad
    \sum_{i} \big | {(\underline{\Sigma}^{1/2} G \bfk{u})_i (\underline{\Sigma}^{1/2} \Pi \bfk{v})_i} \big | = \oprec (1).
    \label{eqn-cauchy-3}
\end{equation}
Here we point out that in the above bounds the two spectral parameters of $G$ are allowed to be different. These bounds remain valid when considering the summation over $\mu \in \llbracket N \rrbracket$ instead of $i \in \llbracket M \rrbracket$. We usually utilize these bounds without explicitly referring to them. The next lemma is a straightforward consequence of the local law, and it is useful in controlling negligible terms that arise in the proof.

\begin{lemma}
\label{lemma-diag-vector-sum}
We have the following estimates uniformly in $z, z_{\star} \in \mathbb{D}$ and deterministic $\bfk{u} \in \mathbb{R}^{M}, \bfk{v} \in \mathbb{R}^{N}$ with $\| \bfk{u} \|, \| \bfk{v} \| \lesssim 1$,
\begin{align*}
    & \sum_{i} \underline{G}_{i i}(z) (\underline{\Sigma}^{1/2} G (z_{\star}) \bfk{v})_i = \oprec (1),
    & \quad & \sum_{i} \underline{G}_{i i}(z) (\underline{\Sigma}^{1/2} \Upsilon (z_{\star}) \bfk{u})_i = \oprec (1), \\
    & \sum_{\mu} \underline{G}_{\mu \mu}(z) (\underline{\Sigma}^{1/2} G (z_{\star}) \bfk{u})_\mu = \oprec (1),
    & \quad & \sum_{\mu} \underline{G}_{\mu \mu}(z) (\underline{\Sigma}^{1/2} \Upsilon (z_{\star}) \bfk{v})_\mu = \oprec (1).
\end{align*}
\end{lemma}

\begin{proof}[Proof of Lemma \ref{lemma-diag-vector-sum}]
We only provide the proof for the first estimate, and the others can be proven in a similar manner.  Consider the decomposition
\begin{align*}
    \sum_{i} \underline{G}_{i i}(z) (\underline{\Sigma}^{1/2} G(z_{\star}) \bfk{v})_i
    = \sum_{i} \underline{\Pi}_{i i}(z) (\underline{\Sigma}^{1/2} G(z_{\star}) \bfk{v})_i
    + \sum_{i} (\underline{G}_{i i}(z) - \underline{\Pi}_{i i}(z)) (\underline{\Sigma}^{1/2} G (z_{\star}) \bfk{v})_i,
\end{align*}
where second sum can be dominated by $\oprec (1)$ using the isotropic law,
\begin{align*}
    & \underline{G}_{i i}(z) - \underline{\Pi}_{i i}(z) = \oprec (N^{-1/2}), \\
    & (\underline{\Sigma}^{1/2} G(z_{\star}) \bfk{v})_i
    = (\underline{\Sigma}^{1/2} \mathbf{e}_i)^\top G(z_{\star}) \bfk{v}
    = (\underline{\Sigma}^{1/2} \mathbf{e}_i)^\top \Upsilon(z_{\star}) \bfk{v}
    = \oprec (N^{-1/2}).    
\end{align*}
Note that we have used $(\underline{\Sigma}^{1/2} \mathbf{e}_i)^\top \Pi(z_{\star}) \bfk{v} = 0$ since $\bfk{v} \in \mathbb{R}^N$. As for the first sum, let $\boldsymbol{\pi}(z) \in \mathbb{R}^M$ denote the diagonal of $\underline{\Pi}_{M}(z)$. Note that $\boldsymbol{\pi}(z)$ is deterministic with $\| \boldsymbol{\pi} (z) \| \lesssim \sqrt{N}$. Consequently, by isotropic law, we have
\begin{equation*}
    \sum_{i} \underline{\Pi}_{i i}(z) (\underline{\Sigma}^{1/2} G(z_{\star}) \bfk{v})_i
    = (\Sigma^{1/2} \boldsymbol{\pi}(z))^\top G (z_{\star}) \bfk{v}
    = (\Sigma^{1/2} \boldsymbol{\pi}(z))^\top \Upsilon (z_{\star}) \bfk{v}
    = \oprec (1),
\end{equation*}
where we also used the boundedness of $\| \Sigma \|$.
\end{proof}

We conclude this subsection with the so-called two-resolvent local law concerning the deterministic approximations of
\begin{equation}
    W_M (z, z_j) := G(z) \underline{\Sigma}_M G(z_j)
    \qand
    W_N (z, z_j) := G (z) \underline{I}_N G (z_j),
    \label{def-double-resolvent}
\end{equation}
which play a crucial role in derivation of (\ref{eqn-estimate-010}).

\begin{proposition}[two-resolvent isotropic law]
\label{prop-double-resolvent}
We have the following estimates uniformly in $z, z_{\star} \in \mathbb{D}$ and deterministic $\bfk{u} \in \mathbb{R}^{M}, \bfk{v} \in \mathbb{R}^{N}$ with $\| \bfk{u} \|, \| \bfk{v} \| \lesssim 1$. 

\begin{enumerate}
    \item For $W_M$, we have
    \begin{align*}
        \bfk{u}^\top W_M (z, z_{\star}) \bfk{u}
        & = \frac{m[z, z_{\star}]}{m(z) m(z_{\star})} \bfk{u}^\top \Pi_M (z) \Sigma \Pi_M (z_{\star}) \bfk{u}
        + \oprec (N^{-1/2}), \\
        \bfk{v}^\top W_M (z, z_{\star}) \bfk{v}
        & = \frac{1}{\sqrt{z z_{\star}}} \left ( \frac{m[z, z_{\star}]}{m(z) m(z_{\star})} - 1 \right ) \bfk{v}^\top \bfk{v}
        + \oprec (N^{-1/2}), \\
        \bfk{u}^\top W_M (z, z_{\star}) \bfk{v} 
        & = \oprec (N^{-1/2}).
    \end{align*}
    \item For $W_N$, we have
    \begin{align*}
        \bfk{u}^\top W_N (z, z_{\star}) \bfk{u}
        & = \sqrt{z z_{\star}} m[z, z_{\star}] \bfk{u}^\top \Pi_M (z) \Sigma \Pi_M (z_{\star}) \bfk{u}
        + \oprec (N^{-1/2}), \\
        \bfk{v}^\top W_N (z, z_{\star}) \bfk{v}
        & = m[z, z_{\star}] \bfk{v}^\top \bfk{v}
        + \oprec (N^{-1/2}), \\
        \bfk{u}^\top W_N (z, z_{\star}) \bfk{v} 
        & = \oprec (N^{-1/2}).
    \end{align*}
\end{enumerate}
\end{proposition}

Importantly, according to Proposition \ref{prop-double-resolvent}, we cannot obtain the deterministic approximation of $W_M(z, z_{\star}) = G(z) \underline{\Sigma}_M G(z_{\star})$ by simply substituting the two resolvents with $\Pi(z)$ and $\Pi(z_{\star})$. We also note that for $\bfk{u} \in \mathbb{R}^{M}$ and $\bfk{v} \in \mathbb{R}^{N}$, one actually has
\begin{align*}
    \bfk{u}^\top W_M(z, z_{\star}) \bfk{u}
    & = \bfk{u}^\top G_M(z) \Sigma G_M(z_{\star}) \bfk{u} \\
    \bfk{v}^\top W_M(z, z_{\star}) \bfk{v}
    & = \frac{1}{\sqrt{z z_{\star}}} \bfk{v}^\top G_N (z) (X^\top \Sigma^2 X) G_N (z_{\star}) \bfk{v}.
\end{align*}
Here the top-left block $G_M(z) \Sigma G_M(z_{\star})$ has been previously encountered in the analysis of the asymptotics of eigenvectors of sample covariance matrices \cite{baiAsymptoticsEigenvectorsLarge2007}, where the authors employed the martingale difference argument to derive its deterministic approximation. To attain the optimal error bound, here we adopt the cumulant expansion approach that was recently developed for establishing the \emph{multi-resolvent local laws} of Wigner matrices \cite{cipolloniEigenstateThermalizationHypothesis2021,cipolloniOptimalMultiresolventLocal2022,cipolloniThermalisationWignerMatrices2022}. The proof is deferred to Section \ref{sec-two-resolvent}, where interested readers can find more details.

\subsection{Term of index (1,0,0)}

In this subsection we focus on the estimation of $\mathcal{T}_{1, 0, 0}$. Recall the definition of $\mathcal{T}_{1, 0, 0}$ in (\ref{def-T-pqr}) and the derivative formula for $\mathcal{Z}_{i \mu}$ provided in Section \ref{subsec-tech-lemmas}. Since $\kappa_2 = 1$, we have
\begin{equation*}
    \mathcal{T}_{1, 0, 0}
    = \frac{1}{N^{1/2}} \sum_{i \mu} 
    \partial_{i \mu} \mathcal{Z}_{i \mu} \cdot \Phi^{\ell + 1} \cdot \mathrm{e}^{\mathrm{i} \Theta}
    = \frac{1}{N^{1/2}} \sum_k s_k z_k \sum_{i \mu} 
    \varphi_k (\Pi_k \dimu G_k \dimu G_k)
    \Phi^{\ell + 1} \mathrm{e}^{\mathrm{i} \Theta}.
\end{equation*}
Compared with (\ref{eqn-term-100}), it suffices to show that, for each $k$, we have
\begin{equation}
    \frac{1}{N^{1/2}} \sum_{i \mu} \varphi_k (\Pi_k \dimu G_k \dimu G_k) 
    = - \sqrt{N} \varphi_k ( {(z_k^{-1} + \Pi_k) G_k} ) + \oprec (N^{-1/2}).
    \label{eqn-summation-100}
\end{equation}
Let us omit the subscript $k$ for the remaining parts of this subsection, as the argument applies to all $k$. Utilizing the self-consistent equation (\ref{eqn-self-consistent}), we have
\begin{align*}
    \frac{1}{z} + \Pi_M 
    & = \frac{1}{z} - \frac{1}{z(I + m \Sigma)}
    = - m \Sigma \Pi_M
    = - \bigg ( {\frac{1}{N} \sum_{\mu} \underline{\Pi}_{\mu \mu}} \bigg ) \Sigma \Pi_M \\
    \frac{1}{z} + \Pi_N 
    & = \left ( \frac{1}{z} + m \right ) I
    = \frac{1 + zm}{zm} \cdot m I
    = - \bigg ( {\frac{1}{N} \sum_{i} \underline{\Pi}_{i i}} \bigg ) \Pi_N.
\end{align*}
Recall the definition of the linear functional $\varphi$ in (\ref{def-functional}). Let us denote the summation in l.h.s. of (\ref{eqn-summation-100}) by $\mathcal{S}$ and decompose it according to the coefficients $\alpha, \beta, \gamma, \delta$ as follows,
\begin{equation*}
    \mathcal{S}
    = \frac{1}{N^{1/2}} \sum_{i \mu} \varphi (\Pi \dimu G \dimu G)
    =: \alpha \mathcal{S}^{\alpha} + \beta \mathcal{S}^{\beta} + \gamma \mathcal{S}^{\gamma} + \delta \mathcal{S}^{\delta}.
\end{equation*}
It is worth noting that there is an abuse of notation as the definition of the summation $\mathcal{S}$ varies from one subsection to another. To avoid any ambiguity, we could introduce notation like $\mathcal{S}_{p,q,r}$ to explicitly link it to $\mathcal{T}_{p,q,r}$. However, this might become cumbersome and lengthy, and thus we adopt the convention of providing the specific definition of $\mathcal{S}$ at the beginning of each subsequent subsection, restricting its applicability only to that particular subsection.

Let us begin with the $\alpha$-component
\begin{align*}
    \mathcal{S}^\alpha
    & = \frac{1}{N^{1/2}} \sum_{i \mu} \mathbf{u}^\top \Pi \dimu G \dimu G \mathbf{r} \\
    & = \frac{1}{N^{1/2}} \sum_{i \mu}
    (\underline{\Sigma}^{1/2} \Pi \mathbf{u})_i \underline{G}_{\mu \mu} (\underline{\Sigma}^{1/2} G \mathbf{r})_i
    + \frac{1}{N^{1/2}} \sum_{i \mu} (\underline{\Sigma}^{1/2} \Pi \mathbf{u})_i \underline{G}_{\mu i} (\underline{\Sigma}^{1/2} G \mathbf{r})_\mu \\
    & =: \mathcal{S}_1^\alpha + \mathcal{S}_2^\alpha,
\end{align*}
where we used the definition of $\dimu$ in (\ref{def-underline-delta}) and the fact that $(\underline{\Sigma}^{1/2} \Pi \mathbf{u})_\mu = 0$. For the first sum, we can activate the average law (\ref{eqn-average-law}) to get
\begin{align*}
    \mathcal{S}_1^\alpha 
    & = \frac{1}{N^{1/2}} \sum_{i} (\underline{\Sigma}^{1/2} \Pi \mathbf{u})_i (\underline{\Sigma}^{1/2} G \mathbf{r})_i \sum_{\mu} \underline{\Pi}_{\mu \mu}
    + \oprec (N^{-1/2}) \\
    & = \frac{1}{N^{1/2}} \sum_{\mu} \underline{\Pi}_{\mu \mu} \cdot (\Sigma \Pi_M \mathbf{u})^\top G \mathbf{r} + \oprec (N^{-1/2}) \\
    & = - \sqrt{N} [ {(z^{-1} + \Pi_M) \mathbf{u}} ]^\top G \mathbf{r} + \oprec (N^{-1/2})
    = - \sqrt{N} \mathbf{u}^\top (z^{-1} + \Pi) G \mathbf{r} + \oprec (N^{-1/2}),
\end{align*}
where the error term in first line can be dominated by $\oprec (N^{-1/2})$ as
\begin{equation*}
    \sum_{\mu} (\underline{G}_{\mu \mu} - \underline{\Pi}_{\mu \mu}) = \oprec (1)
    \qand
    \sum_{i} (\underline{\Sigma}^{1/2} \Pi \mathbf{u})_i (\underline{\Sigma}^{1/2} G \mathbf{r})_i = \oprec (1).
\end{equation*}
Here we note that a further approximation of $\sqrt{N} \mathbf{u}^\top (z^{-1} + \Pi) G \mathbf{r}$ as $\sqrt{N} \mathbf{u}^\top (z^{-1} + \Pi) \Pi \mathbf{r}$ is not applicable. Actually, according to the isotropic law, this approximation will induce an error of order $\oprec (1)$ due to the presence of the prefactor $\sqrt{N}$. We now turn to the second sum $\mathcal{S}_2^\alpha$. By noting $\underline{G} = \underline{G}^\top$ and recalling the definition of $\underline{I}_N$ in (\ref{def-udl-sigma-I}), we have
\begin{align*}
    \mathcal{S}_2^\alpha
    & = \frac{1}{N^{1/2}} \sum_{i \mu} \mathbf{u}^\top \Pi \underline{\Sigma}^{1/2} \mathbf{e}_i \mathbf{e}_i^\top \underline{\Sigma}^{1/2} G \underline{\Sigma}^{1/2}
    \mathbf{e}_\mu \mathbf{e}_\mu^\top \underline{\Sigma}^{1/2} G \mathbf{r} \\
    & = \frac{1}{N^{1/2}} (\Sigma \Pi_M \mathbf{u})^\top 
    G \underline{I}_N G \mathbf{r}
    = \oprec ( N^{-1/2} ),
\end{align*}
where in the last step we used the bounds $\| G \| \prec 1$ and $\| \mathbf{u} \|, \| \mathbf{r} \|, \| \Sigma \|, \| \Pi \| \lesssim 1$. In conclusion, we have proved that
\begin{equation*}
    \mathcal{S}^\alpha = - \sqrt{N} \mathbf{u}^\top (z^{-1} + \Pi) G \mathbf{r} + \oprec (N^{-1/2}).
\end{equation*}
This concludes our proof for the $\alpha$-component, and the same argument applies to the other components. For the $\beta$-component $\mathcal{S}^\beta$, we have
\begin{align*}
    \mathcal{S}^\beta
    & = \frac{1}{N^{1/2}} \sum_{i \mu}
    (\underline{\Sigma}^{1/2} \Pi \mathbf{v})_\mu \underline{G}_{i i} (\underline{\Sigma}^{1/2} G \mathbf{w})_\mu
    + \frac{1}{N^{1/2}} \sum_{i \mu} (\underline{\Sigma}^{1/2} \Pi \mathbf{v})_\mu \underline{G}_{i \mu} (\underline{\Sigma}^{1/2} G \mathbf{w})_i \\
    & = \frac{1}{N^{1/2}} \sum_{i} \underline{\Pi}_{i i} \cdot (\Pi_N \mathbf{v})^\top G \mathbf{w} 
    + \frac{1}{N^{1/2}} (\Pi_N \mathbf{v})^\top G \underline{\Sigma}_M G \mathbf{w} 
    + \oprec (N^{-1/2}) \\
    & = - \sqrt{N} \mathbf{v}^\top (z^{-1} + \Pi) G \mathbf{w} + \oprec (N^{-1/2}).
\end{align*}
For the $\gamma$-component $\mathcal{S}^\gamma$, we can get
\begin{align*}
    \mathcal{S}^\gamma
    & = \frac{1}{N^{1/2}} \sum_{i \mu}
    (\underline{\Sigma}^{1/2} \Pi \mathbf{u})_i \underline{G}_{\mu \mu} (\underline{\Sigma}^{1/2} G \mathbf{w})_i
    + \frac{1}{N^{1/2}} \sum_{i \mu} (\underline{\Sigma}^{1/2} \Pi \mathbf{u})_i \underline{G}_{\mu i} (\underline{\Sigma}^{1/2} G \mathbf{w})_\mu \\
    & = \frac{1}{N^{1/2}} \sum_{i} \underline{\Pi}_{\mu \mu} \cdot (\Sigma \Pi_M \mathbf{u})^\top G \mathbf{w} 
    + \frac{1}{N^{1/2}} (\Sigma \Pi_M \mathbf{u})^\top G \underline{I}_N G \mathbf{w} 
    + \oprec (N^{-1/2}) \\
    & = - \sqrt{N} \mathbf{u}^\top (z^{-1} + \Pi) G \mathbf{w} + \oprec (N^{-1/2}).
\end{align*}
Finally, for the $\delta$-component $\mathcal{S}^\delta$, we have
\begin{align*}
    \mathcal{S}^\delta
    & = \frac{1}{N^{1/2}} \sum_{i \mu}
    (\underline{\Sigma}^{1/2} \Pi \mathbf{v})_\mu \underline{G}_{i i} (\underline{\Sigma}^{1/2} G \mathbf{r})_\mu
    + \frac{1}{N^{1/2}} \sum_{i \mu} (\underline{\Sigma}^{1/2} \Pi \mathbf{v})_\mu \underline{G}_{i \mu} (\underline{\Sigma}^{1/2} G \mathbf{r})_i \\
    & = \frac{1}{N^{1/2}} \sum_{i} \underline{\Pi}_{i i} \cdot (\Pi_N \mathbf{v})^\top G \mathbf{r} 
    + \frac{1}{N^{1/2}} (\Pi_N \mathbf{v})^\top G \underline{\Sigma}_M G \mathbf{r} 
    + \oprec (N^{-1/2}) \\
    & = - \sqrt{N} \mathbf{v}^\top (z^{-1} + \Pi) G \mathbf{r} + \oprec (N^{-1/2}).
\end{align*}
Summarising the estimates above, we conclude the proof of (\ref{eqn-summation-100}) and, consequently, Proposition \ref{prop-self-correction}.

\subsection{Term of index (2,0,0)}

This subsection is devoted to the estimation of $\mathcal{T}_{2, 0, 0}$. Plugging the derivative formula of $\sum_{i \mu} \partial_{i \mu}^2 \mathcal{Z}_{i \mu}$ into the definition (\ref{def-T-pqr}) yields
\begin{align*}
    \mathcal{T}_{2, 0, 0} 
    & = \frac{\kappa_{3}}{2 N} \sum_{i \mu} 
    \partial_{i \mu}^2 \mathcal{Z}_{i \mu}
    \cdot \Phi^{\ell + 1} \cdot \mathrm{e}^{\mathrm{i} \Theta} \\
    & = - \frac{\kappa_{3}}{N} \sum_{k} s_k z_k^{3/2} 
    \sum_{i \mu} \varphi_k (\Pi_k \dimu G_k \dimu G_k \dimu G_k)
    \Phi^{\ell + 1} \mathrm{e}^{\mathrm{i} \Theta}.
\end{align*}
Recall the definition of $\mathcal{L}_k$ in Definition \ref{def-sesquilinear-quantities}. It suffices to show that, for each $k$,
\begin{align*}
    & \ \frac{1}{N} \sum_{i \mu} \varphi_k (\Pi_k \dimu G_k \dimu G_k \dimu G_k) \\
    = & \ \frac{\gamma_k}{N}
    \sum_i \underline{\Pi}_{ii} (z_k) (\underline{\Sigma}^{1/2} \Pi_k \mathbf{u}_k)_i
    \sum_\mu \underline{\Pi}_{\mu \mu} (z_k) (\underline{\Sigma}^{1/2} \Pi_k \mathbf{w}_k)_\mu \\
    & + \frac{\delta_k}{N} \sum_\mu \underline{\Pi}_{\mu \mu} (z_k) (\underline{\Sigma}^{1/2} \Pi_k \mathbf{v}_k)_\mu
    \sum_i \underline{\Pi}_{ii} (z_k) (\underline{\Sigma}^{1/2} \Pi_k \mathbf{r}_k)_i 
    + \oprec (N^{-1/2}).
\end{align*}
Again, let us disregard the subscript $k$ as the argument applies to all $k$. Consider the sum
\begin{equation*}
    \mathcal{S}
    = \frac{1}{N} \sum_{i \mu} \varphi_k (\Pi_k \dimu G_k \dimu G_k \dimu G_k)
    =: \alpha \mathcal{S}^\alpha + \beta \mathcal{S}^\beta + \gamma \mathcal{S}^\gamma + \delta \mathcal{S}^\delta.
\end{equation*}
We now show that $\mathcal{S}^\gamma$ and $\mathcal{S}^\delta$ are the two dominating components, while the other two can be controlled by $\oprec (N^{-1/2})$. By the definition of $\dimu$, we can decompose the $\gamma$-component as follows,
\begin{align*}
    \mathcal{S}^\gamma
    = & \ \frac{1}{N} \sum_{i \mu}
    (\underline{\Sigma}^{1/2} \Pi \mathbf{u})_i \underline{G}_{\mu \mu} 
    \underline{G}_{i i}  (\underline{\Sigma}^{1/2} G \mathbf{w})_\mu \\
    & + \frac{2}{N} \sum_{i \mu}
    (\underline{\Sigma}^{1/2} \Pi \mathbf{u})_i \underline{G}_{\mu \mu} 
    \underline{G}_{\mu i} (\underline{\Sigma}^{1/2} G \mathbf{w})_i 
    + \frac{1}{N} \sum_{i \mu}
    (\underline{\Sigma}^{1/2} \Pi \mathbf{u})_i \underline{G}_{\mu i}
    \underline{G}_{\mu i} (\underline{\Sigma}^{1/2} G \mathbf{w})_\mu \\
    =: & \ \mathcal{S}_1^\gamma + \mathcal{S}_2^\gamma + \mathcal{S}_3^\gamma,
\end{align*}
where we also used $\underline{G}_{i \mu} = \underline{G}_{\mu i}$ and $(\underline{\Sigma}^{1/2} \Pi \mathbf{u})_{\mu} = 0$. For the first sum, we can activate the isotropic law (\ref{eqn-isotropic-law}) to get
\begin{align*}
    \mathcal{S}_1^\gamma 
    & = \frac{1}{N} \sum_{i} \underline{\Pi}_{ii} (\underline{\Sigma}^{1/2} \Pi \mathbf{u})_i
    \sum_{\mu} \underline{\Pi}_{\mu \mu} (\underline{\Sigma}^{1/2} G \mathbf{w})_\mu + \oprec (N^{-1/2}) \\
    & = \frac{1}{N} \sum_{i} \underline{\Pi}_{ii} (\underline{\Sigma}^{1/2} \Pi \mathbf{u})_i
    \sum_{\mu} \underline{\Pi}_{\mu \mu} (\underline{\Sigma}^{1/2} \Pi \mathbf{w})_\mu + \oprec (N^{-1/2}),
\end{align*}
where we implicitly used the elementary bounds of the form
\begin{equation*}
    \Big | {\sum_{\mu} \underline{G}_{\mu \mu} (\underline{\Sigma}^{1/2} G \mathbf{w})_\mu} \Big | 
    \prec \sum_{\mu} \big | {(\underline{\Sigma}^{1/2} G \mathbf{w})_\mu} \big | = \oprec (N^{1/2}).
\end{equation*}
For the second sum, we can use the basic bounds $\underline{G}_{\mu \mu} = \oprec (1)$, $\underline{G}_{\mu i} = \oprec (N^{-1/2})$ to get
\begin{equation*}
    | \mathcal{S}_2^\gamma |
    \prec \frac{1}{N^{1/2}} \sum_{i}
    \big | (\underline{\Sigma}^{1/2} \Pi \mathbf{u})_i (\underline{\Sigma}^{1/2} G \mathbf{w})_i \big |
    = \oprec (N^{-1/2}).
\end{equation*}
Finally, for the third sum, we apply $\underline{G}_{i \mu} = \oprec (N^{-1/2})$ twice to get
\begin{equation*}
    | \mathcal{S}_3^\gamma | 
    \prec \frac{1}{N^{2}} 
    \sum_{i} \big | (\underline{\Sigma}^{1/2} \Pi \mathbf{u})_i \big |
    \sum_{\mu} \big | (\underline{\Sigma}^{1/2} G \mathbf{w})_\mu \big |
    = \oprec (N^{-1}).
\end{equation*}
Summarizing, we have proved that
\begin{equation*}
    \mathcal{S}^\gamma 
    = \frac{1}{N} \sum_{i} \underline{\Pi}_{ii} (\underline{\Sigma}^{1/2} \Pi \mathbf{u})_i
    \sum_{\mu} \underline{\Pi}_{\mu \mu} (\underline{\Sigma}^{1/2} \Pi \mathbf{w})_\mu 
    + \oprec (N^{-1/2}).
\end{equation*}
The same argument applies to the $\delta$-component, which yields
\begin{align*}
    \mathcal{S}^{\delta} 
    & = \frac{1}{N} \sum_{i \mu} (\underline{\Sigma}^{1/2} \Pi \mathbf{v})_\mu \underline{G}_{i i} \underline{G}_{\mu \mu} (\underline{\Sigma}^{1/2} G \mathbf{r})_i + \oprec (N^{-1/2}) \\
    & = \frac{1}{N} \sum_{i} \underline{\Pi}_{i i} (\underline{\Sigma}^{1/2} \Pi \mathbf{r})_i 
    \sum_{\mu} \underline{\Pi}_{\mu \mu} (\underline{\Sigma}^{1/2} \Pi \mathbf{v})_\mu + \oprec (N^{-1/2}).
\end{align*}
To complete the proof of (\ref{eqn-estimate-200}), it remains to show $\mathcal{S}^\alpha, \mathcal{S}^\beta = \oprec (N^{-1/2})$. We only present the proof for $\mathcal{S}^\alpha$, since $\mathcal{S}^\beta$ can be controlled in a similar manner. Again, we decompose the sum into three parts,
\begin{align*}
    \mathcal{S}^\alpha
    = & \ \frac{1}{N} \sum_{i \mu}
    (\underline{\Sigma}^{1/2} \Pi \mathbf{u})_i \underline{G}_{\mu \mu} 
    \underline{G}_{i i}  (\underline{\Sigma}^{1/2} G \mathbf{r})_\mu \\
    & + \frac{2}{N} \sum_{i \mu}
    (\underline{\Sigma}^{1/2} \Pi \mathbf{u})_i \underline{G}_{\mu \mu} 
    \underline{G}_{\mu i} (\underline{\Sigma}^{1/2} G \mathbf{r})_i 
    + \frac{1}{N} \sum_{i \mu}
    (\underline{\Sigma}^{1/2} \Pi \mathbf{u})_i \underline{G}_{\mu i}
    \underline{G}_{\mu i} (\underline{\Sigma}^{1/2} G \mathbf{r})_\mu \\
    =: & \ \mathcal{S}_1^\alpha + \mathcal{S}_2^\alpha + \mathcal{S}_3^\alpha.
\end{align*}
Here $\mathcal{S}_2^\alpha$ and $\mathcal{S}_3^\alpha$ can be dominated by $\oprec (N^{-1/2})$ using an argument similar to that applied to $\mathcal{S}_2^\gamma$ and $\mathcal{S}_3^\gamma$. As for $\mathcal{S}_1^\alpha$, we utilize Lemma \ref{lemma-diag-vector-sum} to handle the summation w.r.t. $\mu$,
\begin{equation*}
    \mathcal{S}_1^\alpha
    \prec \frac{1}{N} \Big | {\sum_{i} \underline{G}_{i i} (\underline{\Sigma}^{1/2} \Pi \mathbf{u})_i} \Big | 
    \Big | {\sum_{\mu} \underline{G}_{\mu \mu} (\underline{\Sigma}^{1/2} G \mathbf{r})_\mu} \Big |
    = \oprec (N^{-1/2}).
\end{equation*}
Therefore, we have concluded the proof of (\ref{eqn-estimate-200}).

\subsection{Term of index (0,1,0)}

In the subsequent subsections, the summations we need to deal with are of the form $\mathcal{S} = \sum_{i \mu} \varphi (A^{i \mu}) \varphi (B^{i \mu})$, where $A^{i \mu}, B^{i \mu}$ are $(M+N) \times (M+N)$ matrices, as exemplified in (\ref{eqn-summation-example}). Utilizing the definition of $\varphi$ in (\ref{def-functional}), we may write
\begin{align*}
    \varphi (A^{i \mu}) 
    & = \alpha \xi_{\alpha}^{i\mu} + \beta \xi_{\beta}^{i\mu} + \gamma \xi_{\gamma}^{i\mu} + \delta \xi_{\delta}^{i\mu}, \\
    \varphi (B^{i \mu}) 
    & = \alpha \zeta_{\alpha}^{i\mu} + \beta \zeta_{\beta}^{i\mu} + \gamma \zeta_{\gamma}^{i\mu} + \delta \zeta_{\delta}^{i\mu}.
\end{align*}
which naturally leads to the following decomposition
\begin{equation*}
    \mathcal{S} = \sum_{a, b \in \{ \alpha, \beta, \gamma, \delta \}} ab \mathcal{S}^{a, b},
    \quad \text{ where } \quad
    \mathcal{S}^{a, b} = \sum_{i \mu} \xi_{a}^{i \mu} \zeta_{b}^{i \mu}.
\end{equation*}
In the following subsections, we will adopt the above formulation as the definition of $\mathcal{S}^{a, b}$. For example, if $\mathcal{S}$ is given by (\ref{eqn-summation-example}), then
\begin{equation*}
    \mathcal{S}^{\alpha, \beta} 
    = \sum_{i \mu} \mathbf{u}_k^\top \Pi_k \dimu \Upsilon_k \mathbf{r}_k
    \cdot \mathbf{v}_j^\top ( G_j \dimu G_j - \Pi_j \dimu \Pi_j ) \mathbf{w}_j.
\end{equation*}
Note that the order of the two superscripts in $\mathcal{S}^{a, b}$ matters, that is, the components $\mathcal{S}^{a, b}$ and $\mathcal{S}^{b, a}$ are generally not identical for $a \not= b$. 

This subsection focuses on estimating $\mathcal{T}_{0, 1, 0}$. Substituting the derivative formula for $\partial_{i \mu} \Phi$ into the definition of $\mathcal{T}_{0, 1, 0}$ as given in (\ref{def-T-pqr}), we obtain
\begin{align*}
    \mathcal{T}_{0, 1, 0}
    & = \frac{1}{N^{1/2}} \sum_{i \mu} 
    \mathcal{Z}_{i \mu} \cdot \partial_{i \mu} \Phi^{\ell + 1} \cdot \mathrm{e}^{\mathrm{i} \Theta} \\
    & = (\ell + 1) \sum_{k j} s_k s_j \sqrt{z_k z_j} \sum_{i \mu} 
    \varphi_k (\Pi_k \dimu \Upsilon_k)
    \varphi_j (G_j \dimu G_j - \Pi_j \dimu \Pi_j)
    \Phi^{\ell} \mathrm{e}^{\mathrm{i} \Theta}.
\end{align*}
It turns out that we need to deal with the sum
\begin{equation}
    \mathcal{S} = \sum_{i \mu} 
    \varphi_k (\Pi_k \dimu \Upsilon_k)
    \varphi_j (G_j \dimu G_j - \Pi_j \dimu \Pi_j).
    \label{eqn-summation-example}
\end{equation}
Recall the convention we made in the beginning of this subsection. We now address the estimation of the components $\mathcal{S}^{a, b}$ for $a, b \in \{ \alpha, \beta, \gamma, \delta \}$. Let us commence with the $(\alpha, \alpha)$-component. Note the $\mathbf{u}^\top \Pi \dimu \Pi \mathbf{r}$ automatically vanishes by the definition of $\dimu$ and the fact that $\mathbf{u}, \mathbf{r} \in \mathbb{R}^M$. Hence, we have
\begin{align*}
    \mathcal{S}^{\alpha, \alpha} 
    = & \ \sum_{i \mu} \mathbf{u}_k^\top \Pi_k \dimu \Upsilon_k \mathbf{r}_k
    \cdot \mathbf{u}_j^\top ( G_j \dimu G_j - \Pi_j \dimu \Pi_j ) \mathbf{r}_j \\
    = & \ \sum_{i \mu} 
    (\underline{\Sigma}^{1/2} \Pi_k \mathbf{u}_k)_i (\underline{\Sigma}^{1/2} G_k \mathbf{r}_k)_\mu
    \cdot (\underline{\Sigma}^{1/2} G_j \mathbf{u}_j)_i (\underline{\Sigma}^{1/2} G_j \mathbf{r}_j)_\mu \\
    & + \sum_{i \mu} 
    (\underline{\Sigma}^{1/2} \Pi_k \mathbf{u}_k)_i (\underline{\Sigma}^{1/2} G_k \mathbf{r}_k)_\mu
    \cdot (\underline{\Sigma}^{1/2} G_j \mathbf{u}_j)_\mu (\underline{\Sigma}^{1/2} G_j \mathbf{r}_j)_i \\
    = & \ \mathbf{u}_k^\top \Pi_k \underline{\Sigma}_M G (z_j) \mathbf{u}_j 
    \cdot \mathbf{r}_k^\top W_N (z_k, z_j) \mathbf{r}_j
    + \mathbf{u}_k^\top \Pi_k \underline{\Sigma}_M G (z_j) \mathbf{r}_j
    \cdot \mathbf{r}_k^\top W_N (z_k, z_j) \mathbf{u}_j.
\end{align*}
Here we recall that $\underline{\Sigma}_M$ and $W_N$ are defined in (\ref{def-udl-sigma-I}) and (\ref{def-double-resolvent}), respectively. Hence, by employing the isotropic laws, namely Propositions \ref{prop-local-law} and \ref{prop-double-resolvent}, to the resolvents $G (z_j)$ and $W_N (z_k, z_j)$, we arrive at the estimate
\begin{align*}
    \mathcal{S}^{\alpha, \alpha} 
    = \sqrt{z_k z_j} m[z_k, z_j] 
    \big ( {\underline{\mathbf{u}}_k^\top \underline{\mathbf{u}}_j 
    \cdot \underline{\mathbf{r}}_k^\top \underline{\mathbf{r}}_j 
    + \underline{\mathbf{u}}_k^\top \underline{\mathbf{r}}_j
    \cdot \underline{\mathbf{r}}_k^\top \underline{\mathbf{u}}_j } \big )
    + \oprec (N^{-1/2}),
\end{align*}
where we used the abbreviations $\underline{\mathbf{u}}_k = \Sigma^{1/2} \Pi_M (z_k) \mathbf{u}_k$ and $\underline{\mathbf{r}}_k = \Sigma^{1/2} \Pi_M (z_k) \mathbf{r}_k$. Similarly, for the $(\beta, \beta)$-component, we have
\begin{align*}
    \mathcal{S}^{\beta, \beta} 
    & = \mathbf{v}_k^\top \Pi_k \underline{I}_N G (z_j) \mathbf{v}_j 
    \cdot \mathbf{w}_k^\top W_M (z_k, z_j) \mathbf{w}_j
    + \mathbf{v}_k^\top \Pi_k \underline{I}_N G (z_j) \mathbf{w}_j
    \cdot \mathbf{w}_k^\top W_M (z_k, z_j) \mathbf{v}_j \\
    & = \frac{m [z_k, z_j] - m_k m_j}{\sqrt{z_k z_j} m_k^2 m_j^2}
    \big ( {\underline{\mathbf{v}}_k^\top \underline{\mathbf{v}}_j 
    \cdot \underline{\mathbf{w}}_k^\top \underline{\mathbf{w}}_j
    + \underline{\mathbf{v}}_k^\top \underline{\mathbf{w}}_j 
    \cdot \underline{\mathbf{w}}_k^\top \underline{\mathbf{v}}_j} \big )
    + \oprec (N^{-1/2}),
\end{align*}
where we used $\underline{\mathbf{v}}_k = m_k \mathbf{v}_k$ and $\underline{\mathbf{w}}_k = m_k \mathbf{w}_k$.

We next proceed to estimate the components $\mathcal{S}^{\gamma, \gamma}, \mathcal{S}^{\gamma, \delta}, \mathcal{S}^{\delta, \gamma}, \mathcal{S}^{\delta, \delta}$. Note that for the $\gamma$-component associated with $\varphi_j (G_j \dimu G_j - \Pi_j \dimu \Pi_j)$, we can decompose it into three parts,
\begin{equation}
    (\underline{\Sigma}^{1/2} G \mathbf{u})_i (\underline{\Sigma}^{1/2} G \mathbf{w})_\mu 
    + (\underline{\Sigma}^{1/2} G \mathbf{u})_\mu (\underline{\Sigma}^{1/2} G \mathbf{w})_i
    - (\underline{\Sigma}^{1/2} \Pi \mathbf{u})_i (\underline{\Sigma}^{1/2} \Pi \mathbf{w})_\mu.
    \label{eqn-decomp-gamma-2}
\end{equation}
This readily leads to a decomposition of the $(\gamma, \gamma)$-component of $\mathcal{S}$,
\begin{equation*}
    \mathcal{S}^{\gamma, \gamma} = \mathcal{S}_1^{\gamma, \gamma} + \mathcal{S}_2^{\gamma, \gamma} + \mathcal{S}_3^{\gamma, \gamma}.
\end{equation*}
Here the first sum reads
\begin{align*}
    \mathcal{S}_1^{\gamma, \gamma} 
    & = \sum_{i \mu}
    (\underline{\Sigma}^{1/2} \Pi_k \mathbf{u}_k)_i (\underline{\Sigma}^{1/2} (G_k - \Pi_k) \mathbf{w}_k)_\mu
    \cdot (\underline{\Sigma}^{1/2} G_j \mathbf{u}_j)_i (\underline{\Sigma}^{1/2} G_j \mathbf{w}_j)_\mu \\
    & = \mathbf{u}_k^\top \Pi_k \underline{\Sigma}_M G (z_j) \mathbf{u}_j
    \cdot \big [ \mathbf{w}_k^\top W_N (z_k, z_j) \mathbf{w}_j
    - \mathbf{w}_k^\top \Pi_k \underline{I}_N G(z_j) \mathbf{w}_j \big ].
\end{align*}
Applying the isotropic laws, we obtain
\begin{equation*}
    \mathcal{S}_1^{\gamma, \gamma}
    = \frac{m [z_k, z_j] - m_k m_j}{m_k m_j} \underline{\mathbf{u}}_k^\top \underline{\mathbf{u}}_j
    \cdot \underline{\mathbf{w}}_k^\top \underline{\mathbf{w}}_j
    + \oprec (N^{-1/2}).
\end{equation*}
As for the second sum, we have
\begin{align*}
    \mathcal{S}_2^{\gamma, \gamma} 
    & = \sum_{i \mu}
    (\underline{\Sigma}^{1/2} \Pi_k \mathbf{u}_k)_i (\underline{\Sigma}^{1/2} (G_k - \Pi_k) \mathbf{w}_k)_\mu
    \cdot (\underline{\Sigma}^{1/2} G_j \mathbf{u}_j)_\mu (\underline{\Sigma}^{1/2} G_j \mathbf{w}_j)_i \\
    & = \mathbf{u}_k^\top \Pi_k \underline{\Sigma}_M G (z_j) \mathbf{w}_j
    \cdot \big [ \mathbf{w}_k^\top W_N (z_k, z_j) \mathbf{u}_j
    - \mathbf{w}_k^\top \Pi_k \underline{I}_N G(z_j) \mathbf{u}_j \big ]
    = \oprec (N^{-1}).
\end{align*}
Finally, for the third sum, it can be controlled as follows,
\begin{align*}
    \mathcal{S}_3^{\gamma, \gamma} 
    & = - \sum_{i \mu}
    (\underline{\Sigma}^{1/2} \Pi_k \mathbf{u}_k)_i (\underline{\Sigma}^{1/2} (G_k - \Pi_k) \mathbf{w}_k)_\mu
    \cdot (\underline{\Sigma}^{1/2} \Pi_j \mathbf{u}_j)_i (\underline{\Sigma}^{1/2} \Pi_j \mathbf{w}_j)_\mu \\
    & = - \mathbf{u}_k^\top \Pi_k \underline{\Sigma}_M \Pi_j \mathbf{u}_j
    \cdot \mathbf{w}_k^\top (G_k - \Pi_k) \underline{I}_N \Pi_j \mathbf{w}_j
    = \oprec (N^{-1/2}).
\end{align*}
Up to this point, we have established
\begin{equation*}
    \mathcal{S}^{\gamma, \gamma}
    = \frac{m [z_k, z_j] - m_k m_j}{m_k m_j} \underline{\mathbf{u}}_k^\top \underline{\mathbf{u}}_j
    \cdot \underline{\mathbf{w}}_k^\top \underline{\mathbf{w}}_j
    + \oprec (N^{-1/2}).
\end{equation*}
The estimation of the components $\mathcal{S}^{\gamma, \delta}, \mathcal{S}^{\delta, \gamma},$ and $\mathcal{S}^{\delta, \delta}$ follows a similar approach. The results are summarized as follows,
\begin{align*}
    \mathcal{S}^{\gamma, \delta}
    & = \frac{m [z_k, z_j] - m_k m_j}{m_k m_j} \underline{\mathbf{u}}_k^\top \underline{\mathbf{r}}_j
    \cdot \underline{\mathbf{w}}_k^\top \underline{\mathbf{v}}_j
    + \oprec (N^{-1/2}), \\
    \mathcal{S}^{\delta, \gamma}
    & = \frac{m [z_k, z_j] - m_k m_j}{m_k m_j} \underline{\mathbf{v}}_k^\top \underline{\mathbf{w}}_j
    \cdot \underline{\mathbf{r}}_k^\top \underline{\mathbf{u}}_j
    + \oprec (N^{-1/2}), \\
    \mathcal{S}^{\delta, \delta}
    & = \frac{m [z_k, z_j] - m_k m_j}{m_k m_j} \underline{\mathbf{v}}_k^\top \underline{\mathbf{v}}_j
    \cdot \underline{\mathbf{r}}_k^\top \underline{\mathbf{r}}_j
    + \oprec (N^{-1/2}).
\end{align*}
Let us recall the definition of $\mathcal{V}^{(0,1,0)}_{kj}$ as stated in Definition \ref{def-sesquilinear-quantities}. To conclude the proof of (\ref{eqn-estimate-010}), we still need to demonstrate that all the remaining components of $\mathcal{S}$ not mentioned above can be dominated by $\oprec (N^{-1/2})$. The computation for these components follows a similar approach as described above, although they ultimately do not contribute to $\mathcal{V}^{(0,1,0)}_{kj}$. For example, we have
\begin{align*}
    \mathcal{S}^{\alpha, \beta}
    = & \ \mathbf{u}_k^\top \Pi_k \underline{\Sigma}_M G (z_j) \mathbf{v}_j 
    \cdot \mathbf{r}_k^\top W_N (z_k, z_j) \mathbf{w}_j \\
    & + \mathbf{u}_k^\top \Pi_k \underline{\Sigma}_M G (z_j) \mathbf{w}_j
    \cdot \mathbf{r}_k^\top W_N (z_k, z_j) \mathbf{v}_j
    = \oprec (N^{-1}). \\
    \mathcal{S}^{\gamma, \alpha}
    = & \ \mathbf{u}_k^\top \Pi_k \underline{\Sigma}_M G (z_j) \mathbf{u}_j
    \cdot \big [ \mathbf{w}_k^\top W_N (z_k, z_j) \mathbf{r}_j
    - \mathbf{w}_k^\top \Pi_k \underline{I}_N G(z_j) \mathbf{r}_j \big ] \\
    & + \mathbf{u}_k^\top \Pi_k \underline{\Sigma}_M G (z_j) \mathbf{r}_j
    \cdot \big [ \mathbf{w}_k^\top W_N (z_k, z_j) \mathbf{u}_j
    - \mathbf{w}_k^\top \Pi_k \underline{I}_N G(z_j) \mathbf{u}_j \big ]
    = \oprec (N^{-1/2}).
\end{align*}
Similarly, we have $\mathcal{S}^{\beta, \alpha} = \oprec (N^{-1})$ and $\mathcal{S}^{\gamma, \beta}, \mathcal{S}^{\delta, \alpha}, \mathcal{S}^{\delta, \beta} = \oprec (N^{-1/2})$. For $\mathcal{S}^{\alpha, \gamma}$, we can use the decomposition (\ref{eqn-decomp-gamma-2}) again to write $\mathcal{S}^{\alpha, \gamma} = \mathcal{S}_1^{\alpha, \gamma} + \mathcal{S}_2^{\alpha, \gamma} + \mathcal{S}_3^{\alpha, \gamma}$, where
\begin{align*}
    \mathcal{S}_1^{\alpha, \gamma}
    & = \mathbf{u}_k^\top \Pi_k \underline{\Sigma}_M G (z_j) \mathbf{u}_j
    \cdot \mathbf{r}_k^\top W_N (z_k, z_j) \mathbf{w}_j
    = \oprec (N^{-1/2}), \\
    \mathcal{S}_2^{\alpha, \gamma}
    & = \mathbf{u}_k^\top \Pi_k \underline{\Sigma}_M G (z_j) \mathbf{w}_j
    \cdot \mathbf{r}_k^\top W_N (z_k, z_j) \mathbf{u}_j
    = \oprec (N^{-1/2}), \\
    \mathcal{S}_3^{\alpha, \gamma}
    & = - \mathbf{u}_k^\top \Pi_k \underline{\Sigma}_M \Pi_j \mathbf{u}_j
    \cdot \mathbf{r}_k^\top G (z_k) \underline{I}_N \Pi_j \mathbf{w}_j
    = \oprec (N^{-1/2}).
\end{align*}
A similar argument leads to $\mathcal{S}^{\beta, \gamma}, \mathcal{S}^{\alpha, \delta}, \mathcal{S}^{\alpha, \delta} = \oprec (N^{-1/2})$. Therefore, we have concluded the proof of (\ref{eqn-estimate-010}).

\subsection{Term of index (1,2,0)}

In this subsection, we turn to the estimation of $\mathcal{T}_{1, 2, 0}$. As per our definition, we have
\begin{align*}
    \mathcal{T}_{1, 2, 0}
    & = \frac{\kappa_{4}}{2 N^{3/2}} \sum_{i \mu} 
    \partial_{i \mu} \mathcal{Z}_{i \mu} \cdot \partial_{i \mu}^2 \Phi^{\ell + 1} \cdot \mathrm{e}^{\mathrm{i} \Theta} \\
    & = \frac{\kappa_{4} (\ell + 1)}{N^{3/2}} \sum_{i \mu} 
    \partial_{i \mu} \mathcal{Z}_{i \mu} \cdot \partial_{i \mu}^2 \Phi \cdot \Phi^{\ell} \mathrm{e}^{\mathrm{i} \Theta}
    + \frac{\kappa_{4} (\ell + 1) \ell}{2 N^{3/2}} \sum_{i \mu} 
    \partial_{i \mu} \mathcal{Z}_{i \mu} \cdot (\partial_{i \mu} \Phi)^2 \cdot \Phi^{\ell-1} \mathrm{e}^{\mathrm{i} \Theta}.
\end{align*}
We claim that the second sum above can be controlled by $\oprec (N^{-1/2})$. Actually, using the elementary bound $\| G \| \prec 1$, one can easily deduce 
\begin{equation*}
    \partial_{i \mu} \mathcal{Z}_{i \mu}
    = \sum_{k} s_k z_k \varphi_k (\Pi_k \dimu G_k \dimu G_k)
    = \oprec (1).
\end{equation*}
Consequently, it suffices to examine
\begin{equation*}
    \frac{1}{N^{3/2}} \sum_{i \mu} (\partial_{i \mu} \Phi)^2
    \prec \frac{1}{N^{1/2}} \sum_{k j} \sum_{i \mu} \big | 
    \varphi_k (G_k \dimu G_k \dimu G_k) 
    \varphi_j (G_j \dimu G_j \dimu G_j) \big | .
\end{equation*}
As before, we can expand the summands using the definition of $\varphi$. Given that both $| \mathcal{J} |$ and the number of components constituting $\varphi$ are finite, we only need to focus on the summation over $i$ and $\mu$, where the summands take the form of the product of two $\bfk{u}^\top G \dimu G \dimu G \bfk{v}$'s. Here, the resolvent entries $\underline{G}_{i i}, \underline{G}_{\mu \mu}, \underline{G}_{i \mu}$ can be simply bounded by $\oprec (1)$. Moreover, since there are exactly two factors indexed by a single  $x$ (e.g. $(\underline{\Sigma}^{1/2} G \bfk{u})_{i}$) and two indexed by a single $\mu$ (e.g. $(\underline{\Sigma}^{1/2} G \bfk{u})_{\mu}$), we can invoke the argument mentioned in (\ref{eqn-cauchy-3}) to control the summation by $\oprec (1)$. Therefore, our claim follows. To summarize, we have
\begin{equation*}
    \mathcal{T}_{1, 2, 0}
    = \frac{\kappa_{4} (\ell + 1)}{N} \sum_{kj} s_k s_j z_k z_j
    \sum_{i \mu} \varphi_k \big ( {\Pi_k (\dimu G_k)^2} \big )
    \varphi_j \big ( {G_j (\dimu G_j)^2} \big )
    \Phi^{\ell} \mathrm{e}^{\mathrm{i} \Theta} + \oprec (N^{-1/2}).
\end{equation*}
We now proceed to estimate the sum
\begin{equation*}
    \mathcal{S} = \frac{1}{N} \sum_{i \mu} \varphi_k (\Pi_k \dimu G_k \dimu G_k)
    \varphi_j (G_j \dimu G_j \dimu G_j).
\end{equation*}
The sum $\mathcal{S}$ consists of two non-vanishing components, $\mathcal{S}^{\alpha, \alpha}$ and $\mathcal{S}^{\beta, \beta}$.  Let us estimate $\mathcal{S}^{\alpha, \alpha}$ firstly. By definition, the $\alpha$-component of $\varphi_k (\Pi_k \dimu G_k \dimu G_k)$ consists of two terms,
\begin{equation}
    \underbrace{(\underline{\Sigma}^{1/2} \Pi_k \mathbf{u}_k)_i \underline{G}_{\mu \mu} (z_k) (\underline{\Sigma}^{1/2} G_k \mathbf{r}_k)_i}_{(1)}
    + \underbrace{(\underline{\Sigma}^{1/2} \Pi_k \mathbf{u}_k)_i \underline{G}_{\mu i} (z_k)
    (\underline{\Sigma}^{1/2} G_k \mathbf{r}_k)_\mu}_{(2)}.
    \label{def-comp-120-aa1}
\end{equation}
While the $\alpha$-component of $\varphi_j (G_j \dimu G_j \dimu G_j)$ consists of four terms,
\begin{align}
\begin{split}
    & \ \underbrace{(\underline{\Sigma}^{1/2} G_j \mathbf{u}_j)_i \underline{G}_{\mu \mu} (z_j) (\underline{\Sigma}^{1/2} G_j \mathbf{r}_j)_i}_{(1)} 
    + \underbrace{(\underline{\Sigma}^{1/2} G_j \mathbf{u}_j)_i \underline{G}_{\mu i} (z_j) (\underline{\Sigma}^{1/2} G_j \mathbf{r}_j)_\mu}_{(2)} \\
    + & \ \underbrace{(\underline{\Sigma}^{1/2} G_j \mathbf{u}_j)_\mu \underline{G}_{i \mu} (z_j) (\underline{\Sigma}^{1/2} G_j \mathbf{r}_j)_i}_{(3)} 
    + \underbrace{(\underline{\Sigma}^{1/2} G_j \mathbf{u}_j)_\mu \underline{G}_{i i} (z_j) (\underline{\Sigma}^{1/2} G_j \mathbf{r}_j)_\mu}_{(4)}.
\end{split} \label{def-comp-120-aa2}    
\end{align}
Note that (\ref{def-comp-120-aa1}) and (\ref{def-comp-120-aa2}) naturally lead to the following decomposition of $\mathcal{S}^{\alpha, \alpha}$,
\begin{equation*}
    \mathcal{S}^{\alpha, \alpha} 
    = \mathcal{S}_{1,1}^{\alpha, \alpha} + \mathcal{S}_{1,2}^{\alpha, \alpha} + \mathcal{S}_{1,3}^{\alpha, \alpha} + \mathcal{S}_{1,4}^{\alpha, \alpha} 
    + \mathcal{S}_{2,1}^{\alpha, \alpha} + \mathcal{S}_{2,2}^{\alpha, \alpha} + \mathcal{S}_{2,3}^{\alpha, \alpha} + \mathcal{S}_{2,4}^{\alpha, \alpha} 
\end{equation*}
where the subscripts indicate which terms of (\ref{def-comp-120-aa1}) and (\ref{def-comp-120-aa2}) are used for the summand. For example, we have
\begin{equation*}
    \mathcal{S}^{\alpha, \alpha}_{2, 3} = \frac{1}{N} \sum_{i \mu} 
    (\underline{\Sigma}^{1/2} \Pi_k \mathbf{u}_k)_i \underline{G}_{\mu i} (z_k)
    (\underline{\Sigma}^{1/2} G_k \mathbf{r}_k)_\mu
    \cdot (\underline{\Sigma}^{1/2} G_j \mathbf{u}_j)_\mu \underline{G}_{i \mu} (z_j) (\underline{\Sigma}^{1/2} G_j \mathbf{r}_j)_i.
\end{equation*}
We note that any sum involving the subscripts $2$ or $3$ can be dominated by $\oprec (N^{-1})$ using the isotropic law $\underline{G}_{i \mu} = \underline{G}_{\mu i} = \oprec (N^{-1/2})$. For example, we can use (\ref{eqn-cauchy-2}-\ref{eqn-cauchy-3}) to get
\begin{equation*}
    \mathcal{S}_{1,2}^{\alpha, \alpha} 
    \prec \frac{1}{N^{3/2}} \sum_{i} \big | {
    (\underline{\Sigma}^{1/2} \Pi_k \mathbf{u}_k)_i
    (\underline{\Sigma}^{1/2} G_k \mathbf{r}_k)_i 
    (\underline{\Sigma}^{1/2} G_j \mathbf{u}_j)_i } \big |
    \sum_{\mu} \big | {(\underline{\Sigma}^{1/2} G_j \mathbf{r}_j)_\mu} \big |
    = \oprec (N^{-1}).
\end{equation*}
While for $\mathcal{S}_{1,4}^{\alpha, \alpha}$, one can use $\underline{G}_{i i}, \underline{G}_{\mu \mu} = \oprec (1)$ to get
\begin{equation*}
    \mathcal{S}_{1,4}^{\alpha, \alpha} 
    \prec \frac{1}{N} \sum_{i} \big | {((\underline{\Sigma}^{1/2} \Pi_k \mathbf{u}_k)_i (\underline{\Sigma}^{1/2} G_k \mathbf{r}_k)_i} \big |
    \sum_{\mu} \big | {(\underline{\Sigma}^{1/2} G_j \mathbf{u}_j)_\mu 
    (\underline{\Sigma}^{1/2} G_j \mathbf{r}_j)_\mu} \big |
    = \oprec (N^{-1}).
\end{equation*}
Therefore, by applying the isotropic law, we have
\begin{align*}
    & \ \mathcal{S}^{\alpha, \alpha} 
    = \mathcal{S}_{1,1}^{\alpha, \alpha} + \oprec (N^{-1}) \\
    = & \ \frac{1}{N} \sum_{i} (\underline{\Sigma}^{1/2} \Pi_k \mathbf{u}_k)_i  
    (\underline{\Sigma}^{1/2} G_k \mathbf{r}_k)_i
    (\underline{\Sigma}^{1/2} G_j \mathbf{u}_j)_i 
    (\underline{\Sigma}^{1/2} G_j \mathbf{r}_j)_i
    \sum_{\mu} \underline{G}_{\mu \mu} (z_k) \underline{G}_{\mu \mu} (z_j) + \oprec (N^{-1}) \\
    = & \ \frac{1}{N} \sum_{i} (\underline{\mathbf{u}}_k)_i  
    (\underline{\mathbf{r}}_k)_i
    (\underline{\mathbf{u}}_j)_i  
    (\underline{\mathbf{r}}_j)_i
    \sum_{\mu} \underline{\Pi}_{\mu \mu} (z_k) \underline{\Pi}_{\mu \mu} (z_j) + \oprec (N^{-1/2}).
\end{align*}
A similar argument leads to
\begin{equation*}
    \mathcal{S}^{\beta, \beta} 
    = \frac{1}{N} 
    \sum_{\mu} (\underline{\mathbf{v}}_k)_\mu 
    (\underline{\mathbf{w}}_k)_\mu 
    (\underline{\mathbf{v}}_j)_\mu 
    (\underline{\mathbf{w}}_j)_\mu
    \sum_{i} \underline{\Pi}_{i i} (z_k) \underline{\Pi}_{i i} (z_j)
    + \oprec (N^{-1/2}).
\end{equation*}
Recall our definition of $\mathcal{V}^{(1,2,0)}_{kj}$ in Definition \ref{def-sesquilinear-quantities}. To complete the proof of (\ref{eqn-estimate-120}), we need to demonstrate that all the remaining components of $\mathcal{S}$ can be dominated by $\oprec (N^{-1/2})$. To streamline the presentation, we choose three representative parts of $\mathcal{S}^{\alpha, \gamma}$ to illustrate the idea. Note that the $\gamma$-component of $\varphi (G \dimu G \dimu G)$ consists of four parts,
\begin{align*}
    & \ \underbrace{(\underline{\Sigma}^{1/2} G_j \mathbf{u}_j)_i \underline{G}_{\mu \mu} (z_j) (\underline{\Sigma}^{1/2} G_j \mathbf{w}_j)_i}_{(1)}
    + \underbrace{(\underline{\Sigma}^{1/2} G_j \mathbf{u}_j)_i \underline{G}_{\mu i} (z_j) (\underline{\Sigma}^{1/2} G_j \mathbf{w}_j)_\mu}_{(2)} \\
    + & \ \underbrace{(\underline{\Sigma}^{1/2} G_j \mathbf{u}_j)_\mu \underline{G}_{i \mu} (z_j) (\underline{\Sigma}^{1/2} G_j \mathbf{w}_j)_i}_{(3)}
    + \underbrace{(\underline{\Sigma}^{1/2} G_j \mathbf{u}_j)_\mu \underline{G}_{i i} (z_j) (\underline{\Sigma}^{1/2} G_j \mathbf{w}_j)_\mu}_{(4)}
\end{align*}
Combined with (\ref{def-comp-120-aa1}), we obtain the decomposition 
\begin{equation*}
    \mathcal{S}^{\alpha, \gamma} 
    = \mathcal{S}_{1,1}^{\alpha, \gamma} + \mathcal{S}_{1,2}^{\alpha, \gamma} + \mathcal{S}_{1,3}^{\alpha, \gamma} + \mathcal{S}_{1,4}^{\alpha, \gamma} 
    + \mathcal{S}_{2,1}^{\alpha, \gamma} + \mathcal{S}_{2,2}^{\alpha, \gamma} + \mathcal{S}_{2,3}^{\alpha, \gamma} + \mathcal{S}_{2,4}^{\alpha, \gamma}.
\end{equation*}
It is evident that in any part, each summand precisely involves four factors indexed by a single $i$ or $\mu$, e.g. $(\underline{\Sigma}^{1/2} G_j \mathbf{u}_j)_i$. Let us refer to them as vector entries below to facilitate our discussion. Additionally, each summand also involves two resolvent entries, e.g. $\underline{G}_{\mu \mu} (z_j)$, that have two subscripts. We note that the total number of the $i$'s (or $\mu$'s) as subscript inside each summand must be four. We now examine three distinct cases based on the subscripts of the vector entries.

\textbf{Case 1}: Among the four vector entries, they all share the same subscript $i$ (or $\mu$). We claim that at least one of these entries can be dominated by $\oprec (N^{-1/2})$. In fact, the only scenarios where all the vector entries cannot be dominated by $\oprec (N^{-1/2})$ occur in $\mathcal{S}^{\alpha, \alpha}$ and $\mathcal{S}^{\beta, \beta}$, which are precisely the sums contributing to the variance term. In this case, the resolvent entry $\underline{G}_{\mu \mu}$ (or $\underline{G}_{i i}$, respectively) appear twice, and we can easily dominate them by $\oprec (1)$. For instance, in the summation $\mathcal{S}^{\alpha, \gamma}_{1,1}$ we can utilize $(\underline{\Sigma}^{1/2} G \mathbf{w})_i = \oprec (N^{-1/2})$ to obtain  
\begin{align*}
    \mathcal{S}^{\alpha, \gamma}_{1,1}
    & = \frac{1}{N} \sum_{i \mu} (\underline{\Sigma}^{1/2} \Pi \mathbf{u})_i \underline{G}_{\mu \mu} (\underline{\Sigma}^{1/2} G \mathbf{r})_i
    \cdot (\underline{\Sigma}^{1/2} G \mathbf{u})_i \underline{G}_{\mu \mu} (\underline{\Sigma}^{1/2} G \mathbf{w})_i \\
    & \prec \frac{1}{N^{1/2}} \sum_{i} \big | (\underline{\Sigma}^{1/2} \Pi \mathbf{u})_i (\underline{\Sigma}^{1/2} G \mathbf{r})_i
    (\underline{\Sigma}^{1/2} G \mathbf{u})_i\big |
    = \oprec (N^{-1/2}).
\end{align*}

\textbf{Case 2}: There are three vector entries with subscript $i$ (or $\mu$) and one vector entry indexed by $\mu$ (or $i$, respectively). In this case the resolvent entry $\underline{G}_{\mu i} = \underline{G}_{i \mu}$ appears exactly once, and we can dominate it by $\oprec (N^{-1/2})$. For example,
\begin{align*}
    \mathcal{S}^{\alpha, \gamma}_{1,2}
    & = \frac{1}{N} \sum_{i \mu} (\underline{\Sigma}^{1/2} \Pi \mathbf{u})_i \underline{G}_{\mu \mu}
    (\underline{\Sigma}^{1/2} G \mathbf{r})_i
    \cdot (\underline{\Sigma}^{1/2} G \mathbf{u})_i \underline{G}_{\mu i}
    (\underline{\Sigma}^{1/2} G \mathbf{w})_\mu \\
    & \prec \frac{1}{N^{3/2}} \sum_{i} \big | {(\underline{\Sigma}^{1/2} \Pi \mathbf{u})_i
    (\underline{\Sigma}^{1/2} G \mathbf{r})_i
    (\underline{\Sigma}^{1/2} G \mathbf{u})_i} \big | 
    \sum_{\mu} \big | {(\underline{\Sigma}^{1/2} G \mathbf{w})_\mu} \big | 
    = \oprec (N^{-1}).
\end{align*}

\textbf{Case 3}: There are two vector entries with subscript $i$ and another two with subscript $\mu$. In this case we dominate the two resolvent entries simply by $\oprec (1)$. As for the summation of vector entries, we can easily dominated it by $\oprec (1)$ as demonstrated in (\ref{eqn-cauchy-3}). For example,
\begin{align*}
    \mathcal{S}^{\alpha, \gamma}_{1,4}
    & = \frac{1}{N} \sum_{i \mu} (\underline{\Sigma}^{1/2} \Pi \mathbf{u})_i \underline{G}_{\mu \mu} (\underline{\Sigma}^{1/2} G \mathbf{r})_i
    \cdot (\underline{\Sigma}^{1/2} G \mathbf{u})_\mu \underline{G}_{i i}
    (\underline{\Sigma}^{1/2} G \mathbf{w})_\mu \\
    & \prec \frac{1}{N} \sum_{i} \big | {(\underline{\Sigma}^{1/2} \Pi \mathbf{u})_i (\underline{\Sigma}^{1/2} G \mathbf{r})_i} \big | 
    \sum_{\mu} \big | {(\underline{\Sigma}^{1/2} G \mathbf{u})_\mu (\underline{\Sigma}^{1/2} G \mathbf{w})_\mu} \big | 
    = \oprec (N^{-1}).
\end{align*}
The argument presented above can be extended to all other parts of $\mathcal{S}^{\alpha, \gamma}$ as well as to all the remaining components not explicitly discussed. Due to space limitations, we omit the details for brevity. With this, we have completed the proof of (\ref{eqn-estimate-120}).

\subsection{Term of index (1,0,1)}

We now turn our attention to the final nontrivial component $\mathcal{T}_{1, 0, 1}$. By definition, we can write
\begin{align*}
    \mathcal{T}_{1, 0, 1}
    & = \frac{\kappa_{3}}{N} \sum_{i \mu} \partial_{i \mu} \mathcal{Z}_{i \mu}
    \cdot \Phi^{\ell + 1} \cdot \partial_{i \mu} \mathrm{e}^{\mathrm{i} \Theta} \\
    & = - \frac{\mathrm{i} \kappa_{3}}{N^{1/2}} \sum_{k j} s_k t_j z_k z_j^{1/2} 
    \sum_{i \mu} \varphi_k (\Pi_k \dimu G_k \dimu G_k) 
    \varphi_j (\Pi_j \dimu \Pi_j)
    \Phi^{\ell + 1} \mathrm{e}^{\mathrm{i} \Theta}.
\end{align*}
It turns out that we need to examine the sum
\begin{equation}
    \mathcal{S}
    = \frac{1}{N^{1/2}} \sum_{i \mu} \varphi_k (\Pi_k \dimu G_k \dimu G_k)
    \varphi_j (\Pi_j \dimu \Pi_j).
    \label{eqn-summation-101}
\end{equation}
It is important to note that  we have $\mathcal{S}^{a, \alpha} = 0$ and $\mathcal{S}^{a, \beta} = 0$ for all $a \in \{ \alpha, \beta, \gamma, \delta \}$ since the $\alpha$ and $\beta$-components of $\varphi (\Pi \dimu \Pi)$ automatically vanish,
\begin{equation*}
    \mathbf{u}^\top \Pi \dimu \Pi \mathbf{r} = 0
    \qand 
    \mathbf{v}^\top \Pi \dimu \Pi \mathbf{w} = 0.
\end{equation*}
On the other hand, the $\gamma$-component of $\varphi_j (\Pi_j \dimu \Pi_j)$ is given by
\begin{equation*}
    \mathbf{u}_j^\top \Pi_j \dimu \Pi_j \mathbf{w}_j
    = (\underline{\Sigma}^{1/2} \Pi_j \mathbf{u}_j)_i (\underline{\Sigma}^{1/2} \Pi_j \mathbf{w}_j)_{\mu}.
\end{equation*}
Recalling our decomposition of the $\alpha$-component of $\varphi_k (\Pi_k \dimu G_k \dimu G_k)$ as presented in (\ref{def-comp-120-aa1}), we can thus decompose the $(\alpha, \gamma)$-component of $\mathcal{S}$ as $\mathcal{S}^{\alpha, \gamma} = \mathcal{S}^{\alpha, \gamma}_{1} + \mathcal{S}^{\alpha, \gamma}_2$. For the first part, we can employ the isotropic law to obtain
\begin{align*}
    \mathcal{S}_1^{\alpha, \gamma} 
    & = \frac{1}{N^{1/2}} \sum_{i \mu}
    (\underline{\Sigma}^{1/2} \Pi_k \mathbf{u}_k)_i \underline{G}_{\mu \mu} (z_k) (\underline{\Sigma}^{1/2} G_k \mathbf{r}_k)_i
    \cdot (\underline{\Sigma}^{1/2} \Pi_j \mathbf{u}_j)_i (\underline{\Sigma}^{1/2} \Pi_j \mathbf{w}_j)_{\mu} \\
    & = \frac{1}{N^{1/2}} \sum_{i} (\underline{\mathbf{u}}_k)_i (\underline{\mathbf{r}}_k)_i (\underline{\mathbf{u}}_j)_i
    \sum_{\mu} \underline{\Pi}_{\mu \mu} (z_k) (\underline{\mathbf{w}}_j)_\mu + \oprec (N^{-1/2}).
\end{align*}
On the other hand, using $\underline{G}_{\mu i}, (\underline{\Sigma}^{1/2} G \mathbf{r})_\mu = \oprec (N^{-1/2})$, one can easily deduce
\begin{align*}
    \mathcal{S}_2^{\alpha, \gamma} 
    & = \frac{1}{N^{1/2}} \sum_{i \mu}
    (\underline{\Sigma}^{1/2} \Pi_k \mathbf{u}_k)_i \underline{G}_{\mu i} (z_k) (\underline{\Sigma}^{1/2} G_k \mathbf{r}_k)_\mu
    \cdot (\underline{\Sigma}^{1/2} \Pi_j \mathbf{u}_j)_i (\underline{\Sigma}^{1/2} \Pi_j \mathbf{w}_j)_{\mu} \\
    & \prec \frac{1}{N^{3/2}} 
    \sum_{i} \big | {(\underline{\Sigma}^{1/2} \Pi_k \mathbf{u}_k)_i (\underline{\Sigma}^{1/2} \Pi_j \mathbf{u}_j)_i} \big | 
    \sum_{\mu} \big | {(\underline{\Sigma}^{1/2} \Pi_j \mathbf{w}_j)_{\mu}} \big |
    = \oprec (N^{-1}). 
\end{align*}
Summarizing the above, we have proved
\begin{equation*}
    \mathcal{S}^{\alpha, \gamma} = \frac{1}{N^{1/2}} \sum_{i} (\underline{\mathbf{u}}_k)_i (\underline{\mathbf{r}}_k)_i (\underline{\mathbf{u}}_j)_i
    \sum_{\mu} \underline{\Pi}_{\mu \mu} (z_k) (\underline{\mathbf{w}}_j)_\mu + \oprec (N^{-1/2}).
\end{equation*}
By employing a similar argument, one can also establish that
\begin{align*}
    \mathcal{S}^{\alpha, \delta} & = \frac{1}{N^{1/2}} \sum_{i} (\underline{\mathbf{u}}_k)_i (\underline{\mathbf{r}}_k)_i (\underline{\mathbf{r}}_j)_i
    \sum_{\mu} \underline{\Pi}_{\mu \mu} (z_k) (\underline{\mathbf{v}}_j)_\mu + \oprec (N^{-1/2}), \\
    \mathcal{S}^{\beta, \gamma} & = \frac{1}{N^{1/2}} \sum_{\mu} (\underline{\mathbf{v}}_k)_\mu (\underline{\mathbf{w}}_k)_\mu (\underline{\mathbf{w}}_j)_\mu
    \sum_{i} \underline{\Pi}_{i i} (z_k) (\underline{\mathbf{u}}_j)_i + \oprec (N^{-1/2}), \\
    \mathcal{S}^{\beta, \delta} & = \frac{1}{N^{1/2}} \sum_{\mu} (\underline{\mathbf{v}}_k)_\mu (\underline{\mathbf{w}}_k)_\mu (\underline{\mathbf{v}}_j)_\mu
    \sum_{\mu} \underline{\Pi}_{i i} (z_k) (\underline{\mathbf{r}}_j)_i + \oprec (N^{-1/2}).
\end{align*}
Recall our definition of $\mathcal{W}_{kj}$ in Definition \ref{def-sesquilinear-quantities}. To conclude the proof of (\ref{eqn-estimate-101}), it remains to demonstrate
\begin{equation*}
    \mathcal{S}^{\gamma, \gamma}, \mathcal{S}^{\gamma, \delta}, \mathcal{S}^{\delta, \gamma}, \mathcal{S}^{\delta, \delta} = \oprec (N^{-1/2}).
\end{equation*}
This task is straightforward, and we only provide a brief explanation of how to control $\mathcal{S}^{\gamma, \gamma}$ as an example. Note that the $\gamma$-component of $\varphi_k (\Pi_k \dimu G_k \dimu G_k)$ consists of two parts, explicitly given by
\begin{equation*}
    (\underline{\Sigma}^{1/2} \Pi_k \mathbf{u}_k)_i \underline{G}_{\mu \mu} (z_k) (\underline{\Sigma}^{1/2} G_k \mathbf{w}_k)_i
    + (\underline{\Sigma}^{1/2} \Pi_k \mathbf{u}_k)_i \underline{G}_{\mu i} (z_k)
    (\underline{\Sigma}^{1/2} G_k \mathbf{w}_k)_\mu.
\end{equation*}
Again, this naturally leads to a decomposition $\mathcal{S}^{\gamma, \gamma} = \mathcal{S}^{\gamma, \gamma}_{1} + \mathcal{S}^{\gamma, \gamma}_2$. Applying the isotropic bounds $(\underline{\Sigma}^{1/2} G \mathbf{w})_i, \underline{G}_{\mu i} = \oprec (N^{-1/2})$ to these two parts respectively, we obtain
\begin{equation*}
    \mathcal{S}^{\gamma, \gamma}_{1} = \oprec (N^{-1/2})
    \qand 
    \mathcal{S}^{\gamma, \gamma}_{2} = \oprec (N^{-1}).
\end{equation*}
Therefore, we have finished our proof of (\ref{eqn-estimate-101}).

\subsection{Negligible terms}

To complete the proof of Proposition \ref{prop-estimation-terms}, we still need to prove that $\mathcal{T}_{p, q, r} = \oprec (N^{-1/2})$ for $1 \leq p + q + r \leq 3$ while
\begin{equation*}
    (p, q, r) \notin
    \{ (1, 0, 0), (2, 0, 0), (0, 1, 0), (1, 2, 0), (1, 0, 1) \}.
\end{equation*}
We only present the proof for $\mathcal{T}_{1, 1, 0}$, as the proofs for other terms are similar or even easier. Additionally, we omit the subscripts $k$ and $j$ here, as they are not relevant to the argument. To prove $\mathcal{T}_{1, 1, 0} = \oprec (N^{-1/2})$, it suffices to show that
\begin{equation*}
    \mathcal{S}
    = \frac{1}{N^{1/2}} \sum_{i \mu}
    \varphi (\Pi \dimu G \dimu G)
    \varphi (G \dimu G - \Pi \dimu \Pi)
    = \oprec (N^{-1/2}).
\end{equation*}
Compared with the summation (\ref{eqn-summation-101}) associated with the non-negligible term $\mathcal{T}_{1, 0, 1}$ discussed earlier, the primary distinction lies in the second factor of the summand. In the case of $\mathcal{T}_{1, 0, 1}$, the second factor is given by $\varphi (\Pi \dimu \Pi)$, whose $\gamma$ and $\delta$-components are of constant order. In contrast, for $\mathcal{T}_{1, 1, 0}$, the second factor takes the form of $\varphi (G \dimu G - \Pi \dimu \Pi)$, which can be dominated by $\oprec (N^{-1/2})$, according to the isotropic law (\ref{eqn-isotropic-law}). This allows us to establish the negligibility of $\mathcal{T}_{1, 1, 0}$. 

We take $\mathcal{S}^{\alpha, \gamma}$ as an example, and the other components can be handled in a similar manner. Note that we can decompose the $\gamma$-component of $\varphi (G \dimu G - \Pi \dimu \Pi)$ as follows,
\begin{align*}
    & \ \mathbf{u}^{\top}( G\dimu G-\Pi\dimu\Pi)\mathbf{w} \\
    = & \ \mathbf{u}^{\top}(G-\Pi)\dimu G\mathbf{w}+\mathbf{u}^{\top}\Pi\dimu(G-\Pi)\mathbf{w} \\
    = & \ \underbrace{(\Sigma^{1/2} \Upsilon \mathbf{u})_{i}(\underline{\Sigma}^{1/2}G\mathbf{w})_{\mu}}_{(1)}
    +\underbrace{(\underline{\Sigma}^{1/2} \Upsilon\mathbf{u})_{\mu}(\underline{\Sigma}^{1/2}G\mathbf{w})_{i}}_{(2)}
    + \underbrace{(\underline{\Sigma}^{1/2}\Pi\mathbf{u})_{i}(\underline{\Sigma}^{1/2} \Upsilon \mathbf{w})_{\mu}}_{(3)}.
\end{align*}
Combined with (\ref{def-comp-120-aa1}), we may write 
\begin{equation*}
    \mathcal{S}^{\alpha, \gamma} = \mathcal{S}_{1,1}^{\alpha, \gamma} + \mathcal{S}_{1,2}^{\alpha, \gamma} + \mathcal{S}_{1,3}^{\alpha, \gamma} + \mathcal{S}_{2,1}^{\alpha, \gamma} + \mathcal{S}_{2,2}^{\alpha, \gamma} + \mathcal{S}_{2,3}^{\alpha, \gamma}.
\end{equation*}
For $\mathcal{S}_{1,1}^{\alpha, \gamma}$, we can use the bounds $(\underline{\Sigma}^{1/2} \Upsilon \mathbf{u})_i = \oprec (N^{-1/2})$ and $\underline{G}_{\mu \mu} = \oprec (1)$ to get
\begin{equation*}
    \mathcal{S}_{1,1}^{\alpha, \gamma} = \frac{1}{N}
    \sum_{i} \big | { (\underline{\Sigma}^{1/2} \Pi \mathbf{u})_i (\underline{\Sigma}^{1/2} G \mathbf{r})_i} \big |
    \sum_{\mu} \big | {(\underline{\Sigma}^{1/2} G \mathbf{w})_\mu} \big |
    = \oprec (N^{-1/2}).
\end{equation*}
The same strategy applies for $\mathcal{S}_{1,2}^{\alpha, \gamma}$. While for $\mathcal{S}_{1,3}^{\alpha, \gamma}$, we need to activate Lemma \ref{lemma-diag-vector-sum} to get
\begin{align*}
    \mathcal{S}_{1,3}^{\alpha, \gamma} = \frac{1}{N^{1/2}}
    \sum_{i} (\underline{\Sigma}^{1/2} \Pi \mathbf{u})_i (\underline{\Sigma}^{1/2} G \mathbf{r})_i (\underline{\Sigma}^{1/2} \Pi \mathbf{u})_i 
    \sum_{\mu} \underline{G}_{\mu \mu} (\underline{\Sigma}^{1/2} \Upsilon \mathbf{w})_\mu
    = \oprec (N^{-1/2}).
\end{align*}
As for the sums with the first subscript equal to $2$, the argument is indeed simpler. In conclusion, we have demonstrated that $\mathcal{T}_{1, 1, 0} = \oprec (N^{-1/2})$.

\subsection{Reminder term}
\label{subsec-reminder}

We conclude this section by presenting the proof of Proposition \ref{prop-estimation-error}. Recall that $\mathcal{R} = \sum_{i \mu} \mathcal{R}^{i \mu}$, and according to (\ref{eqn-error-imu-bound}), we have
\begin{equation*}
    | \mathcal{R}^{i \mu} | 
    \leq \sum_{p + q + r = 4} \left [ 
    N^{-2} \sup_{|t| \leq N^{-1/2 + \varepsilon}} 
    \big |\mathcal{D}_{p, q, r}^{i \mu} (t) \big |
    + N^{-L} \sup_{t \in \mathbb{R}} 
    \big |\mathcal{D}_{p, q, r}^{i \mu} (t) \big | \right ].
\end{equation*}
for arbitrary small $\varepsilon > 0$ and arbitrary large $L > 0$. Also recall our definition of $X(t)$ following (\ref{eqn-error-imu-bound}). Let $G(t) \equiv G(t, z)$ be the Green function built on $X(t)$. Then, we have the following resolvent expansion formula
\begin{equation}
    G(t) - G = \sqrt{z} (x_{i \mu} - t) G(t) \Delta^{i \mu} G.
    \label{eqn-resolvent-expansion-error}
\end{equation}
Fix any small $\varepsilon > 0$. Note that with high probability we have $\| G \| \leq N^{\varepsilon}$ and $|x_{i \mu}| \leq N^{-1/2 + \varepsilon}$. For any $|t| \leq N^{-1/2 + \varepsilon}$, we may use (\ref{eqn-resolvent-expansion-error}) twice to deduce firstly $\| G(t) \| \leq N^{2\varepsilon}$ and then $\| G(t) - G \| \leq N^{-1/2 + 5 \varepsilon}$. In particular, this implies
\begin{equation*}
    \big | {\mathcal{D}_{p, q, r}^{i \mu} (t) - \mathcal{D}_{p, q, r}^{i \mu}} \big | \leq N^{-1/2 + C_{\ell} \varepsilon}
\end{equation*}
with high probability for some large constant $C_{\ell} > 0$ depending only on $\ell$. On the other hand, we can employ the argument presented in the previous subsections once again to demonstrate that
\begin{equation*}
    \frac{1}{N^{p+q+r}} \sum_{i \mu} \big | {\mathcal{D}_{p, q, r}^{i \mu}} \big |
    = \oprec (N^{-1/2}),
    \quad \forall p + q + r = 4.
\end{equation*}
Details are omitted here for the sake of simplicity. To summarize, we have obtained
\begin{equation*}
    \sum_{i \mu} \sum_{p + q + r = 4}
    N^{-2} \sup_{|t| \leq N^{-1/2 + \varepsilon}} 
    \big | \mathcal{D}_{p, q, r}^{i \mu} (t) \big |
    = \oprec (N^{-1/2}).
\end{equation*}
Now, it remains to provide a bound for the second term in (\ref{eqn-error-imu-bound}). Note that for Proposition \ref{prop-estimation-error} we have assumed $\{ z_k \} \subset \mathbb{D}_1$. In particular, we can utilize the elementary deterministic bound $\| G(t) \| \leq N^{C_1}$. By examining the definition of $\mathcal{D}_{p, q, r}^{i \mu} (t)$, we obtain the uniform bound
\begin{equation*}
    \big | {\mathcal{D}_{p, q, r}^{i \mu} (t)} \big | 
    \leq N^{\hat{C}_{\ell}},
    \quad \forall t \in \mathbb{R},
\end{equation*}
where the constant $\hat{C}_{\ell} > 0$ depends only on $\ell$ and our choice of $C_1$. On the other hand, the constant $L > 0$ in (\ref{eqn-error-imu-bound}) can be arbitrary large. For example, we may choose $L = \hat{C}_{\ell} + 3$. Then
\begin{equation*}
    \sum_{i \mu} \sum_{p + q + r = 4}
    N^{-L} \sup_{t \in \mathbb{R}} 
    \big | \mathcal{D}_{p, q, r}^{i \mu} (t) \big |
    \lesssim N^{-1}.
\end{equation*}
In conclusion, we have established $\mathcal{R} = \oprec (N^{-1/2})$ for $\{ z_k \} \subset \mathbb{D}_1$.

\section{Double-resolvent isotropic law} 
\label{sec-two-resolvent}

This section is dedicated to the proof of Proposition \ref{prop-double-resolvent}. For simplicity, we only provide the proof for $W \equiv W_M (z, z_{\star})$, as the proof for $W_N (z, z_{\star})$ is much easier. Let $z, z_{\star} \in \mathbb{D}$ be arbitrary spectral parameters. In the following we consistently use the abbreviations
\begin{equation*}
    H \equiv H (z), G \equiv G (z)
    \qand 
    H_{\star} \equiv H (z_{\star}), G_{\star} \equiv G (z_{\star}).
\end{equation*}
Given any $(M+N) \times (M+N)$ matrix $A$, we use the following notion for its blocks,
\begin{equation*}
    A = \begin{bmatrix}
    A_\alpha & A_\gamma \\ A_\delta & A_\beta
    \end{bmatrix},
    \quad \text{ such that } A_\alpha \in \mathbb{R}^{M \times M},
    A_\beta \in \mathbb{R}^{N \times N}.
\end{equation*}
For example, the blocks of $W$ are given by
\begin{align*}
    & W_{\alpha} = G_M \Sigma G_{M \star},
    & \quad & 
    W_{\gamma} = \frac{1}{\sqrt{z_{\star}}} G_M \Sigma Y G_{N \star}, \\
    & W_{\delta} = \frac{1}{\sqrt{z}} G_N Y^\top \Sigma G_{M \star},
    & \quad & 
    W_{\beta} = \frac{1}{\sqrt{z z_{\star}}} G_N Y^\top \Sigma Y G_{N \star}.
\end{align*}
For any square matrix $A$, we use the notation $\langle A \rangle := N^{-1} \operatorname{tr} A$. It is important to note that the normalization factor is consistently $N^{-1}$ regardless of the dimension of matrix $A$.

Our proof of Proposition \ref{prop-double-resolvent} is a standard adaptation of the arguments presented in \cite{cipolloniOptimalMultiresolventLocal2022} to the case of sample covariance matrices. Notably, the analysis here is considerably simplified since we only need to consider product of two resolvents within the spectral domain $\mathbb{D}$, which satisfies $\operatorname{dist} (\mathbb{D}, \operatorname{supp} \varrho) \gtrsim 1$. The only complication is that $\Pi$, the deterministic approximation of $G$, is anisotropic in our case. Nonetheless, the overall approach remains consistent with the techniques developed in \cite{cipolloniOptimalMultiresolventLocal2022}. We provide the details in Section \ref{subsec-double-resol-general}.

On the other hand, if one only aims at the distribution of an individual spiked eigenvalue, then it suffices to prove Proposition \ref{prop-double-resolvent} with $z = z_{\star}$. In this specific case, we can indeed present a more concise and algebraic proof. See Section \ref{subsec-double-resol-same} for more details.

\subsection{General case} 
\label{subsec-double-resol-general}

Recall that our objective is to establish that $W$ can be well approximated by $\Xi := \operatorname{diag} (\Xi_{\alpha}, \Xi_{\beta})$, where
\begin{equation}
    \Xi_{\alpha} := \frac{m [z, z_\star]}{m m_{\star}} \Pi_M \Sigma \Pi_{M \star}
    \qand
    \Xi_{\beta} := \frac{1}{\sqrt{z z_{\star}}} \left ( \frac{m [z, z_\star]}{m m_{\star}} - 1 \right ) I_N.
    \label{def-Xi}
\end{equation}
To begin, let us adapt the notion of \emph{second-order renormalisation} from \cite{cipolloniOptimalMultiresolventLocal2022} to the context of sample covariance matrices. Given a matrix-valued function $g(H, H_{\star})$, we define
\begin{equation}
    \mathcal{E} [H g(H, H_{\star})] 
    := H g(H, H_{\star}) - \frac{\sqrt{z}}{N} \sum_{i \mu} \dimu [\partial_{i\mu} g(H, H_{\star})],
    \label{def-second-renormalization}
\end{equation}
where the summation essentially corresponds to the second-order term in cumulant expansion of $H g(H, H_{\star})$. In particular, if $X$ is Gaussian, then according to Lemma \ref{lemma-cumulant-expansion} we have $\mathbb{E} \mathcal{E} [H g(H, H_{\star})] = 0$. Given any $(M+N) \times (M+N)$ matrix $A$, let us define
\begin{equation}
    \mathcal{S} [A] := \frac{1}{N} \sum_{i\mu} \dimu A \dimu
    = \frac{1}{N} \begin{bmatrix} (\operatorname{tr} A_{\beta}) \Sigma & \Sigma A_{\delta}^\top \\ A_{\gamma}^\top \Sigma & \operatorname{tr} \Sigma A_{\alpha} \end{bmatrix}.
    \label{def-calX}
\end{equation}
We note that the operator $\mathcal{S}$ is linear and can be decomposed into the summation of its diagonal and off-diagonal parts,
\begin{equation*}
    \mathcal{S} [A]
    = \begin{bmatrix} \langle A_{\beta} \rangle \Sigma & 0 \\ 
    0 & \langle A_{\alpha} \Sigma \rangle \end{bmatrix}
    + \frac{1}{N} \begin{bmatrix} 0 & \Sigma A_{\delta}^\top \\ A_{\gamma}^\top \Sigma & 0 \end{bmatrix}
    := \mathcal{S}_{\mathrm{d}} [A] + \mathcal{S}_{\mathrm{o}} [A].
\end{equation*}
Now, by the derivative formula $\partial_{i\mu} G = -\sqrt{z} G \dimu G$, we have
\begin{equation*}
    \mathcal{E} [H G] 
    = H G + \frac{z}{N} \sum_{i\mu} \dimu G \dimu G
    = I + z G + z \mathcal{S} [G] G.
\end{equation*}
On the other hand, using the self-consistent equation (\ref{eqn-self-consistent}), it can be easily verified that
\begin{equation*}
    I + z \Pi + z \mathcal{S} [\Pi] \Pi = 0.
\end{equation*}
This equation is commonly referred to as the \emph{Matrix Dyson Equation} in the literature (see e.g. \cite{erdosMatrixDysonEquation2019}). We now combine the two equations above to obtain
\begin{align}
\begin{split}
    \Pi \mathcal{E} [H G] 
    & = \Pi + z \Pi G + z \Pi \mathcal{S} [G] G \\
    & = \Pi - (I + z \Pi \mathcal{S} [\Pi]) G + z \Pi \mathcal{S} [G] G
    = \Pi - G + z \Pi \mathcal{S} [G - \Pi] G.     
\end{split} \label{eqn-expansion-G}
\end{align}
On the other hand, by (\ref{def-second-renormalization}) and (\ref{def-calX}) we have
\begin{align}
\begin{split}
    \mathcal{E} & [H W]
    = \mathcal{E} [H G \underline{\Sigma}_M G_{\star}] \\
    & = H G \underline{\Sigma}_M G_{\star} - \frac{\sqrt{z}}{N} \sum_{i\mu} \dimu (\partial_{i\mu} G) \underline{\Sigma}_M G_{\star}
    + \frac{\sqrt{z z_{\star}}}{N} \sum_{i\mu} \dimu G \underline{\Sigma}_M G_{\star} \dimu G_{\star} \\
    & = \mathcal{E} [H G] \underline{\Sigma}_M G_{\star} + \sqrt{z z_{\star}} \mathcal{S} [W] G_{\star} \\
    & = \mathcal{E} [H G] \underline{\Sigma}_M G_{\star} + \sqrt{z z_{\star}} \mathcal{S}_{\mathrm{d}} [W] G_{\star} + \sqrt{z z_{\star}} \mathcal{S}_{\mathrm{o}} [W] G_{\star}. 
\end{split} \label{eqn-expansion-W}
\end{align}
Multiplying both sides of (\ref{eqn-expansion-W}) with $\Pi$ and combining it with (\ref{eqn-expansion-G}), we obtain
\begin{equation}
    W = \Pi \underline{\Sigma}_M \Pi_{\star} 
    + \sqrt{z z_{\star}} \Pi \mathcal{S}_{\mathrm{d}} [W] \Pi_{\star}
    - \Pi \mathcal{E} [H W]
    + \mathcal{Y},
    \label{eqn-intermediate-estimate-W}
\end{equation}
where $\mathcal{Y} = \sum_{\ell = 1}^4 \mathcal{Y}_\ell$ encompasses the error terms that are relatively easy to manage,
\begin{align*}
    & \mathcal{Y}_1 := \Pi \underline{\Sigma}_M ( G_{\star} - \Pi_{\star} ),
    & \quad & \mathcal{Y}_2 := z \Pi \mathcal{S} [G - \Pi] W,  \\
    & \mathcal{Y}_3 := \sqrt{z z_{\star}} \Pi \mathcal{S}_{\mathrm{d}} [W] ( G_{\star} - \Pi_{\star} ),
    & \quad & \mathcal{Y}_4 := \sqrt{z z_{\star}} \Pi \mathcal{S}_{\mathrm{o}} [W] G_{\star}.
\end{align*}
Actually, by leveraging the local laws for a single resolvent and employing the trivial bounds $\| H \|, \| G \| \prec 1$, it is not difficult to verify that for all $1 \leq \ell \leq 4$,
\begin{equation*}
    \langle {\mathcal{Y}_\ell \underline{\Sigma}_M} \rangle,
    \langle {\mathcal{Y}_\ell \underline{I}_N} \rangle = \oprec (N^{-1})
    \qand
    \bfk{u}^\top \mathcal{Y}_\ell \bfk{v} = \oprec (N^{-1/2}),
\end{equation*}
where the isotropic estimate holds uniformly in deterministic $\bfk{u}, \bfk{v} \in \mathbb{R}^{M + N}$ with bounded norms. Now, multiplying both sides of (\ref{eqn-intermediate-estimate-W}) with $\underline{\Sigma}_M$ or $\underline{I}_N$ and taking trace, we obtain
\begin{align*}
    \langle W_{\alpha} \Sigma \rangle & = \langle \Pi_M \Sigma \Pi_{M \star} \Sigma \rangle 
    + \sqrt{z z_{\star}} \langle W_{\beta} \rangle \langle \Pi_M \Sigma \Pi_{M \star} \Sigma \rangle 
    - \langle \Pi \mathcal{E} [H W] \underline{\Sigma}_M \rangle + \oprec (N^{-1}), \\
    \langle W_{\beta} \rangle & =
    \sqrt{z z_{\star}} m m_{\star} \langle W_{\alpha} \Sigma \rangle 
    - \langle \Pi \mathcal{E} [H W] \underline{I}_N \rangle + \oprec (N^{-1}).
\end{align*}
Combining the two equations, we obtain
\begin{equation}
    \langle W_{\alpha} \Sigma \rangle
    - \frac{\langle \Pi_M \Sigma \Pi_{M \star} \Sigma \rangle}{1 - z z_{\star} m m_{\star} \langle \Pi_M \Sigma \Pi_{M \star} \Sigma \rangle} 
    = \langle \mathcal{E} [H W] B \rangle + \oprec (N^{-1}),
    \label{eqn-intermediate-estimate-W-alpha}
\end{equation}
where the deterministic matrix $B$ is given by
\begin{equation*}
    B = - \frac{1}{1 - z z_{\star} m m_{\star} \langle \Pi_M \Sigma \Pi_{M \star} \Sigma \rangle}  \begin{bmatrix} \Sigma \Pi_M & 0 \\ 0 & \sqrt{z z_{\star}} \langle \Pi_M \Sigma \Pi_{M \star} \Sigma \rangle \Pi_N \end{bmatrix}.
\end{equation*}
Using the self-consistent equation (\ref{eqn-self-consistent}), it is not difficult to obtain
\begin{equation*}
    \langle \Pi_M \Sigma \Pi_{M \star} \Sigma \rangle
    = \frac{1}{z z_{\star}} \left ( \frac{1}{m m_{\star}} - \frac{1}{m[z, z_{\star}]} \right ),
\end{equation*}
On the one hand, this implies $\| B \| \lesssim 1$. On the other hand, upon inspecting the definition of $\Xi_\alpha$ in (\ref{def-Xi}), we can rewrite (\ref{eqn-intermediate-estimate-W-alpha}) as
\begin{equation}
    \langle {(W - \Xi) \underline{\Sigma}_M} \rangle
    = \langle {(W_{\alpha} - \Xi_{\alpha}) \Sigma} \rangle
    = \langle \mathcal{E} [H W] B \rangle + \oprec (N^{-1}).
    \label{eqn-representation-average}
\end{equation}
It is worth noting that a similar equation can be established for $\langle {(W - \Xi) \underline{I}_N} \rangle = \langle {W_{\beta} - \Xi_{\beta}} \rangle$. We now claim that (\ref{eqn-representation-average}) can be used to obtain the following averaged laws:
\begin{equation}
    \langle (W - \Xi) \underline{\Sigma}_M \rangle = \oprec (N^{-1})
    \qand
    \langle (W - \Xi) \underline{I}_N \rangle = \oprec (N^{-1}).
    \label{eqn-two-resolvent-average}
\end{equation}
Assuming that (\ref{eqn-two-resolvent-average}) has already been established, we can apply it back to the r.h.s. of (\ref{eqn-intermediate-estimate-W}). As a result, for any $\bfk{u}, \bfk{v} \in \mathbb{R}^{M + N}$ with bounded norms, we obtain
\begin{equation}
    \bfk{u}^\top (W - \Xi) \bfk{v}
    = \bfk{u}^\top \Pi \mathcal{E} [H W] \bfk{v} + \oprec (N^{-1/2}),
    \label{eqn-representation-isotropic}
\end{equation}
where we utilized the fact that the definition of $\Xi$ satisfies
\begin{align*}
    \Xi_{\alpha} & = \Pi_M \Sigma \Pi_{M \star} + \sqrt{z z_{\star}} \langle \Xi_{\beta} \rangle \Pi_M \Sigma \Pi_{M \star}, \\
    \Xi_{\beta} & = \sqrt{z z_{\star}} m m_{\star} \langle \Xi_{\alpha} \Sigma \rangle I_N.
\end{align*}
Now with (\ref{eqn-representation-isotropic}), one can proceed as before to prove the isotropic estimate
\begin{equation}
    \bfk{u}^\top (W - \Xi) \bfk{v}
    =  \oprec (N^{-1/2}),
    \label{eqn-two-resolvent-isotropic}
\end{equation}
which concludes the proof of Proposition \ref{prop-double-resolvent}. It is essential to emphasize that the proof of (\ref{eqn-two-resolvent-isotropic}) using (\ref{eqn-representation-isotropic}) is completely analogous to the proof of (\ref{eqn-two-resolvent-average}) using (\ref{eqn-representation-average}). For the sake of simplicity, we only provided the proof of (\ref{eqn-two-resolvent-average}) to illustrate the idea.

\begin{proof}[Proof of (\ref{eqn-two-resolvent-average})]
Let us concentrate on the left estimate of (\ref{eqn-two-resolvent-average}). It suffices to prove that, for each fixed $n \in \mathbb{N}$, we have
\begin{equation*}
    \mathbb{E} \big | {\langle (W - \Xi) \underline{\Sigma}_M \rangle} \big |^{2 n} \prec N^{-2 n}.
\end{equation*}
Due to the technical reasons mentioned in Section \ref{sec-recurrence-equation}, we again make the assumption that $| {\operatorname{Im} z} | \wedge | {\operatorname{Im} z_{\star}} | \geq N^{-C}$ for a sufficiently large constant $C > 0$. This assumption allows us to freely utilize the deterministic bound $\| G \| \leq N^C$. As before, one can employ a simple continuity argument to extend the result to all $z, z_{\star} \in \mathbb{D}$. Details are omitted here. Using (\ref{eqn-representation-average}), we can write
\begin{align}
\begin{split}
    \mathbb{E} \big | {\langle (W - \Xi) \underline{\Sigma}_M \rangle} \big |^{2 n}
    = & \ \oprec (N^{-2n}) 
    + \oprec (N^{-1}) \mathbb{E} \big | {\langle (W - \Xi) \underline{\Sigma}_M \rangle} \big |^{2n-1} \\
    & + \mathbb{E} \big [ {\langle {\mathcal{E} [H W] B} \rangle 
    \langle {(W - \Xi) \underline{\Sigma}_M} \rangle^{n-1} \langle {(W^* - E^*) \underline{\Sigma}_M} \rangle^n} \big ].    
\end{split} \label{eqn-high-moment-1}
\end{align}
For the second line of (\ref{eqn-high-moment-1}), we recall that
\begin{equation*}
    \langle {\mathcal{E} [H W] B} \rangle =
    \sum_{i \mu} \Big [ {x_{i \mu}  \langle {\dimu W B} \rangle 
    - \frac{1}{N} \langle {\dimu (\partial_{i \mu} W) B} \rangle} \Big ] .
\end{equation*}
Therefore, by applying Lemma \ref{lemma-cumulant-expansion} and the Leibniz rule, we can deduce that, up to some large constant independent of $N$, the second line of (\ref{eqn-high-moment-1}) is bounded by
\begin{align} 
\begin{split}
    & \ N^{-1} \mathbb{E} \bigg [ 
    {\Big | \sum_{i \mu} \langle {\dimu W B} \rangle 
    \cdot \partial_{i \mu} \langle {W \underline{\Sigma}_M} \rangle \Big | 
    \cdot \big | \langle {(W - \Xi) \underline{\Sigma}_M} \rangle \big |^{2n-2}} \bigg ] \\
    + & \ N^{-1} \mathbb{E} \bigg [ 
    {\Big | \sum_{i \mu} \langle {\dimu W B} \rangle 
    \cdot \partial_{i \mu} \langle {W^* \underline{\Sigma}_M} \rangle \Big | 
    \cdot \big | \langle {(W - \Xi) \underline{\Sigma}_M} \rangle \big |^{2n-2}} \bigg ] \\
    + & \ \sum_{p+ \mathfrak{d} (\mathcal{I}) + \mathfrak{d} (\mathcal{I}^*) = 2}^{100n} N^{- (p + \mathfrak{d} (\mathcal{I}) + \mathfrak{d} (\mathcal{I}^*) + 1) /2}   
    \mathbb{E} \left [ \mathcal{D}_{p, \mathcal{I}, \mathcal{I}^*} 
    \big | \langle {(W - \Xi) \underline{\Sigma}_M} \rangle \big |^{2n- 1 - |\mathcal{I}| - |\mathcal{I}^*|} \right ] \\
    + & \ \oprec (N^{-2n})
\end{split} \label{eqn-high-moment-2}
\end{align}
where the summation in the third line of (\ref{eqn-high-moment-2}) is taken over multisets $\mathcal{I}, \mathcal{I}^* \subset \mathbb{N}_+$ (allowed to be empty). Additionally, we use the notation $\mathfrak{d} (\mathcal{I}) := \sum_{q \in \mathcal{I}} q$ and define
\begin{equation*}
    \mathcal{D}_{p, \mathcal{I}, \mathcal{I}^*}
    := \sum_{i \mu} \big | {\partial_{i \mu}^p \langle {\dimu W B} \rangle} \big | 
    \prod_{q \in \mathcal{I}} \big | {\partial_{i \mu}^q \langle {W \underline{\Sigma}_M} \rangle} \big | 
    \prod_{r \in \mathcal{I}^*} \big | {\partial_{i \mu}^{r} \langle {W^* \underline{\Sigma}_M} \rangle} \big |.
\end{equation*}
We note that the upper limit $100n$ in the summation is casual, as long as it depends only on $n$ and ensures that the reminder terms arising from Lemma \ref{lemma-cumulant-expansion} can be controlled by $\oprec (N^{-2n})$. The control of the remainder terms is analogous to the one presented in Section \ref{subsec-reminder}, and therefore, it is omitted here.

Now let us examine the first line of (\ref{eqn-high-moment-2}). The exact definition of $B$ is irrelevant here, but we recall that $\| B \| \lesssim 1$. Also recall the trivial norm bound $\| G \| \prec 1$. As a result, we have
\begin{align*}
    & \ \Big | {\sum_{i \mu} \langle {\dimu W B} \rangle 
    \cdot \partial_{i \mu} \langle {W \underline{\Sigma}_M} \rangle} \Big | \\
    \lesssim & \ \Big | {\sum_{i \mu} \langle {\dimu W B} \rangle 
    \langle {\dimu W \underline{\Sigma}_M G} \rangle} \Big | 
    + \Big | {\sum_{i \mu} \langle {\dimu W B} \rangle 
    \langle {\dimu  W^\top \underline{\Sigma}_M G_{\star}} \rangle} \Big | \\
    = & \ N^{-1} \big | {\big \langle {\underline{\Sigma}_M (W B + B^\top W^\top) \underline{I}_N (W \underline{\Sigma}_M G + G \underline{\Sigma}_M W^\top)} \big \rangle} \big |  \\
    & + N^{-1} \big | {\big \langle {\underline{\Sigma}_M (W B + B^\top W^\top) \underline{I}_N (W^\top \underline{\Sigma}_M G_{\star} + G_{\star} \underline{\Sigma}_M W)} \big \rangle} \big |  
    = \oprec (N^{-1}).
\end{align*}
The same argument also applies to the second line of (\ref{eqn-high-moment-2}). Therefore, we can conclude that the first two lines of (\ref{eqn-high-moment-2}) can be controlled by
\begin{equation*}
     \oprec (N^{-2}) \mathbb{E} \big | {\langle (W - \Xi) \underline{\Sigma}_M \rangle} \big |^{2 n - 2}
     + \oprec (N^{-2n}).
\end{equation*}
Next, let us turn to the third line of (\ref{eqn-high-moment-2}). By utilizing the trivial norm bounds, it is easy to verify that for any fixed $p \geq 0$ and $q, r \geq 1$, we have
\begin{equation*}
    \big | {\partial_{i \mu}^p \langle {\dimu W B} \rangle} \big |,
    \big | {\partial_{i \mu}^q \langle {W \underline{\Sigma}_M} \rangle} \big |, 
    \big | {\partial_{i \mu}^{r} \langle {W^* \underline{\Sigma}_M} \rangle} \big |
    \prec N^{-1}.
\end{equation*}
This readily yields the naive bound 
\begin{equation*}
    \mathcal{D}_{p, \mathcal{I}, \mathcal{I}^*} \prec N^2 \cdot N^{-1 - |\mathcal{I}| - |\mathcal{I}^*|}.
\end{equation*}
With this naive bound, we deduce that the summands in the third line of (\ref{eqn-high-moment-2}) can be stochastically dominated by
\begin{equation*}
    N^{-2p} + N^{- (p + \mathfrak{d} (\mathcal{I}) + \mathfrak{d} (\mathcal{I}^*) - 3) /2}
    N^{-|\mathcal{I}| - |\mathcal{I}^*| - 1} 
    \mathbb{E} \big | \langle {(W - \Xi) \underline{\Sigma}_M} \rangle \big |^{2n- 1 - |\mathcal{I}| - |\mathcal{I}^*|}.
\end{equation*}
which is already sufficient when $p+ \mathfrak{d} (\mathcal{I}) + \mathfrak{d} (\mathcal{I}^*) \geq 3$. For the summands with $p+ \mathfrak{d} (\mathcal{I}) + \mathfrak{d} (\mathcal{I}^*) = 2$, we need to handle the summation more carefully. Let us remark that in this case, at least one of the following holds: (i) $p = 0$ or $2$; (ii) one of $\mathcal{I}$ or $\mathcal{I}^*$ equals $\{ 1 \}$.

For any $(M+N) \times (M+N)$ matrix $A$, we denote $\underline{A} := \uSig^{1/2} A \uSig^{1/2}$. Note that by the Cauchy–Schwarz inequality, if $\| A \| \prec 1$, then
\begin{equation}
    \frac{1}{N} \sum_{i \mu} | \underline{A}_{i \mu} | 
    \lesssim \Big ( {\sum_{i \mu} \big | \underline{A}_{i \mu} \big |^2} \Big )^{1/2}
    = \big ( {N \langle {\underline{\Sigma}_M A \underline{I}_N A^*} \rangle} \big )^{1/2}
    \prec N^{1/2}.
    \label{eqn-cauchy-bound}
\end{equation}
Using the derivative formula $\partial_{i\mu} G = -\sqrt{z} G \dimu G$ as well as the Leibniz rule, it is evident that the derivative $\partial_{i \mu}^p \langle {\dimu W B} \rangle$ can be expressed as a sum of products of the $ii$, $\mu \mu$, $i \mu$, or $\mu i$-th entries of certain matrices with underline. For example,
\begin{align*}
    \langle {\dimu W B} \rangle
    = & \ N^{-1} \underline{WB}_{i \mu} + \underline{WB}_{\mu i} \\
    \partial_{i \mu} \langle {\dimu W B} \rangle
    = & -\sqrt{z} \langle {\dimu G \dimu W B} \rangle
    - \sqrt{z_{\star}} \langle {\dimu W \dimu G_{\star} B} \rangle \\
    = & - \frac{\sqrt{z}}{N} \big ( {\underline{G}_{i i} \underline{WB}_{\mu \mu} + 
    + \underline{G}_{\mu \mu} \underline{WB}_{i i}
    + \underline{G}_{i \mu} \underline{WB}_{\mu i}
    + \underline{G}_{\mu i} \underline{WB}_{i \mu}} \big ) \\
    & - \frac{\sqrt{z_{\star}}}{N} \big ( {\underline{W}_{i i} \underline{G_{\star} B}_{\mu \mu} + 
    + \underline{W}_{\mu \mu} \underline{G_{\star} B}_{i i}
    + \underline{W}_{i \mu} \underline{G_{\star} B}_{\mu i}
    + \underline{W}_{\mu i} \underline{G_{\star} B}_{i \mu}} \big ).
\end{align*}
Importantly, for each product in the sum, both the number of $i$ indices and the number of $\mu$ indices should be equal to $(p+1)$. Consequently, when $p$ is even, each product must contain at least one entry indexed by $i \mu$. In this case, we can bound the remaining entries trivially by $\oprec (1)$ and activate (\ref{eqn-cauchy-bound}) to obtain
\begin{equation*}
    \sum_{i \mu} \big | {\partial_{i \mu}^p \langle {\dimu W B} \rangle} \big | = \oprec (N^{1/2}),
    \quad \text{ if } p \in 2 \mathbb{N}.
\end{equation*}
This strategy also applies to $\partial_{i \mu}^q \langle {W \underline{\Sigma}_M} \rangle$ (resp. $\partial_{i \mu}^r \langle {W \underline{\Sigma}_M} \rangle$) when $q$ (resp. $r$) is odd. In particular, we have
\begin{equation*}
    \sum_{i \mu} \big | {\partial_{i \mu} \langle {W \underline{\Sigma}_M} \rangle} \big |,
    \sum_{i \mu} \big | {\partial_{i \mu} \langle {W^* \underline{\Sigma}_M} \rangle} \big |
    = \oprec (N^{1/2}).
\end{equation*}
Consequently, when $p+ \mathfrak{d} (\mathcal{I}) + \mathfrak{d} (\mathcal{I}^*) = 2$, we always have
\begin{equation*}
    \mathcal{D}_{p, \mathcal{I}, \mathcal{I}^*} \prec N^{1/2} \cdot N^{-|\mathcal{I}| - |\mathcal{I}^*|}.
\end{equation*}
With this improved bound, all the summands in the third line of (\ref{eqn-high-moment-2}) can be controlled by
\begin{equation*}
    \oprec (N^{-2p}) + 
    \oprec (N^{-|\mathcal{I}| - |\mathcal{I}^*| - 1})
    \mathbb{E} \big | \langle {(W - \Xi) \underline{\Sigma}_M} \rangle \big |^{2n- 1 - |\mathcal{I}| - |\mathcal{I}^*|}.
\end{equation*}
Now, summarising the estimates for (\ref{eqn-high-moment-1}) and (\ref{eqn-high-moment-2}) above, we deduce that 
\begin{equation*}
    \mathbb{E} \big | {\langle (W - \Xi) \underline{\Sigma}_M \rangle} \big |^{2 n}
    = \oprec (N^{-2n}) 
    + \sum_{\ell = 1}^{2n} \oprec (N^{-\ell}) \mathbb{E} \big | {\langle (W - \Xi) \underline{\Sigma}_M \rangle} \big |^{2n-\ell}.
\end{equation*}
This concludes our proof upon inspecting the definition of stochastic domination.
\end{proof}

\subsection{Special case: two identical spectral parameters} 
\label{subsec-double-resol-same}

In this subsection we present an alternative proof of Proposition \ref{prop-double-resolvent} when $z = z_{\star}$. It suffices to show that, for deterministic $\bfk{u}, \bfk{v} \in \mathbb{R}^{M + N}$ with $\| \bfk{u} \|, \| \bfk{v} \| \lesssim 1$,
\begin{equation*}
    \bfk{u}^\top W \bfk{v}
    = \bfk{u}^\top \begin{bmatrix}
    ({m^\prime}/{m^2}) \Pi_M \Sigma \Pi_M & 0 \\
    0 & z^{-1} ( {m^\prime}/{m^2} - 1 )
    \end{bmatrix} \bfk{v} + \oprec (N^{-1/2}).
\end{equation*}
For $t \geq 0$ and $z \in \mathbb{D}_0$, consider the parameterized resolvents
\begin{equation*}
    \hat{G} (t) 
    \equiv \hat{G} (t, z) 
    := \left ( H - z - \begin{bmatrix} t \Sigma & 0 \\ 0 & 0 \end{bmatrix} \right )^{-1},
\end{equation*}
so that $W = \partial_{t} \hat{G} (0)$. By introducing $\hat{\Sigma} (t) \equiv \hat{\Sigma}(t, z) = {z\Sigma}/(z + t \Sigma)$, we can write
\begin{equation*}
    \hat{G} (t) 
    = \begin{bmatrix}
    z / (z + t \Sigma) & 0 \\ 0 & I
    \end{bmatrix}^{1/2}
    \begin{bmatrix}
    -z & \sqrt{z} \hat{\Sigma}(t)^{1/2} X \\
    \sqrt{z} X^\top \hat{\Sigma}(t)^{1/2} & -z
    \end{bmatrix}^{-1}
    \begin{bmatrix}
    z / (z + t \Sigma) & 0 \\ 0 & I
    \end{bmatrix}^{1/2}.
\end{equation*}
In particular, the middle matrix is the Green function of the (linearized) sample covariance matrix $\hat{\Sigma}(t)^{1/2} X$, whose deterministic approximation can be deduced from Proposition \ref{prop-local-law} as long as $\hat{\Sigma}(t)$ fulfills Assumptions \ref{assumption-spectrum} and \ref{assumption-regularity}. Actually, by $\| \Sigma \| \lesssim 1$ and $z \asymp 1$ for $z \in \mathbb{D}_0$, it is not hard to show that, for some sufficiently small constant $c > 0$, Assumptions \ref{assumption-spectrum} and \ref{assumption-regularity} are satisfied by all matrices in $\{ \hat{\Sigma}(t, z) : t \in [0, c], z \in \mathbb{D}_0 \}$ after some appropriate renaming of $\tau$.
Therefore, we may define $\hat{m} (t) \equiv \hat{m} (t, z)$ as the unique solution in $\mathbb{C}_+$ to the equation
\begin{equation}
    \frac{1}{\hat{m}(t)} 
    = -z + \frac{1}{N} \operatorname{tr} \frac{\hat{\Sigma}(t)}{1 + \hat{m}(t) \hat{\Sigma}(t)} 
    = -z + \frac{1}{N} \operatorname{tr} \frac{z \Sigma}{z + t \Sigma + z \hat{m}(t) \Sigma},
    \label{eqn-self-consistent-hat}
\end{equation}
so that by Proposition \ref{prop-local-law}, we have
\begin{align}
    \bfk{u}^\top \hat{G} (t) \bfk{v}
    & = \bfk{u}^\top \begin{bmatrix}
    z / (z + t \Sigma) & 0 \\ 0 & I
    \end{bmatrix}^{1/2}
    \begin{bmatrix} 
    \hat{\Pi}_M (t) & 0 \\
    0 & \hat{\Pi}_N (t)
    \end{bmatrix}
    \begin{bmatrix}
    z / (z + t \Sigma) & 0 \\ 0 & I
    \end{bmatrix}^{1/2} \bfk{v} + \oprec (N^{-1/2}) \nonumber \\
    & = \bfk{u}^\top \begin{bmatrix}
    - (z + t \Sigma + z \hat{m} (t) \Sigma)^{-1} & 0 \\
    0 & \hat{m} (t)
    \end{bmatrix} \bfk{v} + \oprec (N^{-1/2}),
    \label{eqn-local-law-Ghat}
\end{align}
where we used the notations
\begin{equation*}
    \hat{\Pi}_M (t)
    := - \frac{1}{z [ 1 + \hat{m} (t) \hat{\Sigma} (t) ]}
    = - \frac{(z + t \Sigma)/z}{z + t \Sigma + z \hat{m} (t) \Sigma}
    \quad \text{ and } \quad
    \hat{\Pi}_N (t) = \hat{m} (t) I.
\end{equation*}
By Cauchy’s integral formula for derivatives, the estimate (\ref{eqn-local-law-Ghat}) remains valid after taking derivative on both sides w.r.t. the parameter $t$. Now let us derive the expression for $\partial_t \hat{m} (0)$. Recall $m^\prime \equiv \partial_z m (z)$. By taking derivative on both sides of the self-consistent equation (\ref{eqn-self-consistent}) w.r.t. $z$, we can obtain
\begin{equation*}
    m^\prime
    = \frac{m^2}{1 - N^{-1} z^2 m^2 \operatorname{tr} \Pi_M \Sigma \Pi_M \Sigma}.
\end{equation*}
On the other hand, we have obviously $\hat{m} (0) = m$. By taking derivative on both sides of (\ref{eqn-self-consistent-hat}) w.r.t. the parameter $t$ and setting $t = 0$, we can get
\begin{equation*}
    \partial_t \hat{m}(0) = \frac{N^{-1} z m^2 \operatorname{tr} \Pi_M \Sigma \Pi_M \Sigma}{1 - N^{-1} z^2 m^2 \operatorname{tr} \Pi_M \Sigma \Pi_M \Sigma} 
    = \frac{1}{z} \left ( \frac{m^\prime}{m^2} - 1 \right ).
\end{equation*}
It follows that
\begin{equation*}
    \bfk{u}^\top \partial_{t} \hat{G} (0) \bfk{v}
    = \bfk{u}^\top \begin{bmatrix}
    ({m^\prime}/{m^2}) \Pi_M \Sigma \Pi_M & 0 \\
    0 & z^{-1} ( {m^\prime}/{m^2} - 1 )
    \end{bmatrix} \bfk{v} + \oprec (N^{-1/2}),
\end{equation*}
where we used
\begin{equation*}
    \partial_t \left [ - \frac{1}{z + t \Sigma + z \hat{m} (t) \Sigma} \right ]_{t = 0}
    = \left [ \frac{(1 + z \partial_t \hat{m}(t)) \Sigma}{(z + t \Sigma + z \hat{m} (t) \Sigma)^2} \right ]_{t = 0}
    = \frac{m^\prime}{m^2} \Pi_M \Sigma \Pi_M.
\end{equation*}
We therefore conclude the proof by recalling $\bfk{u}^\top W \bfk{v} = \bfk{u}^\top \partial_{t} \hat{G} (0) \bfk{v}$.

\end{document}